\documentclass[12pt]{article}
\usepackage{amssymb,latexsym}
\title{Non-\sa{} perturbations of
\sa{} \op{}s in 2 dimensions I.}
\author{Michael Hitrik\\MSRI \\1000 Centennial Drive \\ Berkeley, 
CA 94720\\ USA\\hitrik@msri.org \and 
Johannes Sj\"ostrand\\Centre de Math\'ematiques\\Ecole Polytechnique\\FR 
91120 Palaiseau\\France\\johannes@math.polytechnique.fr}
\date{}

\def\wrtext#1{\relax\ifmmode{\leavevmode\hbox{#1}}\else{#1}\fi}
\def\abs#1{\left|#1\right|}
\def\begeq{\begin{equation}}
\def\endeq{\end{equation}}
\def\Remark{\vskip 2mm \noindent {\em Remark}}

\def\ekv#1#2{\begeq\label{#1}#2\endeq}
\def\eekv#1#2#3{\begin{eqnarray}\label{1}#2\\ #3\nonumber\end{eqnarray}}

\def\iint{\int\hskip -2mm\int}

\def\an{analytic}
\def\asy{asymptotic}
\def\bdd{bounded}
\def\bdy{boundary}
\def\coef{coefficient}

\def\canform{canonical transformation}

\def\ctf{canonical transformation}
\def\diffeo{diffeomorphism}

\def\dop{differential operator}

\def\ev{eigenvalue}
\def\e{equation}
\def\fu{function}
\def\fy{family}
\def\F{Fourier}
\def\fop{Fourier integral operator}
\def\fourior{Fourier integral operator}

\def\hol{holomorphic}
\def\indep{independent}
\def\lhs{left hand side}
\def\mfld{manifold}
\def\ml{microlocal}
\def\neigh{neighborhood}
\def\nondeg{non-degenerate}
\def\op{operator}

\def\pb{problem}

\def\plsh{plurisubharmonic}

\def\pro{proposition}
\def\Pro{Proposition}
\def\pop{pseudodifferential operator}

\def\pseudor{pseudodifferential operator}
\def\res{resonance}
\def\rhs{right hand side}
\def\sa{selfadjoint}

\def\sop{Schr\"odinger operator}
\def\st{strictly}
\def\stpsh{\st{} plurisubharmonic}

\def\sufly{sufficiently}
\def\tf{transformation}
\def\Th{Theorem}
\def\th{theorem}
\def\tf{transform}

\def\ufly{uniformly}

\def\wrt{with respect to}

\def\Re{{\rm Re\,}}
\def\Im{{\rm Im\,}}

\newcommand{\eps}{\epsilon}
\def\part#1{\frac{\partial}{\partial #1}}
\def\half{\frac{1}{2}}

\def\norm#1{||\,#1\,||}

\newcommand{\real}{\mbox{\bf R}}
\newcommand{\comp}{\mbox{\bf C}}
\newcommand{\z}{\mbox{\bf Z}}
\newcommand{\nat}{\mbox{\bf N}}

\newcommand{\dist}{\mbox{\rm dist\,}}
\newcommand{\Spec}{\mbox{\rm Spec\,}}

\renewcommand{\Re}{\mbox{\rm Re\,}}
\renewcommand{\Im}{\mbox{\rm Im\,}}

\renewcommand{\exp}{\mbox{\rm exp\,}}
\newcommand{\supp}{\mbox{\rm supp}}

\newtheorem{dref}{Definition}[section]
\newtheorem{lemma}[dref]{Lemma}
\newtheorem{theo}[dref]{Theorem}
\newtheorem{prop}[dref]{Proposition}

\newenvironment{proof}{\vspace{.3cm}\noindent{{\em Proof:}}}{\hfill$\Box$}
\begin{document}
\maketitle

\begin{abstract}
This is the first in a series of works devoted to small non-\sa{} perturbations 
of selfadjoint $h$-\pop{}s in dimension 2. In the present work we treat 
the case when the classical flow of the unperturbed part is periodic and 
the strength $\epsilon $ of the perturbation is $\gg h$ (or sometimes only $\gg h^2$) 
and bounded from above by $h^\delta $ for some $\delta >0$. We get a 
complete \asy{} description of all \ev{}s in certain rectangles 
$[-1/C,1/C]+i\epsilon [F_0-1/C,F_0+1/C]$.  
\end{abstract}



\section{Introduction}\label{section0}
\setcounter{equation}{0}

\par In \cite{MeSj}, A. Melin and the second author observed that for a wide
and stable class of non-\sa{} \op{}s in dimension 2 and in the
semi-classical limit ($h\to 0$), it is possible to describe all \ev{}s
individually in an $h$-\indep{} domain in ${\bf C}$, by means of a
Bohr-Sommerfeld quantization condition. This result is quite remarkable
since the corresponding conclusion in the \sa{} case seems to be
possible only in dimension 1 or under strong (and unstable) assumptions
of complete integrability. The underlying reason for this result is the
absence of small denominators which allows us to avoid the usual
trouble with exceptional sets in the KAM \th{}. 

\par As a next step, the second author noticed (\cite{Sj}) that for
non-\sa{} \op{}s of the form $P(x,hD_x)+i\epsilon Q(x,hD_x)$ it is
possible to find a similar result, when $P$ is selfadjoint, $\epsilon
>0$ small and fixed and the classical bicharacteristic flow is periodic
on each real energy surface. (Again, it is important that we are in
dimension 2.) The method is similar to the one in \cite{MeSj} and uses
non-linear Cauchy-Riemann \e{}s, now in an "$\epsilon $-degenerate"
form. (See also \cite{Sj2} for a different extension.)

\par It soon became quite clear that we run into a fairly vast program,
and that logically one should start with even smaller perturbations,
say $\epsilon ={\cal O}(h^\delta )$, for some $\delta >0$. The present
work is planned to be the first in a series, devoted to small
perturbations of \sa{} \op{}s in dimension 2. In addition to the
challenge of doing plenty of things in dimension 2, that can usually
only be done in dimension 1, we have been motivated by recent progress
around the damped wave equation (\cite{Le}, \cite{AsLe}, \cite{Sj3}, \cite{Hi}), as well as
the problem of barrier top \res{}s for the semi-classical \sop{}
(\cite{KaKe})  where more complete results than the corresponding ones for
\ev{}s of potential wells (\cite{Sj4}, \cite{BaGrPa}, \cite{Po}) seem possible. One
long term goal of this series is to get improved results on the
distribution of
\res{}s for
\st{} convex obstacles in ${\bf R}^3$. See \cite{SjZw2} (and references given
there) for a first result on Weyl \asy{}s for the real parts inside
certain bands. In the case of analytic obstacles, much more can probably
be said, especially in dimension 3 (and 2).

\par Let $M$ denote ${\bf R}^2$ or a compact real-analytic \mfld{} of
dimension 2.

\par When $M={\bf R}^2$, let 
\begeq
\label{0.1}
P_\epsilon =P(x,hD_x,\epsilon ;h)
\endeq
be the Weyl quantization on ${\bf R}^2$ of a symbol $P(x,\xi ,\epsilon
;h)$ depending smoothly on $\epsilon \in{\rm neigh\,}(0,{\bf R})$ with
values in the space of \hol{} \fu{}s of $(x,\xi )$ in a tubular
\neigh{} of ${\bf R}^4$ in ${\bf C}^4$, with 
\begeq\label{0.2}
\vert P(x,\xi ,\epsilon ;h)\vert \le Cm(\Re (x,\xi ))
\endeq 
there. Here $m$ is assumed to be an order \fu{} on ${\bf R}^4$, in the
sense that $m>0$ and 
\begeq\label{0.3}
m(X)\le C_0\langle X-Y\rangle ^{N_0}m(Y),\ X,Y\in{\bf R}^4.
\endeq
We also assume that 
\begeq\label{0.4}
m\ge 1.\endeq
We further assume that 
\begeq\label{0.5}
P(x,\xi ,\epsilon ;h)\sim \sum_{j=0}^\infty  p_{j,\epsilon }(x,\xi)h^j,\ h\to 0, \endeq
in the space of such \fu{}s. We make the ellipticity assumption
\begeq\label{0.6}
\vert p_{0,\epsilon }(x,\xi )\vert \ge {1\over C}m(\Re (x,\xi )),\
\vert (x,\xi )\vert \ge C, \endeq
for some $C>0$.

\par When $M$ is a compact \mfld{}, we let 
\begeq\label{0.7}
P_\epsilon =\sum_{\vert \alpha \vert \le m}a_{\alpha ,\epsilon
}(x;h)(hD_x)^\alpha ,\endeq
be a \dop{} on $M$, such that for every choice of local coordinates,
centered at some point of $M$, $a_{\alpha ,\epsilon }(x;h)$ is a smooth
\fu{} of $\epsilon $ with values in the space of \bdd{} \hol{} \fu{}s
in a complex \neigh{} of $x=0$. We further assume that 
\begeq\label{0.7.5}
a_{\alpha ,\epsilon }(x;h)\sim \sum_{j=0}^\infty  a_{\alpha ,\epsilon
,j}(x)h^j,\ h\to 0, \endeq
in the space of such \fu{}s. The semi-classical principal symbol in
this case is given by
\begeq\label{0.8}
p_{0,\epsilon }(x,\xi )=\sum a_{\alpha ,\epsilon ,0}(x)\xi ^\alpha ,
\endeq
and we make the ellipticity assumption
\begeq\label{0.9}
\vert p_{0,\epsilon}(x,\xi )\vert \ge {1\over C}\langle \xi \rangle ^m,\ (x,\xi
)\in T^*M,\,\vert \xi \vert \ge C,\endeq
for some large $C>0$. (Here we assume that $M$
has been equipped with some Riemannian metric, so that $\vert \xi \vert
$ and $\langle \xi \rangle =(1+\vert \xi \vert ^2)^{1/2}$ are
well-defined.)

\par Sometimes, we write $p_\epsilon $ for $p_{0,\epsilon }$ and
simply $p$ fo $p_{0,0}$. Assume 
\begeq\label{0.10}
P_{\epsilon =0} \hbox{ is formally \sa{}.}
\endeq
In the case when $M$ is compact, we let the under\-lying Hil\-bert
spa\-ce be $L^2(M,\mu (dx))$ for some po\-si\-ti\-ve real-ana\-ly\-tic den\-si\-ty 
$\mu(dx)$ on $M$. 

\par Under these assumptions, $P_\epsilon $ will have discrete spectrum
in some fixed \neigh{} of $0\in{\bf C}$, when $h>0,\epsilon \ge 0$
are \sufly{} small, and the spectrum in this region will be contained
in a band $\vert \Im z\vert \le {\cal O}(\epsilon )$. The purpose of
this work and later ones in this series, is to give detailed \asy{}
results about the distribution of individual \ev{}s inside such a band. 

\par Assume for simplicity that (with $p=p_{\epsilon =0}$)
\begeq\label{0.11}
p^{-1}(0)\cap T^*M\hbox{ is connected.}\endeq
Let $H_p=p'_\xi \cdot {\partial \over \partial x}-p'_x\cdot {\partial
\over \partial \xi }$ be the Hamilton field of $p$. In this work, we
will always assume that for $E\in{\rm neigh\,}(0,{\bf R})$:
\begin{eqnarray}\label{0.12}
\hbox{The }H_p\hbox{-flow is periodic on }p^{-1}(E)\cap T^*M\hbox{
with}\\ \hbox{period }T(E)>0 \hbox{ depending \an{}ally on }E.\nonumber
\end{eqnarray}
Let $q={1\over i}{({\partial \over \partial \epsilon })}_{\epsilon
=0}p_\epsilon $, so that 
\begeq\label{0.13}
p_\epsilon =p+i\epsilon q+{\cal O}(\epsilon ^2m),\endeq
in the case $M={\bf R}^2$ and $p_\epsilon =p+i\epsilon q+{\cal
O}(\epsilon ^2\langle \xi \rangle ^m)$ in the \mfld{} case. Let 
\begeq\label{0.14}
\langle q\rangle ={1\over T(E)}\int_{-T(E)/2}^{T(E)/2}q\circ \exp tH_p
dt\hbox{ on }p^{-1}(E)\cap T^*M.\endeq
Notice that $p,\langle q\rangle $
are in involution; $0=H_p\langle q\rangle =:\{ p,\langle q\rangle \} $.
In section \ref{section2}, we shall see how to reduce ourselves to the case when
\begeq\label{0.15}
p_\epsilon =p+i\epsilon \langle q\rangle +{\cal O}(\epsilon ^2),\endeq
near $p^{-1}(0)\cap T^*M$. An easy consequence of this is that the
spectrum of $P_\epsilon $ in $\{z\in {\bf C}; \vert \Re z\vert <\delta
\}$ is confined to $]-\delta ,\delta [+i\epsilon ]\langle \Re q\rangle
_{{\rm min},0}-o(1),\langle \Re q\rangle
_{{\rm max},0}+o(1)[$, when $\delta ,\epsilon ,h\to 0$, where $\langle
\Re q\rangle _{{\rm min},0}=\min_{p^{-1}(0)\cap T^*M}\langle \Re
q\rangle $ and similarly for $\langle q\rangle _{{\rm max},0}$. We will
mainly think about the case when $\langle q\rangle $ is real-valued but
we will work under the more general assumption that 
\begeq\label{0.16}
\Im \langle q\rangle \hbox{ is an \an{} \fu{} of }p\hbox{ and }\Re
\langle q\rangle ,\endeq
in the region of $T^*M$, where $\vert p\vert \le 1/{\cal O}(1)$.

\par Let $F_0\in [\langle \Re q\rangle _{{\rm min},0},\langle \Re
q\rangle _{{\rm max},0}]$. The purpose of the present work is to
determine all \ev{}s in a rectangle
\begeq\label{0.17}
]-{1\over {\cal O}(1)},{1\over {\cal O}(1)}[+i\epsilon ]F_0-{1\over
{\cal O}(1)},F_0+{1\over {\cal O}(1)}[, \endeq
for 
\begeq\label{0.18}
h\ll \epsilon \le {\cal O}(h^\delta ),
\endeq
where $\delta >0$ is any fixed number. (When the subprincipal symbol of
$P$ is zero, we can treat even smaller values of $\epsilon $: $h^2\ll
\epsilon \le {\cal O}(h^\delta )$.) We will achieve this in the following
two cases, under the general assumption that 
\begeq\label{0.19}
T(0)\hbox{ is the minimal period of every }H_p\hbox{-trajectory in
}\Lambda _{0,F_0},\endeq
where
\begeq\label{0.20}
\Lambda _{0,F_0}:=\{ \rho \in T^*M;\, p(\rho )=0,\, \Re \langle
q\rangle (\rho )=F_0\} ,\endeq
in the following two cases:
\smallskip
\par\noindent I) The first case is when 
\begeq\label{0.21}
dp,\, d\Re \langle q\rangle \hbox{ are linearly \indep{} at every
point of }\Lambda _{0,F_0}.\endeq
This implies that every connected component of $\Lambda _{0,F_0}$ is a
two-dimensional Lagrangian torus. For simplicity, we shall assume that
there is only one such component. Notice that in view of (\ref{0.19}), the
space of closed orbits in $p^{-1}(0)\cap T^*M$;
$$\Sigma :=(p^{-1}(0)\cap T^*M)/\sim ,$$
where $\rho \sim \mu $ if $\rho =\exp tH_p \mu $ for some $t\in{\bf
R}$, becomes a 2-dimensional symplectic \mfld{} near the image of
$\Lambda _{0,F_0}$, and (\ref{0.21}) simply means that $\Re \langle q\rangle $, viewed
as a \fu{} on $\Sigma $, has non-vanishing differential along the image of $\Lambda_{0,F_0}$. 
The image of
$\Lambda _{0,F_0}$ is just a closed curve. The main results in this case
are \Th{}s \ref{Th6.2}, \ref{Th6.4} and they show that the \ev{}s form a distorted lattice.
\smallskip

\par\noindent II) The second case is when $F_0\in\{ \langle \Re
q\rangle _{{\rm min},0},\langle \Re q\rangle _{{\rm max},0}\}$. In this
case, we again view $\langle \Re q\rangle $ as a smooth \fu{} on
$\Sigma $ near the image of $\Lambda _{0,F_0}$ and assume that
\begin{eqnarray}\label{0.22}
\hbox{The Hessian of }\langle \Re q\rangle \hbox{ is non-degenerate
(positive }\\ \hbox{or negative) at every point }\rho \in \Sigma ,
\hbox{ with }\langle \Re q\rangle (\rho )=F_0.\nonumber\end{eqnarray}
The main results in this case are given by \Th{}s \ref{Th6.6}, 
\ref{Th6.7} which tell us
that the \ev{}s form a distorted half-lattice.\smallskip

\par\noindent III) The third natural case would be when $F_0$ is a
critical value of $\Re \langle q\rangle $ corresponding to a saddle
point. We hope to study this case in the near future.\smallskip

\par The \an{}ity assumptions are introduced, because the optimal
spaces are deformations of the usual $L^2$-space obtained by adding
exponential weights with exponents that are ${\cal O}(\epsilon )$, and
there are closely related \fop{}s with complex phase some of which have
associated complex \ctf{}s that are $\epsilon $-perturbations of the
identity. When $\epsilon \sim h^\delta$, $0<\delta <1$, 
appropriate Gevrey type assumptions would probably suffice, but in the
case
$\epsilon \sim h$ we seem to need analyticity assumptions at
one point, even though standard
$C^\infty$--\ml{} analysis would suffice for most of the steps. At the opposite
extreme, $\epsilon$ small but
\indep{} of $h$, the analyticity assumptions seem necessary, and in order
to avoid technicalities, we have chosen to assume \an{}ity \indep{}ly of
the size of $\epsilon $.

\par In the \sa{} case there have been many works about operators whose
associated classical flow is periodic (\cite{We}, \cite{Co}, \cite{BoGu}, \cite{HeRo},
\cite{Do}, \cite{Iv}), and we follow one of the main ideas in those works, namely
to use some sort of averaging procedure in order to reduce the
dimension by one unit, so that in our case, we come down to a one-dimensional problem. 
The implementation of this is more complicated in
our case because of the need to work in modified exponentially weighted
spaces (after suitable FBI-\tf{}s). It should also be pointed out that
in the case when $\epsilon $ is small but \indep{} of $h$ (\cite{Sj}), this
does not seem to work and the problem remains two-dimensional. The same
seems to be the case (for the whole scale of $\epsilon $) in other
situations, when the $H_p$-flow is completely integrable without being
periodic, or more generally when the energy surface $p^{-1}(0)\cap T^*M$
contains certain invariant Lagrangian tori. We intend to treat such
situations later in this series.

\par The plan of the paper is the following:
\smallskip
\par\noindent In section \ref{section1}, we reexamine the Egorov \th{} in a form
suitable for us, and complete some observations of \cite{HeSj2} about the
two term version of this result. 
\smallskip
\par\noindent In section \ref{section2} we perform dimension reduction by averaging.
\smallskip
\par\noindent In section \ref{section3} we make a complete reduction in the torus
case (I) and determine the corresponding quasi-\ev{}s.
\smallskip
\par\noindent In section \ref{section4} we do the analogous work in the extreme
case (II).
\smallskip
\par\noindent In section \ref{section5} we justify the earlier computations by
treating an auxiliary global (Grushin) problem, and we obtain the two
main results. 
\smallskip
\par\noindent In section \ref{section6}, we give a first application to barrier top
\res{}s. 
\smallskip
\par\noindent In the appendix, we review some standard facts about
FBI-\tf{}s on \mfld{}s.
\smallskip 
\par The next work(s) in this series (in addition to \cite{Sj}) will remain
in the case when the classical flow of the unperturbed part is
periodic. We intend to study the saddle point case (III), and the case
when
$\langle q\rangle $ vanishes. 

\vskip 4mm
\noindent
{\it Acknowledgements.} We would like to thank Anders Melin and 
Maciej Zworski for useful discussions. The first author gratefully acknowledges 
the support of the Swedish Foundation for International Cooperation 
in Research and Higher Education 
(STINT) as well as of the MSRI postdoctoral fellowship. 

\section{Quantization of \canform{}s between non-simply connected domains in 
phase space}\label{section1}
\setcounter{equation}{0}

\par We first give an affirmative answer to a question asked in
appendix a of \cite{HeSj2}. Let $\kappa :{\rm neigh\,}((y_0,\eta
_0),T^*{\bf R}^n)\to {\rm neigh\,}((x_0,\xi _0),T^*{\bf R}^n)$
be an analytic \canform{} and consider a corresponding \fop{}
\begeq\label{1.1}
Uu(x)=h^{-{n+N\over 2}}\iint e^{i\phi (x,y,\theta
)/h}a(x,y,\theta ;h)u(y)dyd\theta ,\endeq 
with $a=a_0+{\cal O}(h)$, a classical symbol in $S^{0,0}$ (see the
appendix), and
$\phi$ \nondeg{} phase \fu{} in the sense of H\"ormander~\cite{Ho} (without the
homogeneity requirement in $\theta $) which generates the graph
of $\kappa $. (Since we work \ml{}ly, $\phi ,a$ are assumed to be
defined near a fixed point $(x_0,y_0,\theta _0)$ with $\phi
'_\theta (x_0,y_0,\theta _0)=0$, $(x_0,\xi _0)=(x_0,\phi
_x'(x_0,y_0,\theta _0))$, $(y_0,\eta _0)=(y_0,-\phi
_y'(x_0,y_0,\theta _0))$.) We require
$U$ to be unitary:
\begeq\label{1.2}
U^*U=1, \hbox{ microlocally near }(y_0,\eta _0),\endeq
and we are interested in the improved Egorov property:
\begin{eqnarray}\label{1.3}
\hbox{If }PU=UQ,\hbox{ where }P=P^w,\, Q=Q^w\hbox{ are
$h$-pseudodifferential}\\ \hbox{\op{}s of order 0, then
}P\circ
\kappa =Q+{\cal O}(h^2).\nonumber\end{eqnarray}
In appendix a of \cite{HeSj2}, it was shown that such $U$'s exist
and we shall answer the question raised there, by
establishing the following \pro{}. (We learned from C. Fefferman that Jorge
Silva has obtained essentially the same result in the framework
of classical \fop{}s.) 
\begin{prop}\label{Prop1.1}  Within the class of \op{}s
satisfying {\rm (\ref{1.1})} and {\rm (\ref{1.2})}, the pro\-per\-ty {\rm (\ref{1.3})} 
is equivalent to:
\begeq\label{1.4}
{{a_0}_\vert}_{C_\phi }\hbox{ has constant argument.}\endeq
Here $\phi$ is defined in some open set ${\cal D}(\phi )\subset {\bf R}^{2n+N}$ and 
$$
C_\phi =\{ (x,y,\theta )\in{\cal D}(\phi);\, \phi '_\theta (x,y,\theta )=0\}.
$$
\end{prop}

\begin{proof} We first consider the special case of
\pop{}s, i.e. the case when $\kappa $ is the identity. Then
$a_0$ is the principal symbol and (\ref{1.3}) implies that
$\vert a_0\vert =1$ (after inserting an additional factor $(2\pi
)^{-n}$ in front of the integral and taking the standard phase
$\phi =(x-y)\cdot \theta $). Write 
$$U^{-1}PU=P+U^{-1}[P,U].$$
We see that (\ref{1.3}) holds iff $\{ p,a_0\}=0$ for all $p$, i.e. iff
$a_0={\rm Const.}$ The proposition follows in the case of \pop{}s
since we also know in general that the property (\ref{1.4}) is invariant
under changes of $(\phi ,a)$ in the representation of the given
\op{}. 

\par When $\phi $ is quadratic and $a$ is constant, we have a
metaplectic \op{} and $\kappa $ is linear. In that case, we know
that (\ref{1.3}) holds, and using the special case of $h$-\pop{}s, we see
that we have equivalence between (\ref{1.3}) and (\ref{1.4}) in the case when
$\kappa $ is linear. 

\par Consider a smooth deformation of \canform{}s 
$[0,1]\ni t\mapsto \kappa _t$, with a deformation field
$H_{a(t)}$, so that $\partial _t\kappa_t (\rho )=H_{a(t)}(\kappa
_t(\rho ))$ where $a(t)=a(t,x,\xi )$ is smooth and \indep{} of
$h$. Let $A(t)=a^w(x,hD_x)$ and consider a corresponding family
of \fop{}s $U(t)$ associated to $\kappa _t$:
\begeq\label{1.5}
hD_tU(t)+A(t)\circ U(t)=0.\endeq
Since $A(t)$ are \sa{}, unitarity of $U(t)$ is conserved under
the flow of (\ref{1.5}). Let $U(t)$ be such a unitary \fy{}.
\begin{prop}
\label{Prop1.2} 
We have {\rm (\ref{1.3})} for one value of
$t$ iff we have it for all values of $t$.
\end{prop}

\begin{proof} Suppose we have (\ref{1.3}) for $U(0)$. From
(\ref{1.5}) we get 
$$hD_t(U(t)^{-1})=U(t)^{-1}A(t).$$
Consider a family $P(t)=U(t)PU(t)^{-1}$. Then 
$$hD_tP(t)+[A(t),P(t)]=0,$$
and on the level of Weyl symbols, we get 
$$\partial _tP(t)+\{ a(t),P(t)\} ={\cal O}(h^2),$$
or in other words,
$$(\partial _t+H_{a(t)})P(t)={\cal O}(h^2).$$
This means that 
$$P(t)\circ (\kappa _t(\rho ))=P(0)\circ \kappa _0+{\cal
O}(h^2)=P(\rho )+{\cal O}(h^2),$$
where we used (\ref{1.3}) for $U(0)$ in the last step. Then $P(t)$
fulfills (\ref{1.3}) for all $t$.
\end{proof}
\medskip

\par On the other hand, if $U(t)$ fulfills (\ref{1.5}), we know, using
that the subprincipal symbol of $A(t)$ is 0, that if we
represent 
$$U(t)=h^{-{n+N\over 2}}\iint e^{{i\over h}\phi _t(x,y,\theta
)}a_t(x,y,\theta ;h) u(y)dyd\theta ,$$
with $\phi _t$, $a_t$ depending smoothly on $t$, then the
argument of ${{a_{t,0}}_\vert }_{C_{\phi _t}}$ is constant along
every curve in $\{ (t,x,\theta );\, (x,\theta )\in C_{\phi
_t}\}$ corresponding to a $H_{a(t)}$-trajectory: $t\mapsto
(\kappa _t(\rho ),\kappa _0(\rho ))$. This can be seen either by
a direct computation leading to a real transport equation for the
leading symbol, (using that 
$$e^{-i\phi (x)/h}\circ a^w(x,hD_x)\circ e^{i\phi
(x)/h}=(a(x,\phi '(x)+hD_x))^w+{\cal O}(h^2),$$
see appendix a in \cite{HeSj2}), or by using H\"ormander's
definition (\cite{Ho}) of the principal symbol of a \fop{}, as well as a
result of Duistermaat-H\"ormander giving a real transport
equation for the principal symbol for the evolution \pb{} (\ref{1.5}).

\par In particular, if ${{a_t}_\vert }_{C_{\phi _t}}$ has
constant argument for one value of $t$, the same holds for all
other values.

\par For a given $U$ associated to $\kappa $, choose $\kappa
_t$ and $U(t)$ as in (\ref{1.5}), so that 
$\kappa _0$ is linear, $U(1)=U$. Then using Proposition \ref{Prop1.2} 
and the above remark, we get the equivalences:
[$U$ satisfies (\ref{1.3}).] $\Leftrightarrow$ [$U(0)$ satisfies (\ref{1.3}).]
$\Leftrightarrow$ [The principal symbol of $U(0)$ has constant
argument.] $\Leftrightarrow$ [The principal symbol of $U$ has
constant argument.] This gives Proposition \ref{Prop1.1}.
\end{proof}

\par Let $X,Y$ be analytic \mfld{}s of dimension $n$ equipped
with  analytic integration densities $L(dx)=L_X(dx)$,
$L(dy)=L_Y(dy)$. Let 
$$\kappa :\, \Omega _Y\to \Omega _X$$
be a \canform{} (and \diffeo{}), \an{} for simplicity, where 
$$\Omega _Y\subset\subset T^*Y,\ \Omega _X\subset\subset T^*X,$$
are connected, open with smooth \bdy{}. We do not assume $\Omega
_{X},\Omega _Y$ to be simply connected, so we may have finitely
many closed cycles
$\gamma _1,...,\gamma _N\subset \Omega _Y$ which generate the
homotopy group of $\Omega _Y$. 

\par Let $S:L^2(X)\to H_\Phi (\widetilde{X})$, $T:L^2(Y)\to H_\Psi
(\widetilde{Y})$ be corresponding FBI-\tf{}s as in the appendix,
where
$\widetilde{X},\widetilde{Y}$ denote tubular complex \neigh{}s of
$X,Y$ and with associated \canform{}s:
$$\kappa _S:T^*X\cap\{ \vert \xi \vert <C\}\to \Lambda _\Phi ,\
\kappa _T:T^*Y\cap\{
\vert \eta \vert <C\}\to \Lambda _\Psi ,$$
where we equip $H_\Phi ,H_\Psi $ with the scalar products that
make $S,T$ unitary, and we can have $C>0$ as large as we like.
Choose $C$ large enough, so that $\kappa _S,\kappa _T$ are
well-defined on $\Omega _X,\Omega _Y$ respectively, and let 
$$\widetilde{\Omega }_X=\pi _x\kappa _S\Omega _X\subset
\widetilde{X},\ \widetilde{\Omega }_Y=\pi _y\kappa _T\Omega _Y
\subset\widetilde{Y}.$$
Let $\widetilde{\kappa }:{\Lambda _\Psi
}\to 
{\Lambda _\Phi }$ be the lift of
$\kappa $, so that $\widetilde{\kappa }=\kappa _S\circ \kappa
\circ \kappa _T^{-1}$. Here $\Lambda _{\Phi ,\Psi }$ are
restricted to $\widetilde{\Omega }_{X,Y}$: $\Lambda _\Psi =\{
(y,{2\over i}\partial _y\Psi) ;\, y\in \widetilde{\Omega }_Y\}$,
$\Lambda _\Phi =\{ (x,{2\over i}\partial _x\Phi );\,
x\in\widetilde{\Omega }_X\} $.

\par We shall define a multivalued "Floquet periodic" \fop{}
$U:L^2(Y)\to L^2(X)$ which is only microlocally defined from
$\Omega _Y$ to $\Omega _X$ and associated to $\kappa $. Requiring
that $U$ be microlocally unitary with the improved Egorov
property, we will see that we can have the Floquet periodicity:
\begeq\label{1.6}
\gamma _* U=e^{i\theta (\gamma )}U,\endeq
where $\gamma $ is a closed loop in $\Omega _Y$  joining some
point $\rho $ to itself, $U$ denotes the operator $U$ as it is
defined near $\rho $ and the \lhs{} of (\ref{1.6}) denotes the \op{}
obtained from $U$ by following the loop $\gamma $. We will then
achieve (\ref{1.6}) with $\theta (\gamma )=h^{-1}S(\gamma )+k(\gamma )\pi
/2$, where $S(\gamma )=\int_{\kappa \circ \gamma }\xi
dx-\int_\gamma  \eta dy$ is the difference of the actions of
$\kappa \circ \gamma $ and $\gamma $, and $k(\gamma )\in{\bf Z}$ is
a "Maslov index", both quantities depending only on the homotopy
class of $\gamma $. (Requiring only the unitarity of $U$, we
could take $\theta (\gamma )=S(\gamma )/h$.) 

\par When discussing the improved property (\ref{1.3}), recall from
\cite{HeSj2} and \cite{SjZw}, that on a \mfld{}
with a preferred positive density, we can define the Weyl symbol
of a 0-th order $h$-\pop{} modulo ${\cal O}(h^2)$ by taking the
ordinary Weyl symbol for some system of local coordinates
$x_1,..,x_n$ for which the preferred density reduces to the
Lebesgue measure. Clearly Proposition \ref{Prop1.1} extends to this situation. 

\par We first notice that if 
$$Vu(x)=h^{-{n+N\over 2}}\iint e^{i\phi (x,y,\theta
)/h}a(x,y,\theta ;h)u(y)dyd\theta $$
is an elliptic \fop{} with leading symbol $a_0(x,y,\theta )$ $\ne 0$ on $C_\phi $, 
then we can obtain $V^*V=1+{\cal O}(h)$ by
multiplying $a_0$ by a positive real-analytic \fu{}. 

\par The same remark applies to 
$$\widetilde{V}:H^{\rm loc}_\Psi (\widetilde{\Omega }_Y)\to
H_\Phi ^{\rm loc}(\widetilde{\Omega }_X),$$
if we represent $\widetilde{V}$ as in \cite{MeSj} by 
\begeq\label{1.7}
\widetilde{V}u(x)=h^{-n}\int e^{i\psi (x,y)/h}b(x,y;h) u(y)
e^{-2\Psi (y)/h}L(dy), \endeq
where $\psi (x,y)$ is the multivalued grad-periodic \fu{} near
$\pi _{x,y}\Gamma $, with 
$$\partial _{\overline{x},y}\psi =0,\,\,\partial
_{\overline{x},y}b =0 \hbox{ near }\pi _{x,y}(\Gamma ),$$
$$\partial _x\psi (x,y)={2\over i}\partial_x \Phi (x),\,\, \partial
_{\overline{y}}\psi (x,y)={2\over i}\partial _{\overline{y}}\Psi
(y)\hbox{ on }\pi _{x,y}(\Gamma ),$$
$$\Phi (x)+\Psi (y)+\Im \psi (x,y)\sim {\rm dist\,}((x,y),\pi
_{x,y}(\Gamma ))^2,$$
where $\Gamma $ denotes the graph of $\widetilde{\kappa }$.

\par Recall that $\Im \psi $  is single-valued, and that 
\begeq\label{1.8}
{\rm var}_{(\widetilde{\kappa}\circ \gamma ,\gamma )}\psi =
\int_{\widetilde{\kappa}\circ \gamma }\xi dx-\int_\gamma  \eta dy, \endeq
is the action difference, when $\gamma $ is a closed curve in
$\Lambda _\Psi $ and $(\widetilde{\kappa}\circ \gamma ,\gamma )$ denotes the
curve $t\mapsto (\widetilde{\kappa}(\gamma(t)),\gamma (t))$. Here we also
identify $\Lambda _\Psi ,\Lambda _\Phi $ with $\widetilde{\Omega
}_Y,\widetilde{\Omega }_X$ whenever so is convenient.

\par Thus after multiplying ${b_\vert}_{\pi _{x,y}(\Gamma )}$ by
a positive real-analytic \fu{}, we may assume that 
\begeq\label{1.9}
V^*V=1+{\cal O}(h).\endeq
In order to have the improved Egorov property, we further need
that locally on $\pi _{x,y}(\Gamma )$:
\begeq\label{1.10}
{\rm arg\,}b_0(x,y)=K(y)+{\rm Const.},\hbox{ (notice that
$x=x(y)$ on $\Gamma $),} \endeq
where $K(y)$ is a grad-periodic \fu{} on $\pi _{x,y}(\Gamma )$,
that we do not try to compute here, but whose existence we infer
from Proposition \ref{Prop1.1} and the computation of $\widetilde{V}$ as $S\circ
V\circ T^{-1}$, with $V$ written microlocally with a real phase as
in (\ref{1.1}). 

\par We can find $b_0$ satisfying (\ref{1.10})
everywhere if we accept that ${{b_0}_\vert}_{\pi _{x,y}(\Gamma )}$
is multivalued. More precisely, $K$ is not globally
well-defined on $\pi _{x,y}(\Gamma )\simeq \Omega _Y$, but
$\omega =dK$ is a well-defined closed real 1-form on $\Omega _Y$
and we can find ${{b_0}_\vert}_{\pi _{x,y}(\Gamma )}$, unique up
to a constant factor of modulus 1, such that (\ref{1.9}), (\ref{1.10}) hold,
though $b_0$ will be multivalued:
\begeq\label{1.11}
\gamma _* b_0=\exp (i\int_\gamma \omega )b_0,\endeq
where $\gamma _*b_0$ denotes the new locally defined symbol
obtained by following $b_0$ around the closed loop $\gamma $ in
$\pi _{x,y}(\Gamma )\simeq \Omega _Y$.\medskip

\begin{prop}\label{Prop1.3} We have $\int_\gamma \omega
=k(\gamma ){\pi \over 2}$ for some integer $k(\gamma )\in {\bf
Z}$, for every closed loop $\gamma \subset \pi _{x,y}(\Gamma
)$.\end{prop}
\begin{proof}  Let $\gamma $ be a closed loop and
cover $\gamma $ by small open topologically trivial sets
$\widetilde{\Omega }_0,\widetilde{\Omega
}_1,...,\widetilde{\Omega }_{N-1}$ with increasing index
corresponding to the orientation of $\gamma $ in the natural way.
Let $\widetilde{\Omega }_N=\widetilde{\Omega }_0$. Let $\Omega
_j$ be the corresponding regions in $\Omega _Y$. In $\Omega _j$,
we represent $V$ by 
\begeq\label{1.12}
V_ju(x)=h^{-{n+N_j\over 2}}\iint _{\theta \in {\bf
R}^{N_j}}e^{i\phi _j(x,y,\theta )/h}a_j(x,y,\theta ;h)
u(y)dyd\theta .\endeq
For a given point in $\Omega _j\cap \Omega _{j+1}$, we have $\phi
_j=\phi _{j+1}$, $a_{j+1}=r_{j+1,j}e^{i\alpha _{j+1,j}\pi
/2}a_j+{\cal O}(h)$ at the corresponding points in $C_{\phi _j}$,
$C_{\phi _{j+1}}$, provided that we require all the fiber-variable
dimensions $N_j$ to have the same parity. (Cf.~\cite{Ho}.) This last property is
easy to achieve since we can always add one fiber-variable. We
conclude that 
$$
\int_\gamma  \omega ={\pi \over 2}(\alpha _{1,0}+\alpha
_{2,1}+..+\alpha _{N,N-1}),
$$
and the proposition follows. \end{proof}
\par Take $V$ as above with $b=b_0$ in (\ref{1.7}), so that (\ref{1.9}),
(\ref{1.10}) hold. Put 
\begeq\label{1.13}
U=V(V^*V)^{-{1\over 2}}.\endeq
Then $\widetilde{U}=SUT^{-1}$ is of the form (\ref{1.7}) with
$b=b_0+{\cal O}(h)$. We have $U^*U=1$ and $U$ satisfies (\ref{1.3}).
Since the unitarization is a local operation which commutes with
multiplication by a constant factor of modulus 1, (\ref{1.11}) becomes
valid also for $b$:
\begeq\label{1.14}
\gamma _* b=e^{ik(\gamma ){\pi \over 2}}b.\endeq
Here we also used Proposition \ref{Prop1.3}.

\par Summing up, we get
\medskip

\begin{theo}\label{Th1.4} Under the assumptions above on
$\kappa $, we can find a microlocally defined multivalued \fop{}
$U$ associated to $\kappa $, and a corresponding lift
$\widetilde{U}=SUT^{-1}$ of the form {\rm (\ref{1.7})}, such that $U$ is
unitary: $U^*U=1+{\cal O}(e^{-1/(Ch)})$, satisfies the
improved Egorov property {\rm (\ref{1.3})}, and 
$$\gamma _*U=e^{i(S(\gamma )/h+k(\gamma )\pi /2)}U,$$
for every closed loop in $\Omega _Y$, where $k(\gamma )\in{\bf
Z}$ and 
$$S(\gamma )=\int_{\kappa \circ \gamma }\xi dx-\int_\gamma \eta
dy.$$\end{theo}

\section{Reduction by averaging along trajectories}\label{section2}
\setcounter{equation}{0}

\par Let $P$, $M$
be as in the introduction. We work in a \neigh{} of
$p^{-1}(0)\cap T^*M$, and recall that $P=P_\epsilon $ has the
semi-classical principal symbol
\begeq\label{2.1}
p_\epsilon =p+i\epsilon q+{\cal O}(\epsilon ^2),\endeq 
in a complex \neigh{} of $p^{-1}(0)\cap T^*M$. Let $G_0$ be an
\an{} \fu{} defined near $p^{-1}(0)\cap T^*M$ such that 
\ekv{2.2}
{H_pG_0=q-\langle q\rangle ,}
where $\langle q\rangle $ is the trajectory average, defined in
(\ref{0.14}). We may take
\ekv{2.3}
{G_0={1\over T(E)}\int_{-T(E)/2}^{T(E)/2}(1_{{\bf
R}_-}(t)(t+{T(E)\over 2})+1_{{\bf
R}_+}(t)(t-{T(E)\over 2}))q\circ \exp tH_p dt,}
on $p^{-1}(E)$. 

\par We replace ${\bf R}^4$ by the new IR-\mfld{} 
\ekv{2.4}
{\Lambda _{\epsilon G_0}=\exp (i\epsilon H_{G_0})({\bf R}^4),}
which is defined in a complex \neigh{} of $p^{-1}(0)\cap T^*M$.
Writing $(x,\xi )=\exp (i\epsilon H_{G_0})(y,\eta )$, and using
$\rho =(y,\eta )$ as real symplectic coordinates on $\Lambda
_{\epsilon G_0}$, we get 
\eekv{2.5}
{{{p_\epsilon }_\vert}_{\Lambda _{\epsilon G_0}}=p_\epsilon (\exp
(i\epsilon H_{G_0})(\rho ))} 
{=\sum_{k=0}^\infty  {(i\epsilon H_{G_0})^k\over k!}(p_\epsilon
)=p+i\epsilon \langle q\rangle +{\cal O}(\epsilon ^2).}

\par Iterating this procedure, or looking more directly for
$G(x,\xi ,\epsilon )$ as an \asy{} sum 
\ekv{2.6}
{G\sim \sum_0^\infty  \epsilon ^kG_k(x,\xi )}
in some complex \neigh{} of $p^{-1}(0)\cap T^*M$, we see that we can
find $G_1$, $G_2 \ldots\,$ such that if 
\ekv{2.7}
{\Lambda _{\epsilon G}=\exp (i\epsilon H_G)({\bf R}^4),}
and we again write $\Lambda _{\epsilon G}\ni (x,\xi )=\exp
(i\epsilon H_G)(y,\eta )$ and parametrize by the real variables
$(y,\eta )$, then 
\ekv{2.8}
{{{p_\epsilon }_\vert}_{\Lambda _{\epsilon G}}=p+i\epsilon
\langle q\rangle +\epsilon ^2q_2+\epsilon ^3q_3+...,}
where $q_j=\langle q_j\rangle $, $j\ge 2$. This means that we can
\tf{} $p_\epsilon $ to $p_\epsilon \circ \exp (i\epsilon H_G)$ in
such a way that we get a new leading symbol which Poisson
commutes with the unperturbed leading symbol. 

As is well-known in the \sa{} case, this construction can be
extended to the level of \op{}s, and we may develop this
globally in another paper. In the present work we will do it only
after a reduction to a torus-like situation.

\par After replacing $p_\epsilon $ by $p_\epsilon \circ \exp
(i\epsilon H_{G_0})$ and correspondingly $P_\epsilon $, by
$U_\epsilon ^{-1}\circ P_\epsilon \circ U_\epsilon $, where
$U_\epsilon $
is the \fop{} 
$U_\epsilon =e^{-{i\over h}i\epsilon G_0(x,hD_x)}=e^{{\epsilon
\over h}G_0(x,hD_x)}$ (defined microlocally near $p^{-1}(0)\cap
T^*M$), we may assume that our operator
$P_\epsilon $ is \ml{}ly defined near $p^{-1}(0)\cap T^*M$ and has
the $h$-principal symbol
\ekv{2.9}
{p_\epsilon =p+i\epsilon \langle q\rangle +{\cal O}(\epsilon ^2).}
This can be done in such a way that $P_{\epsilon =0}$
remains the original unperturbed \op{}.

\par Let $\gamma _0\subset p^{-1}(0)\cap T^*M$ be a closed
$H_p$-trajectory and assume that $T(0)$ is the minimal period of
$\gamma _0$. Let $g:{\rm neigh\,}(0,{\bf R})\to {\bf R}$ be the
\an{} \fu{} defined by 
\ekv{2.10}
{g'(E)={T(E)\over 2\pi },\ g(0)=0.}
Then $H_{g\circ p}=g'(p)H_p$ has a $2\pi $-periodic flow and the
same closed trajectories as $H_p$. Clearly $2\pi $ is the minimal
period of $\gamma _0$ when viewed as a $H_{g\circ p}$-trajectory.
\begin{prop}\label{Prop2.1} \it There exists an \an{} \ctf{}
$\kappa :{\rm neigh\,}(\{ \tau =x=\xi =0\} ,T^*(S_t^1\times {\bf
R}_x))\to {\rm neigh\,}(\gamma _0,T^*M)$, mapping $\{ \tau =x=\xi
=0\}$ onto $\gamma _0$, such that $g\circ p\circ \kappa =\tau
$.\end{prop}

\begin{proof} Fix a point $\rho _0\in \gamma _0$
and choose local symplectic coordinates $(t,\tau ;x,\xi )$
centered at $\rho _0$, with $g\circ p=\tau $. This means that 
\ekv{2.11}
{\{ \xi ,x\} =1,\ \{ t,x\} =\{ t,\xi \}=0}
\ekv{2.12}
{H_\tau t=1,\ H_\tau x=H_\tau \xi =0.}
Now extend the definition of $t,\tau ,x,\xi $ to a full \neigh{}
of $\gamma _0$, by putting $\tau =g\circ p$ and requiring
$t,x,\xi $ to solve (\ref{2.12}). Since the $H_\tau $-flow is $2\pi
$-periodic (with $2\pi $ as the minimal period) near $\gamma _0$,
we see that $x,\xi $ are well-defined singlevalued \fu{}s, while
$t$ becomes multivalued in such a way that it increases by $2\pi $
each time we make a loop in the increasing time direction. (\ref{2.11})
extends to a full \neigh{} of
$\gamma _0$. This is equivalent to the proposition.\end{proof}

\par Notice that
\ekv{2.13}
{p\circ \kappa =f(\tau ),} 
where $f:=g^{-1}$. From (\ref{2.9}) we infer that 
\ekv{2.14}
{p_\epsilon \circ \kappa =f(\tau )+i\epsilon \langle q\rangle
(\tau ,x,\xi )+{\cal O}(\epsilon ^2),}
for a new \fu{} $\langle q\rangle $ which is \indep{} of $t$ (and
obtained from the earlier one by composition with $\kappa $).

If we let the \fop{} $U$ quantize $\kappa $ as in section \ref{section1}, we
get a new \op{} $U^{-1}P_\epsilon U$ with leading semi-classical
symbol $p_\epsilon \circ \kappa $ as in (\ref{2.14}). (Here $P_\epsilon
$ is the new version of $P_\epsilon $; $P_{\epsilon ,{\rm
new}}=U^{-1}P_{\epsilon ,{\rm old}}U$.) 

\par Now write simply $p$, $p_\epsilon $, $P_\epsilon $ for the
transformed objects. Then 
$$
P_\epsilon =P(t,x,hD_{t,x},\epsilon;h)
$$
is the formal Weyl quantization of a symbol $P(t,x,\tau ,\xi
,\epsilon ;h)$ which has an \asy{} expansion (\ref{0.5}) in the space
of \hol{} \fu{}s in a fixed complex \neigh{} of $\{ \tau =x=\xi
=0\}$ in $T^*(\widetilde{S}^1\times {\bf C})$, with
$\widetilde{S}^1=S^1+i{\bf R}$, and we will use the same notation
as in section \ref{section0}. (An exact value of the new symbol 
$P(t,x,\tau,\xi,\epsilon ;h)$ cannot be easily defined, but we know how to
define it mod ${\cal O}(e^{-1/(Ch)})$. We shall ho\-we\-ver avoid 
using the full power of analytic \pop{}s, and content ourselves
with the knowledge of $P$ mod ${\cal O}(h^\infty )$.) 

\par Now look for $G^{(1)}=\epsilon G_1(t,\tau ,x,\xi )+\epsilon
^2G_2(t,\tau ,x,\xi )+\ldots$ such that 
$$p_\epsilon \circ \exp i\epsilon H_{G^{(1)}}=f(\tau )+i\epsilon
\langle q\rangle (\tau ,x,\xi )+{\cal O}(\epsilon ^2)$$
is \indep{} of $t$. Here the \lhs{} can be written
$$\sum_{k=0}^\infty  \frac{1}{k!}(i\epsilon H_{G^{(1)}})^kp_\epsilon ,$$
and we get 
\begin{eqnarray*}&&p_\epsilon +i\epsilon ^2H_{G_1}(f(\tau ))+{\cal
O}(\epsilon ^3)=\\
&&f(\tau )+i\epsilon \langle q\rangle (\tau ,x,\xi )-i\epsilon
^2f'(\tau ){\partial \over \partial t}G_1+{\cal O}(\epsilon
^2)+{\cal O}(\epsilon ^3),\end{eqnarray*}
where the ${\cal O}(\epsilon ^2)$ term is the same as in (\ref{2.14}).
It is clear that we can find $G_1$ so that the $\epsilon ^2$-term
in this expression is \indep{} of $t$. Looking at the ${\cal
O}(\epsilon ^3)$-term we then determine $G_2$ and so on. (In this
construction, we could have applied $\kappa $ at the very
beginning before replacing $q$ by $\langle q\rangle $ by
averaging, and then incorporated $G_0$ into the expression
$G=G_0+\epsilon G_1+\ldots\,$, and as already indicated, this
could also have been done entirely (and in a full \neigh{} of
$p^{-1}(0)\cap T^*M$), before applying $\kappa $.) 

\par After replacing $p_\epsilon $ by $p_\epsilon \circ \exp
i\epsilon H_{G^{(1)}}$, we are now reduced to the case when 
\ekv{2.15}
{p_\epsilon =f(\tau )+i\epsilon \langle q\rangle (\tau ,x,\xi
)+{\cal O}(\epsilon ^2)}
is \indep{} of $t$.

\par Finally we remove the $t$-dependence from the lower order
terms. After conjugating $P_\epsilon $ by a \fop{} $V_\epsilon $,
which quantizes $\exp i\epsilon H_{G^{(1)}}$, we may assume that
$p_\epsilon $ in (\ref{2.15}) is the principal symbol of $P_\epsilon $
(and that it is \indep{} of $t$). Look for an $h$-\pop{}
$A(t,x,hD_{t,x},\epsilon ;h)$ with symbol
\ekv{2.16}
{A(t,x,\tau ,\xi ,\epsilon ;h)\sim \sum_{k=1}^\infty 
a_k(t,x,\tau ,\xi ,\epsilon )h^k, }
such that the full (Weyl) symbol of 
\ekv{2.17}
{e^{{i\over h}A}P_\epsilon e^{-{i\over h}A}=e^{{i\over h}{\rm
ad}_A}P_\epsilon =\sum_{k=0}^\infty  \frac{1}{k!} ({i\over h}{\rm
ad}_A)^kP_\epsilon }
is \indep{} of $t$. Since $A={\cal O}(h)$ we know that ${i\over
h}{\rm ad}_A$ lowers the order in $h$ by one (with the convention
that a symbol $={\cal O}(h^{-j})$ is of order $j$), so (\ref{2.17})
makes sense \asy{}ally. The subprincipal symbol of (\ref{2.17}) is 
$$h(p_{1,\epsilon }(x,\xi )+\{ p_\epsilon ,a_1\}
)=h(p_{1,\epsilon }(x,\xi )+f'(\tau ){\partial \over \partial
t}a_1(t,\tau ,x,\xi ,\epsilon )+{\cal O}(\epsilon )),$$
and we make this \indep{} of $t$ by successively determining the
\coef{}s in the \asy{} series
$$a_1(t,\tau ,x,\xi ,\epsilon )=\sum_{j=0}^\infty  a_{1,j}(t,\tau
,x,\xi )\epsilon ^j.$$
After that we return to (\ref{2.17}) and see that the construction of
$a_2$, $a_3$, ... is essentially the same. 

\par Actually, we do not have to do this construction in 2 steps,
and we can view $\eps G^{(1)}$ above as (a constant factor times) the
leading symbol $a_0={\cal O}(\epsilon ^2)$ in 
\ekv{2.18}
{A\sim \sum_{k=0}^\infty  a_k(t,x,\tau ,\xi ,\epsilon )h^k, }
such that if $P_\epsilon $ denotes the very first \op{} we get on
$S^1\times {\bf R}$, then the \lhs{} of (\ref{2.17})  has a symbol
which is well-defined as an \asy{} series in $(\epsilon ,h)$ and
is \indep{} of $t$. This can be seen by first determining $a_0$
from (\ref{2.17}) (leading to a repetition of what we already did) and
then the other terms. (When $\epsilon $ is small but fixed, the
\pb{} becomes more subtle and the break-up into two steps is more
natural, with the first step being the one containing the new
difficulties.)

\par Summing up the discussion of this section, we have
\begin{prop}\label{Prop2.2} Let $P$, $M$ be as in section
{\rm \ref{section0}}. Let $\gamma _0\subset p^{-1}(0)\cap T^*M$ be a closed
$H_p$-trajectory where $T(0)$ is the minimal period and let
$\kappa $ be the \ctf{} of \Pro{} {\rm \ref{Prop2.1}}. Let $U$ be a
corresponding elliptic \fop{} as in section {\rm \ref{section1}}. Then there exist
$G(x,\xi ,\epsilon )$ (\indep{} of $\gamma _0 ,\kappa ,U$) with
the
\asy{} expansion {\rm (\ref{2.6})} in the space of \hol{} \fu{}s in some fixed
complex \neigh{} of
$p^{-1}(0)\cap T^*M$ and a symbol $A(t,x,\tau ,\xi ,\epsilon ;h)$
as in {\rm (\ref{2.16})}, where 
\ekv{2.19}
{a_k\sim \sum_{j=0}^\infty  a_{k,j}(t,x,\tau ,\xi )\epsilon ^j}
in the space of \hol{} \fu{}s in a fixed complex \neigh{} of
$\tau =x=\xi =0$ in $T^*(\widetilde{S}^1\times {\bf C}))$, such
that if $G,A$ also denote the corresponding Weyl quantizations, the
\op{}
\ekv{2.20}
{\widetilde{P}_\epsilon=e^{{i\over h}A}U^{-1}
e^{-{\eps \over h}G} P_\epsilon e^{{\eps\over h}G}U e^{-{i\over h}A}={\rm
Ad}_{e^{{i\over h}A}U^{-1}e^{-{\eps \over h}G}}P_\epsilon }
has a symbol
 \ekv{2.21}
{\widetilde{P}_\epsilon (x,\tau ,\xi ,\epsilon ;h)\sim
\sum_0^\infty \widetilde{p}_k(x,\tau ,\xi ,\epsilon )h^k}
\indep{} of $t$. Here each $\widetilde{p}_k=\widetilde{p}_k(x,\tau
,\xi ,\epsilon )\sim \sum_{j=0}^\infty 
\widetilde{p}_{k,j}(x,\tau ,\xi )\epsilon ^j$
in the space of \hol{} \fu{}s in a fixed complex \neigh{} of
$\tau ,x,\xi =0$. Moreover 
\ekv{2.22}
{
\widetilde{p}_{0,\epsilon }=f(\tau )+i\epsilon \langle q\rangle
(\tau ,x,\xi )+{\cal O}(\epsilon ^2). }
\end{prop} 

\par If $\langle q\rangle $ has a \nondeg{}
extreme value along $\gamma _0$, then the
proposition is directly applicable (see section \ref{section4}), while in other
situations (such as in section \ref{section3}), it is not
global enough.              

\section{Normal forms and quasi-\ev{}s in the to\-rus case}\label{section3}
\setcounter{equation}{0}

\par Let $P,M,p,q,\langle q\rangle ,\Lambda_{0,F_0}$ be as in
section \ref{section0}. After replacing $q$ by $q-F_0$, we may assume that
$F_0=0$, so we consider
\ekv{3.1}
{\Lambda _{0,0}:\,p=0,\, \Re \langle q\rangle =0.}
Notice that $\Lambda _{0,0}$ is invariant under the $H_p$-flow.
We assume that $T(0)$ is the minimal period for all the closed
trajectories in $\Lambda _{0,0}$ and that
\ekv{3.2}
{
dp,d\langle \Re q\rangle \hbox{ are \indep{} at the points of
}\Lambda _{0,0}, }
so that $\Lambda _{0,0}$ is a Lagrangian \mfld{} and also a union
of tori. Assume for simplicity that $\Lambda _{0,0}$ is
connected, so that it is equal to one single Lagrangian torus. In
this section we work microlocally near $\Lambda _{0,0}$ and
proceed somewhat formally. In section \ref{section5} we follow up with
suitable function spaces and see how to justify the computation
of the spectrum via a global Grushin \pb{}. We have seen that we
can reduce ourselves to the case when 
\ekv{3.3}{p_\epsilon =p+i\epsilon \langle q\rangle +{\cal
O}(\epsilon ^2).}
Assume from now on that $\langle q\rangle $ is real-valued or
more generally that $\langle q\rangle $ is a \fu{} of $p$ and
$\Re \langle q\rangle $. 
We can make a real \canform{}
\ekv{3.4}
{
\kappa :{\rm neigh\,}(\xi =0,T^*{\bf T}^2)\to {\rm
neigh\,}(\Lambda _{0,0},T^*M), \,\, {\bf T}^2=({\bf R}/2\pi {\bf Z})^2,}
such that $p\circ \kappa =p(\xi _1)$, 
$\langle q\rangle \circ \kappa =\langle q\rangle (\xi )$ (with a
slight abuse of notation). 

\par Recall that this can be done in the following way: Let
$\Lambda _{E,F}$ be the Lagrangian torus given by $p=E,\Re \langle{q}\rangle =F$, for
$(E,F)\in {\rm neigh\,}(0,{\bf R}^2)$. Let $\gamma _1(E,F)$ be the
cycle in
$\Lambda_{E,F}$ corresponding to a closed $H_p$-trajectory with
minimal period, and let $\gamma _2(E,F)$ be a second cycle so that
$\gamma _1,\gamma _2$ form a fundamental system of cycles on the
torus
$\Lambda_{E,F}$. Necessarily $\gamma _2$ maps to the simple loop
given by
$\Re \langle q\rangle =F$ in the abstract quotient manifold
$p^{-1}(E)/{\bf R}H_p$. Now it is classical (see~\cite{Ar}) that we
can find a real analytic canonical transformation $\kappa : {\rm
neigh\,}(\eta =0,T^*{\bf T}^2)\ni (y,\eta )\mapsto
(x,\xi )\in {\rm neigh\,}(\Lambda _{0,0},T^*M)$, such that  
$$\eta _j={1\over 2\pi }(\int _{\gamma _j(E,F)}\xi dx-\int_{\gamma
_j(0,0)}\xi dx),$$
where $E,F$ depend on $(x,\xi )$ and are
determined by $(x,\xi )\in \Lambda _{E,F}$, i.e. by
$E=p(x,\xi ),F=\Re \langle q\rangle (x,\xi )$.

\par Let us also recall that this can be done as follows: 
We start by taking a first canonical transformation 
$\kappa _0: {\rm neigh\,} (\xi =0,T^*{\bf T}^2) \to 
{\rm neigh\,}(\Lambda _{0,0},T^*M)$ such that the zero section
is mapped to $\Lambda _{0,0}$. Then using $\kappa
_0$, we can consider
$p,\langle q\rangle$ as living on $T^*{\bf T}^2$.
$\Lambda _{E,F}$ is then given by
$$\xi =\phi _x',\ \phi =\phi _{\rm per}(x,E,F)+{\eta _1}x_1 + {\eta
_2}x_2,\hbox{ with }\det \phi ''_{x,(E,F)}\ne 0,$$ with 
$\eta _j=\eta _j(E,F)$ as
above  (now being the actions$/2\pi $ with respect to
$\xi dx$), and $\phi _{{\rm per}}$ being $(2\pi {\bf Z})^2$-periodic.
Moreover,
$\phi _x'(x,\eta )=0$, $\eta =0$ for
$E=F=0$. It is easy to check, using that our functions are real-valued, that 
$(E,F)\mapsto (\eta _1(E,F),\eta _2(E,F))$ is a local diffeomorphism,
so we can use
$\eta _1,\eta _2$ as new parameters replacing $E,F$, and write
$\phi =\phi (x,\eta )$. Consider
$$\kappa _1:({\partial \phi \over \partial \eta },\eta )\mapsto
(x,{\partial \phi
\over \partial x})$$ which maps the zero section to itself. Then
$\kappa :=\kappa _0\circ \kappa _1$ has the required
properties.

\par Let $U$ be a corresponding \fop{}, implementing $\kappa $, so
that if we denote by $P_\epsilon $ also the conjugated \op{}
$U^{-1}P_\epsilon U$, we have a new \op{} with leading symbol 
\ekv{3.5}
{p_\epsilon =p(\xi _1)+i\epsilon \langle q\rangle (\xi )+{\cal
O}(\epsilon ^2).}
For the conjugated \op{}, we still have the property that
$P_{\epsilon =0}$ is \sa{}. From the assumption (\ref{3.2}) about linear
independence, we get 
\ekv{3.6}
{
\partial _{\xi _1}p(0)\ne 0,\ \partial _{\xi _2}\Re \langle
q\rangle (0)\ne 0. }
As in the preceding section, we can find an $h$-\pop{} $A$ with
symbol
$\sum_{\nu =0}^\infty  h^\nu  a_\nu (x,\xi ,\epsilon )$, 
$a_0={\cal O}(\epsilon^2 )$, such that formally 
\ekv{3.7}
{e^{{i\over h}A}P_\epsilon  e^{-{i\over h}A}=e^{{i\over h}{\rm
ad}_A}(P_\epsilon )=\sum_{k=0}^\infty  {1\over k!}({i\over h}{\rm
ad}_A)^k(P_\epsilon )=:\widetilde{P}_\epsilon ,}
with $\widetilde{P}_{\epsilon }(x,\xi ,\epsilon ;h)$ \indep{} of
$x_1$, and leading symbol 
$$\widetilde{p}_{\epsilon }=p(\xi _1)+i\epsilon \langle
q\rangle (\xi )+{\cal O}(\epsilon ^2)$$
also \indep{} of $x_1$. We recall that the symbol $A(x,\xi
,\epsilon ;h)$ is a formal power series both in $\epsilon $ and
$h$ with \coef{}s all \hol{} in the same complex \neigh{} of $\xi
=0$. This construction can be done in such a way that
$\widetilde{P}_{\epsilon =0}$ is \sa{}.

\par We next look for a further conjugation that eliminates the
$x_2$-dependence in the symbol.
\smallskip

\par\noindent a) We start by considering the general case, when
the subprincipal symbol of $P_{\epsilon =0}$
is not necessarily $0$, so that the complete symbol of
$\widetilde{P}_\epsilon $ takes the form 
\ekv{3.8}
{
\widetilde{P}_\epsilon (x_2,\xi ;h)=\sum_{\nu =0}^\infty  h^\nu
\widetilde{p}_\nu (x_2,\xi ,\epsilon ), }
with 
\ekv{3.9}
{
\widetilde{p}_0(x_2,\xi ,\epsilon )=\widetilde{p}_\epsilon =p(\xi
_1)+i\epsilon \langle q\rangle (\xi )+{\cal O}(\epsilon ^2), }
and $\widetilde{p}_1(x_2,\xi ,0)$ not necessarily identically
equal to $0$. 

\par The easiest case is when $h/\epsilon \le {\cal O}(h^{\delta
_1})$ for some $\delta _1>0$, so that we can consider 
$h/\epsilon $ as an \asy{}ally small parameter. Look for 
\ekv{3.10}
{B(x_2,\xi ,\epsilon ,{h\over \epsilon },h)=\sum_{\nu =0}^\infty 
h^\nu b_\nu (x_2,\xi ,\epsilon ,{h\over \epsilon }),}
with $b_{\nu}={\cal O}(\epsilon +h/\epsilon )$, such that on the
operator level (with $hD_x$ instead of
$\xi
$), 
\ekv{3.11}
{
e^{{i\over h}B}\widetilde{P}_\epsilon e^{-{i\over
h}B}=:\widehat{P}_\epsilon (hD_x,\epsilon ,{h\over \epsilon },h)
}has a symbol \indep{} of $x$. Notice that $B(x_2,hD_x,\epsilon
;h)$ and $p(hD_{x_1})$ commute. On the symbol level we
write
\begin{eqnarray}\label{3.12}
& & \widetilde{P}_\epsilon = p(\xi _1)+
\epsilon (i\langle
q\rangle (\xi )+{\cal O}(\epsilon )+{h\over \epsilon }\widetilde{p}_1(x_2,\xi
,\epsilon )+h{h\over \epsilon }\widetilde{p}_2(x_2,\xi ,\epsilon
)+\ldots )\\
& & = p(\xi _1)+\epsilon (r_0(x_2,\xi ,\epsilon
,{h\over \epsilon })+hr_1(x_2,\xi ,\epsilon ,{h\over \epsilon
})+\ldots),\nonumber\end{eqnarray} 
with 
$$r_0(x_2,\xi ,\epsilon ,{h\over \epsilon })=i\langle q\rangle
(\xi )+{\cal O}(\epsilon )+{h\over \epsilon }\widetilde{p}_1=i\langle
q\rangle (\xi )+{\cal O}(\epsilon )+{\cal O}({h\over \epsilon
}),$$
$$r_1={h\over \epsilon }\widetilde{p}_2(x_2,\xi ,\epsilon ),\, \ldots $$
Notice that $r_j={\cal O}(h/\epsilon )$ for $j\ge 1$. We shall
treat $h/\epsilon $ as an independent parameter.

\par We use this and develop (\ref{3.11}) to get, with ${\rm ad}_bc$
denoting the symbol of  
$$
{\rm ad}_{b(x,hD_x)}c(x,hD_x)=[b(x,hD),c(x,hD)],
$$
$$
p(\xi _1)+\epsilon \sum_{k=0}^\infty \sum_{j_1=0}^\infty
...\sum_{j_k=0}^\infty \sum_{\ell =0}^\infty  h^{\ell
+j_1+..+j_k}{1\over k!}({i\over h}{\rm ad}_{b_{j_1}})..({i\over
h}{\rm ad}_{b_{j_k}})r_\ell=p(\xi _1)+\epsilon \sum_{n=0}^
\infty h^n\widehat{r}_n, $$
with $\widehat{r}_n$ being equal to the sum of all coefficients for
$h^n$
resulting from all the expressions
\ekv{3.13}
{h^{\ell +j_1+..+j_k}{1\over k!}({i\over h}{\rm
ad}_{b_{j_1}})...({i\over h}{\rm ad}_{b_{j_k}})r_\ell ,}
with $\ell+j_1+...+j_k\le n$.

\par The first term is 
$$
\widehat{r}_0=\sum {1\over k!}H_{b_0}^kr_0=r_0 \circ \exp (H_{b_0}),
$$
where we want $\widehat{r}_0$ to be independent of $x_2$ (in
addition to $x_1$). We get with $b_0={\cal O}(\epsilon
+h/\epsilon )$:
\ekv{3.14}{\widehat{r}_0=i\langle q\rangle (\xi )+{\cal O}(\epsilon
+h/\epsilon )-i\partial _{\xi _2}\langle q\rangle \partial
_{x_2}b_0 +{\cal O}((\epsilon ,{h\over \epsilon })^2),}
and using that $\partial _{\xi _2}\langle q\rangle \ne 0$, it is
clear how to construct $b_0={\cal O}(\epsilon +h/\epsilon )$ as a
formal Taylor series in $\epsilon ,h/\epsilon $, so that
$\widehat{r}_0=i\langle q\rangle (\xi )+{\cal O}(\epsilon
+h/\epsilon )$ is \indep{} of $x$ (modulo a term ${\cal
O}(h^\infty )$).

\par Assume for simplicity that the conjugation by $e^{{i\over
h}b_0(x_2,hD_x,\epsilon, h/\epsilon)}$ has already been carried out, so that
we are reduced to the case when $r_0=i\langle q\rangle (\xi
)+{\cal O}(\epsilon +h/\epsilon )$ is \indep{} of $x_2$, and
$r_j={\cal O}(\epsilon +h/\epsilon )$ for $j\ge 1$. Then look for
a new conjugation $\exp {i\over h}{\rm ad}_B$, with $B(x_2,\xi
,\epsilon ,h/\epsilon ;h)=\sum_{\nu =1}^\infty  h^\nu b_\nu
(x_2,\xi ,\epsilon ,{h\over \epsilon })$.
The new expression for the \lhs{} of (\ref{3.11}) becomes
\ekv{3.15}
{p(\xi _1)+\epsilon \sum_{k=0}^\infty \sum_{j_1=1}^\infty
...\sum_{j_k=1}^\infty \sum_{\ell =0}^\infty  h^{\ell
+j_1+..+j_k}{1\over k!}({i\over h}{\rm ad}_{b_{j_1}})..({i\over
h}{\rm ad}_{b_{j_k}})r_\ell=p(\xi _1)+\epsilon \sum_{n=0}^
\infty h^n\widehat{r}_n,}
with $\widehat{r}_n$ equal to the sum of all coefficients for
$h^n$ resulting from the expressions (\ref{3.13}) with $\ell
+j_1+..+j_k\le n$ \it and \rm $j_\nu \ge 1$. Then 
$\widehat{r}_0=r_0$,
$\widehat{r}_1=r_1+H_{b_1}r_0=r_1-H_{r_0}b_1$, ..,
$\widehat{r}_n=r_n-H_{r_0}b_n +s_n$, where $s_n$ only depends on
$b_1,...,b_{n-1}$
and is the sum of all \coef{}s of $h^n$
arising in the expressions (\ref{3.13}) with $\ell +j_1 +..+j_k\le n$,
$j_1,..,j_k,\ell <n$, $j_\nu \ge 1$. 

\par It is therefore clear how to find $b_1,b_2,\ldots $ successively
with $b_j={\cal O}(\epsilon +h/\epsilon )$, such that all the
$\widehat{r}_j$
are \indep{} of $x$
and $={\cal O}(\epsilon +h/\epsilon )$. This completes the proof
of (\ref{3.11}).
\par Summing up the discussion so far, if we do not make
any assumption on the subprincipal symbol of $P_{\epsilon =0}$ and
restrict the attention to $h/\epsilon \le {\cal O}(h^{\delta
_1})$ for some $\delta _1>0$, then we can find 
$$B_0=b_0(x_2,hD_x,\epsilon ,h/\epsilon ),\ b_0={\cal O}(\epsilon
+h/\epsilon ),$$ and 
$$B_1=\sum_{\nu =1}^\infty  b_{\nu}(x_2,hD_x,\epsilon ,h/\epsilon
)h^{\nu},\ b_{\nu}={\cal O}(\epsilon +h/\epsilon ),$$such that 
\ekv{3.18}
{
\widehat{P}_\epsilon :=e^{{i\over h}{\rm ad}_{B_1}}
e^{{i\over h}{\rm ad}_{B_0}}\widetilde{P}_\epsilon }
has a symbol \indep{} of $x$:
\ekv{3.19}{\widehat{P}_\epsilon =p(\xi _1)+\epsilon (r_0(\xi
,\epsilon ,{h\over \epsilon })+hr_1(\xi ,\epsilon ,{h\over
\epsilon })+\ldots),} with
$r_0=i\langle q\rangle (\xi )+{\cal O}(\epsilon +h/\epsilon )$,
and $r_\nu ={\cal O}(\epsilon +h/\epsilon )$ for
$\nu \ge 1$. 

\par Remaining in the general case, without any assumption on the
lower order terms, we now assume merely that $h/\epsilon \le
\delta _0$ for some \sufly{} small $\delta _0>0$. This means that
we can no longer construct $b_0$ by a formal Taylor series in
$h/\epsilon $, and we shall replace $e^{{i\over
h}b_0(x_2,hD_x,\epsilon ,h/\epsilon)}$ by a \fop{}, constructed
directly.

\par Look for $\phi =\phi (x_2,\xi ,\epsilon ,h/\epsilon )$
solving
\ekv{3.20}
{
r_0(x_2,\xi _1,\xi _2+\partial _{x_2}\phi ,\epsilon ,{h\over
\epsilon })=\langle r_0(\cdot ,\xi ,\epsilon ,{h\over \epsilon
})\rangle , }
where $\langle \cdot \rangle $ denotes the average \wrt{} $x_2$.
By the implicit \fu{} theorem, (\ref{3.20}) has a solution with $\partial
_{x_2}\phi $ singlevalued and ${\cal O}(\epsilon +h/\epsilon
)$. If we Taylor expand (\ref{3.20}), we get
$$(\partial _{\xi _2}r_0)(x_2,\xi ,\epsilon ,{h\over \epsilon
})\partial _{x_2}\phi +(r_0(x_2,\xi ,\epsilon ,{h\over \epsilon
})-\langle r_0(\cdot ,\xi ,\epsilon ,{h\over \epsilon })\rangle
)={\cal O}(({h\over \epsilon },\epsilon )^2),$$
and using also that 
$$\partial _{\xi _2}r_0(x_2,\xi ,\epsilon ,{h\over \epsilon
})=i\partial _{\xi _2}\langle q\rangle (\xi )+{\cal O}(\epsilon
+{h\over \epsilon }),$$
we get,
$$\phi =\phi _{{\rm per}}+x_2\zeta _2,$$
with $\zeta _2=\zeta _2(\xi ,\epsilon ,{h\over \epsilon })={\cal
O}((\epsilon ,h/\epsilon )^2)$, and $\phi _{\rm per}={\cal
O}((\epsilon ,h/\epsilon ))$ periodic in $x_2$. Put $\eta =\eta
(\xi ,\epsilon ,h/\epsilon )=(\xi _1,\xi _2+\zeta _2)$, and 
$$\psi (x,\eta ,\epsilon ,{h\over \epsilon })=x\cdot \eta +\phi
_{\rm per},$$
where $\phi _{\rm per}$ is viewed as a \fu{} of $\eta $ rather
than $\xi $. 

\par Consider the \canform{}
$$\kappa :(\psi _\eta ',\eta )\mapsto (x,\psi _x'),$$ which is
$(\epsilon +h/\epsilon )$-close to the identity and can be viewed
as a \fy{} of \tf{}s depending analytically on the parameter $\xi
_1$. With
$\xi =\xi (\eta ,\epsilon ,{h\over \epsilon })$, we have by
construction:
\ekv{3.21}
{
(r_0\circ \kappa )(y,\eta
,\epsilon ,{h\over \epsilon })= \langle r_0(\cdot ,\xi ,\epsilon
,{h\over \epsilon })\rangle =\langle r_0(\cdot ,\eta ,\epsilon
,{h\over \epsilon })\rangle +{\cal O}(\epsilon ^2+({h\over \epsilon
})^2),} and this is a \fu{} of $(y,\eta )$ which is \indep{} of $y$.
Notice that $p(\xi _1)$ is unchanged under composition with 
$\kappa $.

\par We can quantize $\kappa $ as a \fop{} $U$ and after
conjugation by this \op{}, we may assume that we have a new \op{}
$\widetilde{P}_\epsilon $ as in (\ref{3.12}) with
$r_0=i\langle q\rangle (\xi )+{\cal O}(\epsilon +h/\epsilon )$
\indep{} of $x$ and with $r_j={\cal O}(\epsilon +h/\epsilon )$ . 

\par As before, we can then make a further conjugation $e^{{i\over
h}{\rm ad}_{B_1}}$ in order to remove the $x$-dependence
completely and the conclusion is that if we make no assumption on
the subprincipal symbol and restrict the attention to $h/\epsilon
\le \delta _0$, for $\delta _0>0$ small enough, then we can find
a  \fop{},
\ekv{3.22}
{U^{-1}u(x;h)={1\over (2\pi h)^2}\iint e^{{i\over h}(\psi (x,\eta
)-y\cdot \eta )} a(x,\eta;h) u(y)dyd\eta ,}
with $\psi (x,\eta )=x\cdot \eta +\phi _{\rm per}(x_2,\eta
,\epsilon ,h/\epsilon )$, $\phi _{\rm per}={\cal O}(\epsilon
+h/\epsilon )$, and 
$$B_1=\sum_{\nu =1}^\infty  b_{\nu}(x_2,hD_x,\epsilon ,{h\over
\epsilon })h^{\nu},\ b_{\nu}={\cal O}(\epsilon +{h\over \epsilon }),$$
such that 
$$\widehat{P}_\epsilon :=e^{{i\over h}{\rm ad}_{B_1}}{\rm Ad}_U
\widetilde{P}_\epsilon $$
has a symbol \indep{} of $x$ as in (\ref{3.19}), with the same estimates
as there.

\smallskip
\par\noindent b) We now assume that in the original problem,
$P_{\epsilon =0}$ has subprincipal symbol $0$. Then after a first
time averaging, transportation to the torus, and the elimination of the 
$x_1$-dependence, we may assume that 
\ekv{3.23}
{\widetilde{P}(x_2,\xi ,\epsilon ;h)=\sum_{\nu =0}^\infty h^\nu \widetilde{p}_\nu (x_2,\xi
,\epsilon ), }
with $\widetilde{p}_0$ \indep{} of $x$ mod ${\cal O}(\epsilon ^2)$:
\ekv{3.24}
{
\widetilde{p}_0(x_2,\xi ,\epsilon )=p(\xi _1)+i\epsilon \langle q\rangle (\xi
)+{\cal O}(\epsilon ^2), }
\ekv{3.25}
{
\widetilde{p}_1(x_2,\xi ,0)=0.
}
(Recall from section \ref{section1} and the references given there, that the
\canform{}s can be quantized in such a way that Egorov's theorem
holds modulo
${\cal O}(h^2)$.)  In analogy with (\ref{3.12}), we have with
$\widetilde{p}_1(x_2,\xi  ,\epsilon )=\epsilon q_1(x_2,\xi ,\epsilon )$,
\begin{eqnarray}\label{3.26}
\widetilde{P}_\epsilon & = & p(\xi _1)+\epsilon (i\langle q\rangle
(\xi )+{\cal O}(\epsilon )+hq_1(x_2,\xi ,\epsilon )+{h^2\over
\epsilon }\widetilde{p}_2+h{h^2\over \epsilon }\widetilde{p}_3+\ldots)\\
& = & p(\xi _1)+\epsilon (r_0(x_2,\xi ,\epsilon ,{h^2\over \epsilon
})+hr_1(x_2,\xi ,\epsilon ,{h^2\over \epsilon })+h^2r_2+\ldots),
\nonumber\end{eqnarray}
with 
\begin{eqnarray*}
r_0(x_2,\xi ,\epsilon ,{h^2\over \epsilon }) & = & i\langle q\rangle
(\xi )+{\cal O}(\epsilon )+{h^2\over \epsilon }\widetilde{p}_2, \\
r_1(x_2,\xi ,\epsilon ,{h^2\over \epsilon })& = &q_1(x_2,\xi
,\epsilon )+{h^2\over \epsilon }\widetilde{p}_3, \\
 r_2(x_2,\xi ,\epsilon ,{h^2\over\epsilon })& = & {h^2\over \epsilon
}\widetilde{p}_4, \ldots
\end{eqnarray*}

\par We first consider the case when 
\ekv{3.27}
{{h^2 \over \epsilon }\le h^{\delta _1},}
for some fixed $\delta _1>0$. A first conjugation by $e^{{i\over
h}b_0(x_2,hD_x,\epsilon ,{h^2\over \epsilon })}$, with $b_0={\cal
O}(\epsilon +h^2/\epsilon )$, allows us to make $r_0$ \indep{} of
$x_2$, and we still have (\ref{3.26}) with $r_j={\cal O}(1)$ for $j\ge
1$.

\par Then we look for a new conjugation $\exp {i\over h}{\rm
ad}_{B_1}$ with 
\ekv{3.28}
{
B_1(x_2,\xi ,\epsilon ,{h^2\over \epsilon };h)=\sum_{\nu
=1}^\infty  h^\nu  b_\nu (x_2,\xi ,\epsilon ,{h^2\over \epsilon
}). }
The conjugated \op{} (\ref{3.11}) can be expanded as in (\ref{3.15}) and as
after that \e{} it is clear how to get $b_\nu ={\cal O}(1)$ for
$\nu \ge 1$, such that the resulting $\widehat{r}_n$ are
\indep{} of $x_2$, with $\widehat{r}_0(\xi ,\epsilon
,h^2/\epsilon )=i\langle q\rangle (\xi )+{\cal O}(\epsilon
+h^2/\epsilon )$.

\par Summing up the discussion so far, if we assume that the
subprincipal symbol of $P_{\epsilon =0}$ vanishes, and restrict
the attention to the range (\ref{3.27}) for some fixed $\delta _1>0$,
then we can find $B_0=b_0(x_2,hD_{x},\epsilon ,{h^2\over \epsilon })$ with 
$b_0={\cal O}(\epsilon +{h^2\over \epsilon })$
and $B_1(x_2,hD_x,\epsilon ,{h^2\over \epsilon };h)$ with symbol
(\ref{3.28}), and $b_\nu ={\cal O}(1)$, such that 
$$e^{{i\over h}{\rm ad}_{B_1}}e^{{i\over h}{\rm
ad}_{B_0}}\widetilde{P}_\epsilon =\widehat{P}_\epsilon $$
has the symbol
\ekv{3.29}
{
p(\xi _1)+\epsilon (r_0(\xi ,\epsilon ,{h^2\over \epsilon
})+hr_1(\xi ,\epsilon ,{h^2\over \epsilon })+\ldots) }
\indep{} of $x$ and with 
\ekv{3.30}
{
r_0=i\langle q\rangle (\xi )+{\cal O}(\epsilon +{h^2\over
\epsilon }),\ r_\nu ={\cal O}(1),\,\, \nu \ge 1 .}

\par If we replace (\ref{3.27}) by the weaker assumption,
\ekv{3.31}
{
{h^2 \over \epsilon}\le \delta _0,\ \delta _0\ll 1,
}
then again we have to replace the conjugation by $e^{{i\over
h}B_0}$ by that by a \fop{} constructed as earlier: We solve
(\ref{3.20}) (with $h/\epsilon $ replaced by $h^2/\epsilon $) and get
$\partial _{x_2}\phi $ single-valued and ${\cal O}(\epsilon
+h^2/\epsilon )$.

\par Taylor expanding (\ref{3.20}) and using that 
$$\partial _{\xi _2}r_0(x_2,\xi ,\epsilon ,{h^2\over \epsilon
})=i\partial_{\xi_2}\langle q\rangle (\xi )+{\cal O}(\epsilon +{h^2\over \epsilon
}),$$
we get 
$$\phi =\phi _{\rm per}+x_2\zeta _2,$$
with $\zeta _2=\zeta _2(\xi ,\epsilon ,{h^2\over \epsilon
})={\cal O}((\epsilon ,{h^2\over \epsilon })^2)$
and $\phi _{\rm per}={\cal O}(\epsilon +h^2/\epsilon )$ periodic
in $x_2$. Again we put $\eta =\eta (\xi ,\epsilon ,h^2/\epsilon
)=(\xi _1,\xi _2+\zeta _2)$ and 
$$\psi (x,\eta ,\epsilon ,{h^2\over \epsilon })=x\cdot \eta +\phi
_{\rm per}.$$
The \canform{} $\kappa :\, (\psi '_\eta ,\eta )\mapsto (x,\psi
'_x)$ is $(\epsilon +h^2/\epsilon )$-close to the identity and
with $\xi =\xi (\eta ,\epsilon ,h^2/\epsilon )$, we have by
construction
\ekv{3.32}
{
(r_0\circ \kappa )(y,\eta ,\epsilon ,{h^2\over \epsilon
})=\langle r_0(\cdot ,\xi ,\epsilon ,{h^2\over \epsilon
})\rangle =\langle r_0(\cdot ,\eta ,\epsilon ,{h^2\over \epsilon
}\rangle +{\cal O}((\epsilon ,{h^2\over \epsilon })^2), }
which is a \fu{} \indep{} of $y$. Let $U^{-1}$
be the corresponding \fop{} as before. Then after replacing
$\widetilde{P}_\epsilon $ by ${\rm Ad}_U\widetilde{P}_\epsilon $,
we still have (\ref{3.26}), where now $r_0=i\langle q\rangle (\xi )+{\cal
O}(\epsilon +h^2/\epsilon )$ is \indep{} of $x$ and $r_j={\cal
O}(1)$ for $j\ge 1$.

\par We can then make a further conjugation by $e^{{i\over
h}B_1}$ as before, and we get the following conclusion: Assume
that the subprincipal symbol of $P_{\epsilon =0}$ vanishes and
restrict the attention to the range (\ref{3.31}). Then we can find an
elliptic \fop{} $U^{-1}$ of the form (\ref{3.22}) with $\psi $ as above
and $B_1(x_2,hD_x,\epsilon ,h^2/\epsilon ;h)$ with symbol (\ref{3.28}),
and $b_\nu ={\cal O}(1)$, such that 
\ekv{3.33}
{
e^{{i\over h}{\rm ad}_{B_1}}{\rm }Ad_{U} \widetilde{P}_\epsilon
=\widehat{P}_\epsilon (hD_x,\epsilon ,h^2/\epsilon ;h) }
has a symbol $\widehat{P}_\epsilon (\xi ,\epsilon ,h^2/\epsilon
;h)$ of the form (\ref{3.29}), such that (\ref{3.30}) holds. 

\par We finish this section by discussing
what spectral results can be expected from the reductions
above. The first reduction (as in section \ref{section2}) was to conjugate the
original \op{} $P$ by a \fop{} $e^{iG(x,hD,\epsilon )/h}$, with
$G(x,\xi ,\epsilon )\sim \epsilon (G_0(x,\xi )+\epsilon G_1(x,\xi
)+\ldots)$, defined in some complex \neigh{} of $p^{-1}(0)\cap T^*M$,
to achieve that the leading symbol of the conjugated \op{} is of
the form
$p+i\epsilon \langle q\rangle +{\cal O}(\epsilon ^2)$ and Poisson
commutes with $p$. At least formally, the new \op{} also acts on
$L^2(M)$ and we have no Floquet type conditions to worry about.
Geometrically, this corresponds to the fact that a \ctf{}
$\kappa =\exp H_G$ with a single-valued generator $G={\cal
O}(\epsilon )$ preserves actions along closed loops: $\int_{\kappa
\circ \gamma }\xi dx=\int_\gamma \eta dy$, for every closed loop
$\gamma $. 

\par The second reduction was to take $\kappa $ in (\ref{3.4}) and to
conjugate by the inverse of the corresponding \fop{} $U$. Let
$\alpha _1$(=$\gamma _0$) and $\alpha _2$ be the fundamental
cycles in $\Lambda _{0,0}$ given by $\alpha _j=\kappa \circ \beta
_j$, where $\beta _1$, $\beta _2$ are the fundamental cycles in
${\bf T}^2\simeq \{ (x,0)\in T^*{\bf T}^2\} $, given by $x_2=0$
and $x_1=0$ respectively. Put
\ekv{3.34}
{S_j=\int _{\alpha _j} \xi dx,}
so that $S_j$ is the difference of actions, $\int_{\kappa \circ
\beta _j}\xi dx-\int_{\beta _j}\eta dy$, $j=1,2$. Since $\kappa $ is a
\ctf{} we know that if $\beta $ is a closed loop homotopic to
$\beta _j$, then $\int_{\kappa \circ \beta }\xi dx -\int_\beta 
\eta dy=S_j$. 

\par As in \cite{MeSj} or as in \Th{} \ref{Th1.4}, we see (at least
formally) that if we want $Uu$ to be singlevalued on $M$ (possibly
defined only \ml{}ly near $\Lambda _{0,0}$), then $u$ should not
necessarily be periodic on ${\bf R}^2$ (i.e. a \fu{} on ${\bf
T}^2$) but a Floquet periodic \fu{} with 
\ekv{3.35}
{
u(x-\nu )=e^{{i\nu \cdot S\over 2\pi h}+{i\nu \cdot k_0\over 4}}u(x),\
\nu \in (2\pi {\bf Z})^2,\ S=(S_1,S_2),\ k_0\in{\bf Z}^2. }
The conjugated \op{} ${\rm Ad}_{U^{-1}e^{{i\over h}G}}P_\epsilon
$ should therefore act on Floquet periodic \fu{}s as in (\ref{3.35}). 

\par The further conjugations are by \op{}s on the torus that
conserve the property (\ref{3.35}). This is clear from the definitions,
and corresponds to the fact that a \ctf{}: $(y,\eta )\mapsto
(x,\xi )$, generated by $\psi (x,\eta )=x\cdot \eta +\phi _{\rm
per}(x,\eta )$ and close to the identity, conserves actions.
Indeed, on the graph of the \tf{}, we have $\xi dx +yd\eta =d\psi
$, so 
$$\xi dx-\eta dy=d(\psi -y\cdot \eta )=d((x-y)\cdot \eta +\phi
_{\rm per}(x,\eta )),$$
and $(x-y)\cdot \eta +\phi _{\rm per}(x,\eta )$ is singlevalued
on the graph. On the other hand the space of Floquet periodic
\fu{}s as in (\ref{3.35}), equipped with the $L^2$-norm over a
fundamental domain of ${\bf T}^2$, has the ON basis:
\ekv{3.36}
{
e_k(x)=e^{{i\over h}x\cdot (h(k-{k_0\over 4})-{S\over 2\pi })},\
k\in {\bf Z}^2, }
and applying our reductions down to the \op{}
$\widehat{P}_\epsilon $ in the cases (a) and (b) above, we get
formally (in the sense that we do not define the notion of
quasi-\ev{}):
\medskip

\begin{prop}\label{Prop3.1}  Recall that we took $F_0=0$ and
that $S$, $k_0$ are the actions and the Maslov indices in {\rm (\ref{3.34})},
{\rm (\ref{3.35})}.
\par\noindent a) In the general case, $P_\epsilon $ has the
quasi-\ev{}s in $]-{1\over {\cal O}(1)},{1\over {\cal
O}(1)}[+i\epsilon ]-{1\over {\cal O}(1)},{1\over {\cal O}(1)}[$
for $\epsilon ={\cal O}(h^\delta )$, $h/\epsilon \ll 1$:
\ekv{3.37}
{\widehat{P}\left(h(k-{k_0\over 4})-{S\over 2\pi },\epsilon ,{h\over
\epsilon };h\right),\ k\in{\bf Z}^2,}
where $\widehat{P}(\xi ,\epsilon ,{h\over \epsilon };h)$ is
\hol{} in $\xi \in{\rm neigh\,}(0,{\bf C}^2)$, smooth 
in ${h\over \epsilon }, \epsilon \in {\rm neigh\,}(0,{\bf R})$ 
and has the \asy{} expansion {\rm (\ref{3.19})}, when $h\to 0$. 
\smallskip
\par\noindent b) If we assume that $P_{\epsilon =0}$ 
has subprincipal symbol 
$0$, then $P_\epsilon $ has the 
quasi-eigen\-va\-lues in 
$]-{1\over {\cal O}(1)},{1\over {\cal O}(1)}[+i\epsilon ]-{1\over {\cal O}(1)},{1\over {\cal O}(1)}[$
for $\epsilon ={\cal O}(h^\delta )$, $h^2/\epsilon \ll 1$:
\ekv{3.38}
{\widehat{P}\left(h(k-{k_0\over 4})-{S\over 2\pi },\epsilon ,{h^2\over
\epsilon };h\right),\ k\in{\bf Z}^2,}
where $\widehat{P}(\xi ,\epsilon ,h^2/\epsilon ;h)$ is \hol{} in
$\xi \in{\rm neigh\,}(0,{\bf C}^2)$, smooth in $\epsilon$ and $h^2/\epsilon \in {\rm neigh\,}(0,{\bf R})$ 
and has the \asy{} expansion {\rm (\ref{3.29})}, {\rm (\ref{3.30})}, when $h\to 0$.
\end{prop}

\section{Quasi-\ev{}s in the extreme cases}\label{section4}
\setcounter{equation}{0}
\par We make the assumptions of the case II in the introduction
and assume, in order to fix the ideas, that 
\ekv{4.1}
{
0=F_0=\langle \Re q\rangle _{{\rm min},0}.
}
Apply \Pro{} \ref{Prop2.2} and reduce $P_\epsilon $ near
$\gamma _0$ to $\widetilde{P}_\epsilon=
\widetilde{P}(x,hD_{t,x},\epsilon ;h)$ with symbol described in
that \pro{}. Recall that
$\widetilde{P}_\epsilon $ has the leading symbol 
\ekv{4.2}
{\widetilde{p}_\epsilon =f(\tau )+i\epsilon \langle q\rangle
(\tau ,x,\xi )+{\cal O}(\epsilon ^2),}
where $\langle q\rangle (\tau ,x,\xi )$ is equal to the
original averaged function $\langle q\rangle
$, composed with the \ctf{} $\kappa $ of \Pro{} \ref{Prop2.1}. The
assumptions (\ref{0.22}) and (\ref{4.1}) imply that 
\ekv{4.3}
{\Re \langle q\rangle (0,x,\xi )\sim \vert (x,\xi )\vert ^2}
on the real domain. Also recall that we have the assumption
(\ref{0.16}) which with (\ref{4.3}) implies that 
\ekv{4.4}
{
\langle q\rangle (\tau ,x,\xi )=g(\tau ,\Re \langle q\rangle (\tau
,x,\xi )) }
on the real domain, for some \an{} \fu{} $g(\tau ,q)$ with
$g(0,0)=0$, $\Re g(\tau ,q)=q$.

\par We conclude that $(x,\xi )\mapsto i\langle q\rangle (\tau
,x,\xi )+{\cal O}(\epsilon )$, appearing in (\ref{4.2}), has a
\nondeg{} critical point $(x(\tau ,\epsilon ),\xi (\tau ,\epsilon
))={\cal O}(\vert \tau \vert +\epsilon )$ depending \an{}ally on
$\tau ,\epsilon $ and real when $\tau \in{\bf R},\, \epsilon =0$.
After composition with the $(\tau ,\epsilon )$-dependent
(symplectic) translation $(x,\xi )\mapsto (x-x(\tau ,\epsilon ),\xi
-\xi (\tau ,\epsilon ))$ and subtracting the corresponding
critical value, we may assume that the critical point is
$(0,0)$ and hence that 
\ekv{4.5}
{
\widetilde{p}_\epsilon (\tau ,x,\xi )=f(\tau )+i\epsilon q(\tau
,x,\xi ,\epsilon ), }
with 
\ekv{4.6}
{
\Re q(\tau ,x,\xi ,\epsilon )\sim \vert (x,\xi )\vert ^2
}
on the real domain, and 
\ekv{4.7}
{
q(\tau ,x,\xi ,0)=g(\tau ,\Re q(\tau ,x,\xi ,0)),
}
on the real domain, where $g(\tau ,0)=0$, $\Re g(\tau ,q)=q$.

\par We shall next construct a $(\tau ,\epsilon )$-dependent
\ctf{} in the $x,\xi $-variables, which reduces
$\widetilde{p}_\epsilon (\tau ,x,\xi )$ to a \fu{} of $\tau
,\epsilon ,{1\over 2}(x^2+\xi ^2)$. In doing so, we essentially
follow appendix b of \cite{HeSj2}, where the model was $x\xi $ rather
than $p_0:={1\over 2}(x^2+\xi ^2)$. These two quadratic forms are
equivalent up to a constant factor and composition by a linear
complex \ctf{}, so the only difference is that the real domains are
not the same. 

\par Let $p(x,\xi )\sim (x,\xi )^2$ be real and \an{} in a
\neigh{} of (0,0).
\medskip
\begin{lemma}\label{Lemma4.1}. There exists a real and \an{}
\fu{} $f(E)$ defined near $E=0$, with $f(0)=0$, $f'(0)>0$, such
that the Hamilton flow of $f\circ p$ is $2\pi $-periodic, with
$2\pi $ as its minimal period except at $(0,0)$.\end{lemma}
\begin{proof} Consider, first for 
$0<E\ll 1$, the action
$$I(E)=\int_{p^{-1}(E)}\xi dx=E\int_{q_E^{-1}(1)}\eta dy,$$
where $q_E(y,\eta )={1\over E}p(\sqrt{E}(y,\eta ))$, so that
$q_0$ is a positive quadratic form (in the limit $E\to 0$). Then
$q_E$ is an \an{} \fu{} of $\sqrt{E}$ in a \neigh{} of $0$ and
consequently we have the same fact for $I(E)$. If we let $E$
describe a simple closed loop around $0$ in ${\rm neigh\,}(0,{\bf
C})\setminus\{ 0\}$, then
$q_E(y,\eta )$ transforms into $\widetilde{q}_E(y,\eta
)=q_E(-y,-\eta )$ and it follows that $I(E)$ transforms into
itself. It follows that $I(E)$ is \an{} as a \fu{} of $E$. The
period $T(E)$ of the $H_p$-flow is given by $T(E)=I'(E)$ and the
period of the $H_{f\circ p}$-flow is $T(E)/f'(E)$. It suffices to
choose $f$ with $f'(E)=T(E)/2\pi $ and
$f(0)=0$.\end{proof}

\par In the following discussion, we replace $p$ by $f\circ p$,
so that we get a reduction to the case when the $H_p$-flow is
$2\pi $-periodic. After composition with a real linear \ctf{}, we may
assume that $p(x,\xi )=p_0(x,\xi )+{\cal O}((x,\xi )^3)$, even
though that is not really needed for the argument to follow.
Consider the involution $\iota =\exp (\pi H_p)$ with $\iota
^2={\rm id}$. Correspondingly, we have $\iota _0=\exp (\pi
H_{p_0})$, so that $\iota_0(\rho )=-\rho $. Let $f(x,\xi )$ be a
real-valued \an{} \fu{} defined near $(0,0)$ with $df(0,0)\ne 0$,
and put $g={1\over 2}(f-f\circ \iota)$. Then $dg(0)=df(0,0)\ne
0$, and 
\ekv{4.8}
{
g\circ \iota =-g.
}
$\Gamma :=g^{-1}(0)$ is a real curve passing through the origin,
invariant under the action of $\iota$. Let $\Gamma $ also denote
a corresponding complexification. If $g_0,\Gamma _0$
are the corresponding objects for $p_0$, we may assume (though
this is not essential), that $dg(0,0)=dg_0(0,0)$
so that $\Gamma $, $\Gamma _0$ are tangent at $(0,0)$.

\par Since $\Gamma $ is curve, we have ${p_\vert}_{\Gamma
}=q^2$ for some \an{} \fu{} $q$, and similarly
${{p_0}_\vert}_{\Gamma _0}=q_0^2$. (We may assume that
$dq_0=dq\ne 0$ at $0$.) Let $\alpha :\,\Gamma _0\to \Gamma $ be the
\an{} \diffeo{} given by $q\circ \alpha =q_0$, so that
$p\circ \alpha =p_0$ on $\Gamma _0$. For 
\ekv{4.9}
{{\rm neigh\,}((0,0),{\bf C}^2)\ni \rho =\exp tH_{p_0}(\nu ),\
\nu \in\Gamma _0,\, t\in{\bf C},}
we put
\ekv{4.10}
{
\kappa (\rho )=\exp tH_p(\alpha (\nu )).
}
With the precautions taken above, it is easy to see that the
definition of $\kappa (\rho )$ does not depend on how we choose
$\nu \in\Gamma _0$ (unique up to the action of $\iota_0$) and $t$
(unique ${\rm mod\,}(2\pi )$, once $\nu $ has been chosen.) As in
\cite{HeSj2}, we see that some exceptional points $\rho \in{\rm
neigh\,}((0,0),{\bf C}^2)$ cannot be represented as in (\ref{4.9}),
namely the ones $\ne (0,0)$ in the stable outgoing and incoming
complex (Lagrangian) curves for the $iH_{p_0}$-flow, and if $\rho
$ converges to one of these lines, then in general $\vert t\vert
\to \infty $ for the $t$ in (\ref{4.9}), so a priori it is not clear
then that the \rhs{} of (\ref{4.10}) is defined. These difficulties
were analyzed and settled in \cite{HeSj2}, and at this point there is
no difference with our situation, so we conclude that $\kappa $
is a well-defined analytic map in a \neigh{} of $(0,0)$:\medskip

\begin{lemma}\label{Lemma4.2}  With $f,p$ as in Lemma {\rm \ref{Lemma4.1}},
there exists an analytic \ctf{} $\kappa :\,Ê{\rm
neigh\,}((0,0),{\bf R}^2)\to {\rm neigh\,}((0,0),{\bf R}^2)$,
with 
$f\circ p\circ \kappa =p_0$.\end{lemma}

\par If $p$ depends smoothly (\an{}ally) on some real parameters,
and fulfills the assumptions above, then $f,\kappa $ can be
chosen to depend smoothly (\an{}ally) on the same parameters. If
$p=p_\epsilon ={\cal O}((x,\xi )^2)$ is \an{} in $(x,\xi )$,
depends smoothly on
$\epsilon
\in{\rm neigh\,}(0,{\bf R})$ and satisfies the assumptions above
for
$\epsilon =0$, then we get $f_\epsilon (E) ,\kappa _\epsilon
(x,\xi ) $, \hol{} in $E$ and $x,\xi $, depending smoothly on
$\epsilon$ with $f_\epsilon
\circ p_\epsilon \circ \kappa _\epsilon =p_0$, but $f_\epsilon
,\kappa _\epsilon $ are no more necessarily real when $\epsilon
\ne 0$. Clearly $\Im f_\epsilon (E)={\cal O}(\epsilon )$,
$\Im
\kappa _\epsilon (x,\xi )={\cal O}(\epsilon )$ when $E,x,\xi $
are real. In our case the parameters are $\tau ,\epsilon $ and
the above discussion gives:
\medskip
\begin{prop}\label{Prop4.3} For $\widetilde{p}_\epsilon
(\tau ,x,\xi )$ in {\rm (\ref{4.5})}, we can find a \ctf{} $(x,\xi )\mapsto
\kappa _{\tau ,\epsilon }(x,\xi )$ depending analytically on
$\tau $ and smoothly on $\epsilon $ with values in the \hol{}
\ctf{}s: ${\rm neigh\,}((0,0),{\bf C}^2)\to {\rm
neigh\,}((0,0),{\bf C}^2)$, and an \an{} \fu{} $g_\epsilon (\tau
,q)$ depending smoothly on $\epsilon $ such that 
\ekv{4.11}
{
\kappa _{\tau ,\epsilon }(0,0)=(x(\tau ,\epsilon ),\xi (\tau
,\epsilon )), }
\ekv{4.12}
{
\widetilde{p}_\epsilon (\tau ,\kappa _{\tau ,\epsilon }(x,\xi
))=f(\tau )+i\epsilon g_\epsilon (\tau ,{1\over 2}(x^2+\xi ^2)). }
Moreover, $\kappa _{\tau ,0}$ is real when $\tau $ is real and 
\ekv{4.13}
{{\partial \over \partial q}\Re g_\epsilon (0,0)>0.}
\end{prop}

\par As a matter of fact, as in section \ref{section3}, we will apply this
result to a modification of $\widetilde{p}_\epsilon $, containing
also the leading lower order symbol. Before doing so, we recall
how to treat lower order symbols in general for \op{}s with
leading symbol modelled on the 1-dimensional harmonic oscillator
(similarly to what we did in section \ref{section2} and as in \cite{Sj4}). 

\par Consider a formal $h$-\pop{} $Q(x,hD_x;h)$ with symbol 
\ekv{4.14}
{Q(x,\xi ;h)\sim q_0(x,\xi )+hq_1(x,\xi )+...,}
defined in a \neigh{} of $(0,0)\in{\bf R}^2$. As usual,
$q_0,q_1,...$ are supposed to be smooth and we assume 
\ekv{4.15}
{q_0(x,\xi )=g_0(p_0(x,\xi )),}
where $g_0\in C^\infty ({\rm neigh\,}(0,{\bf R}))$ satisfies
$g_0(0)=0$, $g_0'(0)\ne 0$. (We do not assume $g_0$ to be
real-valued.) As in section \ref{section2} we find a smooth function
$a_0(x,\xi )$, defined in a \neigh{} of $(0,0)$, such that 
\ekv{4.16}
{H_{q_0}a_0=q_1-\langle q_1\rangle ,}
where $\langle q_1\rangle $ is the trajectory average ${1\over
2\pi }\int_0^{2\pi }q_1\circ \exp (tH_{p_0})dt$. Adding lower
order corrections, we see that there exists 
\ekv{4.17}
{A(x,\xi ;h)\sim a_0(x,\xi )+ha_1(x,\xi )+\ldots}
with all $a_j$ smooth in some common \neigh{} of $(0,0)$, such
that 
\ekv{4.18}
{e^{iA(x,hD_x;h)}Q(x,hD_x;h)e^{-iA(x,hD_x;h)}=:
\widehat{Q}(x,hD_x;h)}
has a symbol $\widehat{Q}\sim \widehat{q}_0+h\widehat{q}_1+...$,
with $\widehat{q}_0=q_0$ and 
\ekv{4.19}
{H_{q_0}\widehat{q}_j=0,\ \forall j.}
This means that $\widehat{q}_j$ is a smooth \fu{} of $p_0(x,\xi
)$ and as is well-known (and exploited for instance in \cite{Sj4}),
the facts (\ref{4.18}), (\ref{4.19}) can be reformulated by saying that we
have found $A$ as in (\ref{4.17}) such that 
$$e^{iA(x,hD;h)}Q(x,hD;h)e^{-iA(x,hD;h)}=g(p_0(x,hD);h),$$
where $g(E;h)\sim \sum_0^\infty  g_j(E)h^j$ in $C^\infty ({\rm
neigh\,}(0,{\bf R}))$, with $g_0$ as before. When $g_0,q_j$ are
\hol{} in fixed \neigh{}s of $E=0$ and $(x,\xi )=(0,0)$, we get
the corresponding holomorphy for $g_k,\widehat{q}_\ell$.

\par Now return to the \op{} $\widetilde{P}_\epsilon $ of the
beginning of this section. Write the full symbol as 
\ekv{4.20}
{\widetilde{P}_\epsilon (\tau ,x,\xi ,\epsilon
;h)\sim\widetilde{p}_\epsilon (\tau ,x,\xi
)+h\widetilde{p}_1(\tau ,x,\xi ,\epsilon
)+h^2\widetilde{p}_2(\tau ,x,\xi ,\epsilon )+\ldots}
\smallskip
\par\noindent a) Consider first the general case without any
assumptions on the subprincipal symbol, and assume that 
\ekv{4.21}
{h\ll \epsilon <h^\delta ,}
for some fixed $\delta >0$. Following the strategy of section
\ref{section3}, we rewrite (\ref{4.20}) as 
\ekv{4.22}
{\widetilde{P}_\epsilon (\tau ,x,\xi ;h)=f(\tau )+\epsilon
[(i\langle q\rangle (\tau ,x,\xi )+{\cal O}(\epsilon)+
{h\over\epsilon }\widetilde{p}_1(\tau ,x,\xi ))+h{h\over \epsilon
}\widetilde{p}_2+h^2{h\over \epsilon }\widetilde{p}_3+\ldots].}
As before, we now treat $h/\epsilon $ as an additional small
parameter. \Pro{} \ref{Prop4.3} extends to the case when
$\widetilde{p}_\epsilon $ is replaced by $\widetilde{p}_\epsilon
+\epsilon {h\over \epsilon }\widetilde{p}_1$, so we have a \ctf{}
$(x,\xi )\mapsto \kappa _{\tau ,\epsilon ,{h/\epsilon
}}(x,\xi )$ depending \an{}ally on $\tau $ and smoothly on
$\epsilon ,{h\over \epsilon }$, equal to $\kappa _{\tau ,\epsilon
}$ when ${h\over \epsilon }=0$, such that 
$$(\widetilde{p}_\epsilon +\epsilon {h\over \epsilon
}\widetilde{p}_1)(\tau ,\kappa _{\tau ,\epsilon ,{h\over \epsilon
}}(x,\xi ))=f(\tau )+i\epsilon g_{\epsilon ,{h\over \epsilon
}}(\tau ,{1\over 2}(x^2+\xi ^2)),$$
with $g_{\epsilon ,0}=g_\epsilon $ appearing in \Pro{} \ref{Prop4.3}.

\par As in section \ref{section3}, we therefore obtain an elliptic \fop{} 
$U_{\epsilon ,h/\epsilon}$, which is a convolution in $t$, and such that the \F{}
\tf{} \wrt{} $t$, $\widehat{U}_{\epsilon ,h/\epsilon }(\tau )$, is
a 1-dimensional \fop{} in $x$ quantizing $\kappa _{\tau ,\epsilon
,h/\epsilon }$. After conjugation of $\widetilde{P}_\epsilon $
by $U_{\epsilon ,h/\epsilon }$, we get a new \op{}
$\widetilde{P}_\epsilon $ of the same type, with symbol
\ekv{4.23}
{\widetilde{P}(\tau ,x,\xi ,\epsilon ,h/\epsilon ;h)=f(\tau
)+\epsilon [ig_{\epsilon ,{h\over \epsilon }}(\tau ,{1\over
2}(x^2+\xi ^2))+h\widetilde{p}_2+h^2\widetilde{p}_3+\ldots],}
where $\widetilde{p}_2,\widetilde{p}_3,\ldots$
also depend on $h/\epsilon $. 

\par After a further conjugation by $e^{iA(hD_t,x,hD_x,\epsilon
,{h\over \epsilon };h)}$, where each term $A_j$ in the $h$-\asy{}
expansion:
$$A(\tau ,x,\xi ,\epsilon ,{h\over \epsilon };h)\sim A_0(\tau
,x,\xi ,\epsilon ,{h\over \epsilon })+hA_1(\tau
,x,\xi ,\epsilon ,{h\over \epsilon })+\ldots$$
is \hol{} in $\tau ,x,\xi $ in a fixed \neigh{} of
$(0,0,0)\in{\bf C}^3$ and smooth in $\epsilon ,h/\epsilon $, we
get a new \op{} of the form
\ekv{4.24}
{
\widetilde{P}_\epsilon =f(hD_t)+i\epsilon G(hD_t,{1\over
2}(x^2+(hD_x)^2),\epsilon ,{h\over \epsilon };h), }
where 
\ekv{4.25}
{G(\tau ,q,\epsilon ,{h\over \epsilon };h)\sim\sum_0^\infty 
G_j(\tau ,q,\epsilon ,{h\over \epsilon })h^j,}
with $G_j$ \hol{} in $\tau ,q$ in a $j$-\indep{} \neigh{} of
$(0,0)$ and smooth in $\epsilon ,h/\epsilon $. Moreover $G_0$ is
equal to the term $g_{\epsilon ,h/\epsilon }(\tau ,q)$ in
(\ref{4.23}). Recalling that ${1\over 2}(x^2+(hD_x)^2)$ has the \ev{}s
$h({1\over 2}+k_2)$, $k_2\in{\bf N}$, we get the conclusion:
\begin{prop}\label{Prop4.4} Make the assumptions of case {\rm II}
in the introduction, and assume that $F_0=\langle \Re q\rangle
_{{\rm min},0}$ (the case when $F_0$ is a maximum being analogous).
Then in a rectangle $]-{1\over {\cal O}(1)},{1\over {\cal
O}(1)}[+i\epsilon ]F_0-{1\over {\cal O}(1)},F_0+{1\over {\cal
O}(1)}[$, $P_\epsilon $ has the quasi-\ev{}s:
\ekv{4.26}
{f\left(h(k_1-{k_0\over 4})-{S_1\over 2\pi }\right)+
i\epsilon G\left(h(k_1-{k_0\over 4})-{S_1\over
2\pi },h({1\over 2}+k_2),\epsilon ,{h\over \epsilon };h\right),\
(k_1,k_2)\in{\bf Z}\times {\bf N}. }
Here $f(\tau )$ is real-valued with $f(0)=0$, $f'(0)>0$. The function 
$G$ has the properties described in and after {\rm (\ref{4.25})} and $\Re
G_0(0,0,0,0)=F_0$, ${\partial \over \partial q}\Re
G_0(0,0,0,0)>0$. Finally, $k_0$ is a fixed integer.\end{prop}
\par\noindent b) We next consider the case when the subprincipal
symbol of $P_{\epsilon =0}$ vanishes, and assume that 
\ekv{4.27}
{
h^2\ll \epsilon <h^\delta ,
}
for some fixed $\delta >0$. According to the improved Egorov \th{} of section 
\ref{section1}, we know that $\widetilde{p}_1$ in (\ref{4.20})
vanishes for $\epsilon =0$, so we can write ${h\over \epsilon
}\widetilde{p}_1(\tau ,x,\xi ,\epsilon )=h\widehat{p}_1(\tau
,x,\xi ,\epsilon )$ in (\ref{4.22}) and treat this term as a lower
order term, while we now allow ${h^2\over \epsilon
}\widetilde{p}_2$ to be a correction to the leading terms. As in
the corresponding case in section 4, we get $h^2/\epsilon $
as an additional small parameter instead of $h/\epsilon $, and
the same procedure as in case a) now leads to (\ref{4.24}), (\ref{4.25}) with
$h/\epsilon $ replaced by $h^2/\epsilon $.

\begin{prop}\label{Prop4.5} Make the assumptions of \Pro{}
{\rm \ref{Prop4.4}} and assume in addition that the subprincipal symbol
of $P_{\epsilon =0}$ vanishes. Then for $\epsilon $ in the range
{\rm (\ref{4.27})}, $P_\epsilon $ has the quasi-\ev{}s as described in the
preceding \pro{}, with the only difference that "$h/\epsilon $"
in {\rm (\ref{4.26})} should be replaced by "$h^2/\epsilon $".\end{prop}

\section{Global Grushin problem}\label{section5}
\setcounter{equation}{0}

\def\T{{\bf T}}
\def\cantransf{\ctf}
\def\mly{\ml{}ly}
Let $P_{\eps}$ be as in section 1. In sections \ref{section3} 
and \ref{section4} we have constructed microlocal normal forms for 
$P_{\eps}$ near a Lagrangian torus and near a closed $H_p$-trajectory, 
respectively. 
The purpose of this section is to justify the preceding microlocal 
constructions and computations, and to show that the quasi-eigenvalues of 
Proposition \ref{Prop3.1} and Propositions \ref{Prop4.4} and \ref{Prop4.5} 
give, modulo ${\cal O}(h^{\infty})$, all of the true eigenvalues of $P_{\eps}$, 
in suitable regions of the complex plane. This will be achieved by studying an 
auxiliary global Grushin problem, well-posed in a certain $h$-dependent Hilbert space, and 
the first and the main step for us will be to define this space globally. The 
actual setup of the Grushin problem and some of the details of the 
computations will be closely related to the corresponding analysis 
in \cite{MeSj}. 

When constructing the Hilbert space, we shall inspect all the 
steps of the microlocal reductions of sections 
\ref{section2}--\ref{section4}, and implement each 
step of the construction. In doing so, for simplicity, we shall 
concentrate on the case when $M=\real^2$. In view of the results of the 
appendix, 
it will be clear how to extend the following discussion to the case of 
compact real-analytic manifolds. Also, in order to simplify the 
presentation, 
we shall assume throughout the section that the order function $m$, 
introduced in (\ref{0.2}), is equal to $1$. Again, it will 
be 
clear that the discussion below will extend to the case of a general order 
function. Throughout this section we shall assume that $\eps={\cal O}(h^{\delta})$, 
for some fixed $\delta>0$. 

\vskip 2mm
Let $G=G(x,\xi,\eps)$ be as in (\ref{2.6}). We shall introduce an 
IR-manifold $\Lambda_{\eps G}\subset \comp^4$, which in a complex
\neigh{} of $p^{-1}(0)\cap \real^4$ is equal to $\exp(i\eps 
H_G)(\real^4)$, and further away from $p^{-1}(0)\cap \real^4$ agrees with 
the real phase space $\real^4$. The manifold $\Lambda_{\eps G}$ will be $\eps$-close to $\real^4$, 
and when defining it, it will be convenient to work on the FBI transform side. 
We shall use the FBI-Bargmann transform
\begeq
\label{6.1}
Tu(x)=C h^{-3/2} \int e^{i\varphi(x,y)/h}u(y)\,dy, \quad x\in 
\comp^2,\quad C>0, 
\endeq
where $\varphi(x,y)=i/2 (x-y)^2$. Associated to $T$ there is a complex 
linear canonical 
transformation $\kappa_T$, given by 
$$
\comp^4 \ni (y,-\varphi'_y(x,y))\longmapsto (x,\varphi'_x(x,y)) \in 
\comp^4.
$$
It is well known, see \cite{Sj6}, that $\kappa_T$ maps $\real^4$ onto 
$$
\Lambda_{\Phi_0}:=\left\{\left(x,\frac{2}{i}\frac{\partial 
\Phi_0}{\partial x}\right), 
x\in \comp^2\right\},\quad \Phi_0(x)=\frac{(\Im x)^2}{2}.
$$
The IR-manifold $\Lambda_{\eps G}$ has already been defined near $p^{-1}(0)\cap \real^4$, and 
when constructing it globally, we require that the IR-manifold $\kappa_T(\Lambda_{\eps G})$ 
should agree with $\Lambda_{\Phi_0}$ outside a bounded set and that it is 
$\eps$-close to that manifold everywhere. We define therefore $\Lambda_{\eps G}$ so 
that the representation 
\begeq
\label{6.2}
\kappa_T(\Lambda_{\eps G})=
\left\{\left(x,\frac{2}{i}\frac{\partial \Phi}{\partial x}\right), 
x\in \comp^2\right\} =:\Lambda_{\Phi} 
\endeq
holds true. Here the function $\Phi \in C^{\infty}(\comp^2;\real)$ is uniformly 
strictly \plsh{}, and is such that 
$$
\Phi(x)=\Phi_0(x)+\eps g(x,\eps),
$$
with $g(x,\eps)\in C^{\infty}$ in both arguments and with a uniformly 
compact support with respect to $x$. 

Associated to $\Lambda_{\eps G}$ we then introduce the 
corresponding Hilbert space 
$H(\Lambda_{\eps G})$ which agrees with $L^2(\real^2)$ as a space, and 
which we equip with the norm 
$\norm{u}:=\norm{Tu}_{L^2_{\Phi}}$. Here $L^2_{\Phi}=L^2(\comp^2; 
e^{-2\Phi/h}L(dx))$, with $L(dx)$ being the 
Lebesgue measure on $\comp^2$. 

\vskip 2mm
Performing a contour deformation in the integral representation of 
$P_{\eps}$ 
on the FBI-Bargmann transform side, as in \cite{MeSj}, \cite{Sj6}, we see that 
\begeq
\label{6.3}
P_{\eps}={\cal O}(1): H(\Lambda_{\eps G})\rightarrow
H(\Lambda_{\eps G}),
\endeq
and the leading symbol on the FBI transform side is then 
$p_{\eps}\circ \kappa_T^{-1}\bigg|_{\Lambda_{\Phi}}$. Continuing to work 
on the FBI-Bargmann 
transform side, as in section 2 of \cite{MeSj}, we introduce a microlocally 
unitary 
semiclassical \fourior{} 
\begeq
\label{6.4}
e^{\eps G(x,hD_x,\eps)/h}: L^2(\real^2) \rightarrow H(\Lambda_{\eps G}), 
\endeq
microlocally defined near $p^{-1}(0)\cap \real^4$, and associated to the 
complex 
canonical 
transformation $\exp(i\eps H_G): \real^4 \rightarrow \Lambda_{\eps G}$. 
The operator in 
(\ref{6.3}) is then microlocally near $p^{-1}(0)$, unitarily equivalent 
to the operator  
$$
e^{-\eps G(x,hD_x,\eps)/h} P_{\eps} e^{\eps G(x,hD_x,\eps)/h}: L^2 
\rightarrow L^2, 
$$
with the principal symbol  
\begeq
\label{6.5}
p\circ \exp(i\eps H_G)=p+i\eps \langle{q}\rangle +{\cal O}(\eps^2).
\endeq
This averaging procedure allows us therefore to reduce the further analysis 
to an operator $P_{\eps}$, microlocally defined near $p^{-1}(0)\cap \real^4$, which has 
the principal symbol (\ref{6.5}), where $\langle{q}\rangle$, 
as well as the ${\cal O}(\eps^2)$-term, are in involution with $p$. As 
explained in section \ref{section3}, at this stage the operator $P_{\eps}$ acts on 
single-valued functions in $L^2(\real^2)$.

\vskip 4mm
In the first part of this section we shall concentrate on the torus case 
of section \ref{section3}. We assume therefore that $dp$ and $d\Re \langle{q}\rangle$ 
are linearly independent on the set 
\begeq
\label{6.6}
\Lambda_{0,0}:\,\, p=0, \,\, \Re \langle{q}\rangle=0.
\endeq
We recall also the assumption that $T(0)$ is the minimal period of every 
closed $H_p$-trajectory in the Lagrangian torus $\Lambda_{0,0}$, and 
notice 
that in a \neigh{} of $\Lambda_{0,0}$, $p$ and $\Re \langle{q}\rangle$ 
form a completely integrable system. Introduce a new Lagrangian torus 
$\widetilde{\Lambda}_{0,0}\subset \Lambda_{\eps G}$ defined by 
\begeq
\label{6.7}
\widetilde{\Lambda}_{0,0}:\,\, p\circ\exp(-i\eps H_G)=0, 
\,\,\Re\langle{q}\rangle\circ \exp(-i\eps H_G)=0.
\endeq 
In what follows we shall often identify the tori $\Lambda_{0,0}$ and 
$\widetilde{\Lambda}_{0,0}$ by means 
of $\exp(i\eps H_G)$, and we shall continue to write $\Lambda_{0,0}$ for 
$\widetilde{\Lambda}_{0,0}$ when there is no risk of confusion. 
Combi\-ning $\exp(i\eps H_G)$ with the canonical transformation $\kappa$, 
introduced in (\ref{3.4}), and given by the action-angle 
coordinates 
associated with $p$, $\Re \langle{q}\rangle$, we get a smooth canonical 
diffeomorphism
\begeq\label{6.8}
\kappa_{\eps}:{\rm neigh}\left(\xi=0, T^*{\bf T}^2\right) \rightarrow 
{\rm neigh}\left(\Lambda_{0,0},\Lambda_{\eps G}\right),
\endeq
so that $\kappa_{\eps}=\exp(i\eps H_G)\circ\kappa$. As in (\ref{3.34}), we set 
$$
S_j=\int_{\alpha_j} \xi\,dx,\quad j=1,2,
$$
where $\alpha_1$ and $\alpha_2$ are the fundamental cycles in 
$\Lambda_{0,0}$, with $\alpha_1$ 
corresponding to a closed $H_p$-trajectory of the minimal period $T(0)$. 
Introduce also the ``Maslov
indices'' $k_0(\alpha_j)\in \z$, $j=1,2$, of the cycles $\alpha_j$, 
defined as in Proposition \ref{Prop1.3}. Let 
$L^2_{\theta}({\bf T}^2)$ be the subspace of $L^2_{\rm{loc}}(\real^2)$ 
consisting of Floquet periodic functions 
$u(x)$, satisfying 
\begeq
\label{6.9}  
u(x-\nu)=e^{i\theta \cdot \nu} u(x), \quad \nu \in \left(2\pi 
\z\right)^2,\quad \wrtext{where} \,\,\,\theta=\frac{S}{2\pi 
h}+\frac{k_0}{4}. 
\endeq
Here $S=(S_1,S_2)$ and $k_0=(k_0(\alpha_1), k_0(\alpha_2))\in \z^2$. An 
application of Theorem \ref{Th1.4} 
allows us to conclude that there exists a microlocally unitary multivalued 
\fourior{} 
\begeq
\label{6.10}
U: L^2_{\theta}({\bf T}^2)\rightarrow L^2(\real^2),
\endeq
microlocally defined from a \neigh{} of $\xi=0$ in $T^*{\bf T}^2$ to a \neigh{} 
of 
$\Lambda_{0,0}$ in $\real^4$, and associated to $\kappa$ in (\ref{3.4}). Moreover, 
$U$ satisfies the improved Egorov property (\ref{1.3}). The composition 
$e^{\eps G(x,hD_x,\eps)/h}\circ U$ is then 
associated with $\kappa_{\eps}$ in (\ref{6.8}), and 
we have a Egorov's theorem, still with the improved property (\ref{1.3}). 
The operator $P_{\eps}$, acting in $H(\Lambda_{\eps G})$ 
is therefore unitarily equivalent to an 
$h$-\pseudor{} microlocally defined near $\xi=0$, acting in 
$L^2_{\theta}({\bf T}^2)$, and which has the leading symbol 
$$
p(\xi_1)+i\eps \langle{q}\rangle(\xi)+{\cal O}(\eps^2),
$$
independent of $x_1$. We shall continue to write $P_{\eps}$ for the 
conjugated operator on ${\bf T}^2$. 

\vskip 2mm
\noindent
From section \ref{section3} we next recall that there exists an elliptic \pseudor{} of 
the form $e^{iA/h}$, acting on $L^2_{\theta}({\bf T}^2)$, such 
that after a conjugation by it, 
the full symbol of ${P}_{\eps}$ becomes independent of $x_1$. Recall also 
that $A$ is constructed as a formal power series in 
$\eps$ and $h$, with coefficients holomorphic in a fixed complex \neigh{} 
of the zero section 
of $T^*{\bf T}^2$. These formal power series are then realized as 
$C^{\infty}$-symbols, in view of our basic assumption 
$\eps = {\cal O}(h^{\delta})$, $\delta>0$. 

\vskip 2mm
Summing up the discussion so far, we have now achieved that, microlocally 
near 
$\Lambda_{0,0}$, the operator 
$$
P_{\eps}: H(\Lambda_{\eps G}) \rightarrow H(\Lambda_{\eps G})
$$
is equivalent to an operator of the form 
\begeq
\label{6.11}
\widetilde{P}_{\eps}(x_2,\xi,\eps;h) \sim \sum_{\nu=0}^{\infty} h^{\nu} 
\widetilde{p}_{\nu}(x_2,\xi,\eps)
\endeq
acting on $L^2_{\theta}({\bf T}^2)$. Here $\widetilde{p}_{\nu}(x_2,\xi,\eps)$ 
are holomorphic in a $\nu$-independent complex \neigh{} of $\xi=0$, and 
$$
\widetilde{p}_0=p(\xi_1)+i\eps \langle{q}\rangle(\xi)+{\cal O}(\eps^2).
$$
Furthermore, $\widetilde{P}_{\eps=0}$ is selfadjoint. 

\Remark. It follows from the construction together with Theorem \ref{Th1.4} 
that if the subprincipal symbol of $P_{\eps=0}$ vanishes, 
then $\widetilde{p}_1(x_2,\xi,0)=0$. 

\vskip 2mm
We must now implement the final conju\-gation of $\widetilde{P}_{\eps}$, 
which remo\-ves the 
$x_2$-depen\-dence in the full symbol. In doing so, we shall first assume 
that we are in the general case, so that the subprincipal symbol of 
$P_{\eps=0}$ does not necessarily vanish. We shall work under the 
assumption  
\begeq
\label{6.12}
\frac{h}{\eps}\leq \delta_0 \ll 1.
\endeq 
As in section \ref{section3}, we write 
$$
\widetilde{P}_{\eps}=p(\xi_1)+\eps 
\left(r_0\left(x_2,\xi,\eps,\frac{h}{\eps}\right)+
h r_1\left(x_2,\xi,\eps,\frac{h}{\eps}\right)+\ldots\right), 
$$
where
$$
r_0\left(x_2,\xi,\eps,\frac{h}{\eps}\right)=i\langle{q}\rangle(\xi)+{\cal 
O}(\eps)+
\frac{h}{\eps}\widetilde{p}_1(x_2,\xi,\eps), 
$$
and $r_j={\cal O}_j(h/\eps)$, $j\geq 1$. Let us introduce a 
complexification of the 
standard 2-torus, $\widetilde{{\bf T}^2}={\bf T}^2+i\real^2$. From the constructions 
of section \ref{section3} we 
know that there exists a holomorphic \cantransf{}  
\begin{eqnarray}
\label{6.13}
& & \widetilde{\kappa}: {\rm neigh}\left(\Im y=\eta=0, \widetilde{{\bf T}^2}\times 
\comp^2\right) \\ \nonumber
& & \ni (y,\eta) \mapsto (x,\xi) \in {\rm neigh}\left(\Im x=\xi=0, \widetilde{{\bf T}^2}\times \comp^2\right)
\end{eqnarray}
with the generating function of the form 
\begeq
\label{6.14}
\psi\left(x,\eta,\eps,\frac{h}{\eps}\right)=x\cdot\eta+ 
\phi_{\rm per}\left(x_2,\eta,\eps,\frac{h}{\eps}\right), \quad \phi_{\rm 
per}=
{\cal O}\left(\eps+\frac{h}{\eps}\right),
\endeq 
and such that 
\begeq
\label{6.15}
\left(r_0\circ \widetilde{\kappa}\right)(y,\eta,\eps, 
h/\eps)=\langle{r_0(\cdot,\eta,\eps,h/\eps)}\rangle+
{\cal O}\left(\left(\eps,\frac{h}{\eps}\right)^2\right)
\endeq
is independent of $y$ --- see (\ref{3.21}). It follows from 
(\ref{6.14}) that $\widetilde{\kappa}$ is $(\eps+h/\eps)$-close to the 
identity, and has the expression  
$$
(y_1,\eta_1; y_2,\eta_2)\longmapsto (x_1(y_2,\eta), \eta_1; x_2(y_2,\eta), 
\xi_2(y_2,\eta)).
$$
In particular it is true that 
$$
\Im x={\cal O}\left(\eps+\frac{h}{\eps}\right), \quad 
\Im \xi_2={\cal O}\left(\eps+\frac{h}{\eps}\right), \quad \Im \xi_1=0,
$$
on the image of $T^*\T^2$. We introduce now an IR-manifold 
$\widetilde{\Lambda} \subset 
\widetilde{\T^2}\times \comp^2$, which is equal to 
$\widetilde{\kappa}\left(T^*\T^2\right)$ in a complex \neigh{} of the zero 
section of $T^*\T^2$, and 
outside another complex fixed \neigh{} of $\xi=0$, coincides with 
$T^*\T^2$. In the intermediate 
region, we shall construct $\widetilde{\Lambda}$ in such a way that it remains 
an 
$\left(\eps+h/{\eps}\right)$-perturbation of $T^*\T^2$, and such that 
everywhere on 
$\widetilde{\Lambda}$ we have the property 
\begeq
\label{crucial}
(x_1,\xi_1; x_2,\xi_2) \in \widetilde{\Lambda} \Longrightarrow \Im \xi_1=0.
\endeq
When constructing $\widetilde{\Lambda}$ and describing the conjugation of 
$\widetilde{P}_{\eps}$ by a \fourior{} associated to 
$\widetilde{\kappa}$,
it is convenient to work on the FBI transform side. As in section 3 of 
\cite{MeSj}, we notice that 
the FBI-Bargmann transformation introduced in (\ref{6.1}) generates an 
operator from 
$L^2_{\theta}(\T^2)$ to the space of Floquet periodic holomorphic 
functions on $\comp^2$. We continue to 
denote this operator by $T$. Then after the application of the canonical 
transformation $\kappa_T$, associated to 
$T$, the cotangent space $T^*\T^2$ becomes an IR-manifold 
$\Lambda_{\Phi_1} \subset \widetilde{\T^2}\times \comp^2$ given by 
$$
\Lambda_{\Phi_1}: \xi=\frac{2}{i} \frac{\partial \Phi_1}{\partial x}=-\Im 
x, 
\quad \Phi_1(x)=\frac{(\Im x)^2}{2}.
$$
Since $T$ is a convolution operator acting separately in $y_1$ and $y_2$, 
we see that 
$$
\kappa_T(\widetilde{\Lambda})=\Lambda_{\Phi}, \quad \Lambda_{\Phi}: 
\xi=\frac{2}{i} \frac{\partial \Phi}{\partial x},  
$$
where $\Phi$ is an $(\eps+h/\eps)$-perturbation of $\Phi_1$ with the 
property that $\xi_1=(2/i)\frac{\partial \Phi}{\partial x_1}$ is real. It 
follows that $\Phi=\Phi(\Im x_1,x_2)$ is independent of $\Re x_1$. Using a standard 
cutoff function around $\Im x=0$, we modify $\Phi$ away from $\Im x=0$ to 
obtain a strictly plurisubharmonic function $\Phi$ which coincides with $\Phi_1$ 
further away from $\Im x=0$, in such a way that $\Phi$ remains an 
$(\eps+h/\eps)$-perturbation of $\Phi_1$ and is still a function 
independent of $\Re x_1$. We then define the global IR-manifold 
$\widetilde{\Lambda}=\kappa_T^{-1}(\Lambda_{\Phi})$. 

Associated to $\widetilde{\kappa}$, there is a \fourior{} $U^{-1}$ introduced in 
(\ref{3.22}), 
$$
U^{-1}={\cal O}(1): L^2(\T^2) \rightarrow H(\widetilde{\Lambda}),
$$
such that the action of $\widetilde{P}_{\eps}$ on $H(\widetilde{\Lambda})$ 
is microlocally near $\xi=0$ unitarily equivalent to the operator 
$$
U \widetilde{P}_{\eps} U^{-1}: L^2(\T^2)\rightarrow L^2(\T^2),
$$
whose Weyl symbol has the form 
\begeq
\label{6.17}
p(\xi_1)+\eps\left( r_0\left(\xi,\eps,\frac{h}{\eps}\right)+h 
r_1\left(x_2,\xi,\eps,\frac{h}{\eps}\right)+\ldots\right).
\endeq
Here 
$$
r_0=i\langle{q}\rangle(\xi)+{\cal O}\left(\eps+\frac{h}{\eps}\right)
$$
is independent of $x$, and 
$$
r_j={\cal O}\left(\eps +\frac{h}{\eps}\right), \quad j \geq 1. 
$$
The corresponding statement is also true when considering the action on 
$L^2_{\theta}(\T^2)$, since $U^{-1}$ preserves the Floquet property 
(\ref{6.9}).

\vskip 2mm
\noindent
Associated to the IR-deformation $\widetilde{\Lambda}$ on the torus side, 
there is an IR-manifold $\widehat{\Lambda}_{\eps} \subset \comp^4$ which is an 
$(\eps+h/\eps)$-perturbation of $\Lambda_{\eps G}$ near $\Lambda_{0,0}$, 
obtained by replacing $\exp(i\eps H_G)\circ \kappa (T^*\T^2)$ 
there by 
$$
\exp(i\eps H_G)\circ \kappa \circ \widetilde{\kappa}(T^*\T^2)=
\exp(i\eps H_G)\circ \kappa (\widetilde{\Lambda}).
$$ 
In such a way we get a globally defined IR-manifold 
$\widehat{\Lambda}_{\eps}$, which is $(\eps+h/\eps)$-close to 
$\Lambda_{\eps G}$ and agrees with $\real^4$ outside a 
\neigh{} of $p^{-1}(0)\cap \real^4$. Associated with $\widehat{\Lambda}_{\eps}$ we then have a 
Hilbert space $H(\widehat{\Lambda}_{\eps})$, defined similarly to 
$H(\Lambda_{\eps G})$, and 
obtained by modifying the standard weight $\Phi_0(x)$ on the FBI-Bargmann 
transform side. We also get a corres\-pon\-ding new Lagrangian torus 
$\widehat{\Lambda}_{0,0}\subset \widehat{\Lambda}_{\eps}$, with 
the property that \mly{} near $\widehat{\Lambda}_{0,0}$, the original 
operator  
$$  
P_{\eps}: H(\widehat{\Lambda}_{\eps}) \rightarrow
H(\widehat{\Lambda}_{\eps})
$$
is equivalent to an operator on $L^2_{\theta}(\T^2)$, whose complete 
symbol has the form (\ref{6.17}). Taking into account the conjugation by 
an 
elliptic operator $e^{iB_1/h}$ on the torus side, which was constructed 
in section \ref{section3} and which eliminates the $x_2$-dependence also in the terms 
$r_j$ 
with $j\geq 1$, we get the following result. 

\begin{prop}\label{Prop6.1}
We make all the assumptions of case {\rm I} in the introduction, and recall that 
we also take $F_0=0$. Assume that 
$\eps ={\cal O}(h^{\delta})$, $\delta>0$ is such that $h/\eps \leq 
\delta_0$, $0<\delta_0\ll 1$. There exists an 
IR-manifold $\widehat{\Lambda}_{\eps}\subset \comp^4$, and a smooth 
Lagrangian torus 
$\widehat{\Lambda}_{0,0}\subset \widehat{\Lambda}_{\eps}$, such that when 
$\rho \in \widehat{\Lambda}_{\eps}$ is away from a small 
\neigh{} of $\widehat{\Lambda}_{0,0}$ in $\widehat{\Lambda}_{\eps}$, we 
have 
\begeq
\label{mlocal}
\abs{\Re P_{\eps}(\rho,h)} \geq \frac{1}{{\cal O}(1)} \quad \wrtext{or} 
\quad \abs{\Im P_{\eps}(\rho,h)}\geq \frac{\eps}{{\cal O}(1)}.
\endeq
The manifold $\widehat{\Lambda}_{\eps}$ is an 
$\left(\eps+\frac{h}{\eps}\right)$-perturbation of $\real^4$ in the 
natural sense, and it is
equal to $\real^4$ outside a \neigh{} of $p^{-1}(0)\cap \real^4$. We have 
$$
P_{\eps} = {\cal O}(1): H(\widehat{\Lambda}_{\eps}) \rightarrow 
H(\widehat{\Lambda}_{\eps}).
$$
There exists a smooth canonical transformation 
$$
\kappa_{\eps}: {\rm neigh}\,(\widehat{\Lambda}_{0,0}, 
\widehat{\Lambda}_{\eps})\rightarrow {\rm neigh}\,(\xi=0, T^*\T^2), 
$$
such that $\kappa_{\eps}(\widehat{\Lambda}_{0,0})= \T^2 \times \{0\}$. 
Associated to 
$\kappa_{\eps}$, there is a \fourior{} 
$$
U={\cal O}(1): H(\widehat{\Lambda}_{\eps})\rightarrow L_{\theta}^2(\T^2),
$$
which has the following properties: 
\begin{enumerate} 
\item $U$ is concentrated to the graph of $\kappa_{\eps}$ in the sense 
that if $\chi_1\in C_0^{\infty}(\widehat{\Lambda}_{\eps})$, 
$\chi_2\in C_0^{\infty}(T^*\T^2)$, are such that 
$$
\left(\supp \chi_2 \times \supp \chi_1\right)\cap 
\overline{\{(\kappa_{\eps}(y,\eta),y,\eta); 
(y,\eta) \in {\rm neigh}(\widehat{\Lambda}_{0,0}, 
\widehat{\Lambda}_{\eps})\}} = \varnothing, 
$$
then
$$
\chi_2(x,hD_x)\circ U\circ \chi_1(x,hD_x)={\cal O}(h^{\infty}): 
H(\widehat{\Lambda}_{\eps}) \rightarrow
L^2_{\theta}(\T^2).
$$ 
\item The operator $U$ is \mly{} invertible: there exists an operator 
$V = {\cal O}(1): L^2_{\theta}(\T^2) \rightarrow 
H(\widehat{\Lambda}_{\eps})$ such that for every 
$\chi_1 \in C_0^{\infty}({\rm neigh}(\widehat{\Lambda}_{0,0}, 
\widehat{\Lambda}_{\eps}))$, we have 
$$
\left(V U - 1\right)\chi_1(x,hD_x)={\cal O}(h^{\infty}): 
H(\widehat{\Lambda}_{\eps})\rightarrow H(\widehat{\Lambda}_{\eps}).
$$
For every $\chi_2 \in C_0^{\infty}({\rm neigh}(\xi=0, T^*\T^2))$, we have 
$$
\left(U V - 1\right)\chi_2(x,hD_x)={\cal O}(h^{\infty}): 
L^2_{\theta}(\T^2)\rightarrow L^2_{\theta}(\T^2).
$$
\item We have Egorov's theorem: Acting on $L^2_{\theta}(\T^2)$, there 
exists $\widehat{P}\left(hD_x,\eps,\frac{h}{\eps};h\right)$ with the 
symbol 
$$
\widehat{P}\left(\xi,\eps,\frac{h}{\eps};h\right) \sim 
p(\xi_1)+\eps\sum_{j=0}^{\infty} h^j 
r_j\left(\xi,\eps,\frac{h}{\eps}\right), 
\quad \abs{\xi}\leq \frac{1}{{\cal O}(1)}, 
$$
with 
$$
r_0=i\langle{q}\rangle(\xi)+{\cal O}\left(\eps+\frac{h}{\eps}\right), 
$$
and 
$$
r_j={\cal O}_j\left(\eps+\frac{h}{\eps}\right),\quad j\geq 1,
$$
such that $\widehat{P}U = U P_{\eps}$ microlocally, i.e. 
$$
\left(\widehat{P}U - U P_{\eps}\right)\chi_1(x,hD_x)={\cal O}(h^{\infty}), 
\quad \chi_2(x,hD_x)\left(\widehat{P}U - U P_{\eps}\right) = {\cal 
O}(h^{\infty}), 
$$ 
for every $\chi_1$, $\chi_2$ as in 2). 
\end{enumerate}
\end{prop} 

\Remark. The estimate (\ref{mlocal}) holds true thanks to the property 
(\ref{crucial}) of the final deformation, since then the term 
$p(\xi_1)$ does not contribute to the imaginary part of the symbol on the 
torus side. The bound (\ref{mlocal}) will allow us to reduce 
the spectral analysis of $P_{\eps}$ to a small \neigh{} of the Lagrangian 
torus $\widehat{\Lambda}_{0,0}$. 

\vskip 2mm
\noindent
Using Proposition \ref{Prop6.1}, we shall now proceed to describe the spectrum of 
$P_{\eps}$ in a rectangle of the form 
\begeq
\label{6.19}
R_{C,\eps}=\biggl \{z\in \comp; \,\, \abs{\Re z}<\frac{1}{C},\quad \abs{\Im z}< \frac{\eps}{C} \biggr \},
\endeq
for a sufficiently large constant $C>0$. We shall show that the 
eigenvalues in (\ref{6.19}) are given by the 
quasi-eigenvalues of Proposition \ref{Prop3.1}, modulo ${\cal O}(h^{\infty})$. In 
doing so, let us consider the set of the 
quasi-eigenvalues, introduced in (\ref{3.37}), 
$$
\Sigma(\eps,h)=\biggl 
\{\widehat{P}\left(h(k-\theta),\eps,\frac{h}{\eps};h\right); \, k\in \z^2 
\biggr \}\bigcap R_{C,\eps},\quad \theta=\frac{S}{2\pi h}+\frac{k_0}{4}. 
$$
Then the distance between 2 elements of $\Sigma(\eps,h)$ corresponding to 
$k,l\in \z^2$, $k \neq l$, is 
$\geq \eps h\abs{k-l}/{\cal O}(1)$. Introduce
$$
\delta:=1/4 \inf_{k\neq l} 
\dist(\widehat{P}(h(k-\theta),\eps,\frac{h}{\eps};h), 
\widehat{P}(h(l-\theta),\eps,\frac{h}{\eps};h))>0, 
$$ 
and consider the family of open discs 
$$
\Omega_k(h):=\biggl \{z\in R_{C,\eps}; 
\abs{z-\widehat{P}(h(k-\theta),\eps,\frac{h}{\eps};h)}< \delta \biggr 
\},\quad k\in \z^2. 
$$
The sets $\Omega_k(h)$ are then disjoint, and $\dist\left(\Omega_k(h), 
\Omega_l(h)\right) \geq \eps h\abs{k-l}/{{\cal O}(1)}$. As a 
warm-up exercise, we shall first show that $\Spec(P_{\eps})$ in the set 
(\ref{6.19}) is contained in the union of the 
$\Omega_k(h)$. 

\vskip 4mm
When $z\in \comp$ is in the rectangle (\ref{6.19}), let us consider the equation 
\begeq
\label{6.20}
\left(P_{\eps}-z\right)u=v,\quad u\in H(\widehat{\Lambda}_{\eps}). 
\endeq
We notice here that the symbol of 
$$
\Im P_{\eps}=\frac{P_{\eps}-P_{\eps}^*}{2i},
$$
taken in the operator sense in $H(\widehat{\Lambda}_{\eps})$, is 
${\cal O}(\eps)$, and from Proposition \ref{Prop6.1}
we know that away from any fixed \neigh{} of $\widehat{\Lambda}_{0,0}$ in 
$\widehat{\Lambda}_{\eps}$ it is true that 
$\abs{\Im P_{\eps}(\rho,h)}>\eps/C$, provided that 
$\abs{\Re P_{\eps}(\rho,h)}\leq 1/C$, where $C>0$ is sufficiently large. 
Here we are using the same letters for the operators 
and the corresponding (Weyl) symbols, and 
$$
\Re P_{\eps}=\frac{P_{\eps}+P_{\eps}^*}{2}: 
H(\widehat{\Lambda}_{\eps})\rightarrow H(\widehat{\Lambda}_{\eps}). 
$$
We shall also write $p$ to denote the leading symbol of $P_{\eps=0}$, 
acting on $H(\widehat{\Lambda}_{\eps})$. 

Let us introduce a smooth partition of unity on the 
manifold $\widehat{\Lambda}_{\eps}$, 
$$
1=\chi+\psi_{1,+}+\psi_{1,-}+\psi_{2,+}+\psi_{2,-}.
$$
Here $\chi \in C_0^{\infty}(\widehat{\Lambda}_{\eps})$ is such that $\chi=1$ 
near $\widehat{\Lambda}_{0,0}$, and $\supp \, \chi$ is contained in a small 
\neigh{} of $\widehat{\Lambda}_{0,0}$ where $U P_{\eps}=\widehat{P}U$. The functions 
$\psi_{1,\pm} \in C^{\infty}_{0}(\widehat{\Lambda}_{\eps})$
are supported in regions, invariant under the $H_p$--flow, where 
$\pm \Im P_{\eps}>\eps /C$, respectively. Finally $\psi_{2,\pm}\in C^{\infty}_{\rm b}(\widehat{\Lambda}_{\eps})$ 
are such that 
$$
\supp \psi_{2,\pm} \subset \biggl\{ \rho; \pm \Re P_{\eps}(\rho,h) > 1/C\biggr \}. 
$$
Moreover, we arrange so that the functions $\psi_{1,\pm}$ Poisson commute 
with $p$ on $\widehat{\Lambda}_{\eps}$. We shall prove that 
\begeq
\label{6.21}
\norm{(1-\chi)u}\leq {\cal O}\left(\frac{1}{\eps}\right)\norm{v}+ {\cal O}(h^{\infty})\norm{u},
\endeq
where we let $\norm{\cdot}$ stand for the norm in $H(\widehat{\Lambda}_{\eps})$. In doing so, 
we shall first derive a priori estimates for $\psi_{1,+}u$. 

When $N\in \nat$, let 
$$  
\psi_0 \prec \psi_1 \prec \ldots \prec \psi_N,\quad \psi_0:=\psi_{1,+},
$$
be cutoff functions in $C^{\infty}_0(\widehat{\Lambda}_{\eps};[0,1])$, 
supported in an $H_p$--flow invariant region where $\Im P_{\eps}\sim \eps$, and which are 
in involution with $p$. Here standard notation $f\prec g$ means that $\supp\, f$ is contained 
in the interior of the set where $g=1$. It is then true that in the operator norm, 
\begeq
\label{6.21.1}
[P_{\eps},\psi_j]=[P_{\eps=0},\psi_j]+{\cal O}(\eps h)={\cal O}(h^2)+{\cal O}(\eps h)={\cal O}(\eps h), 
\quad 0\leq j\leq N, 
\endeq
since $\eps \geq h$. For future reference we notice that in the case when the 
subprincipal symbol of $P_{\eps=0}$ vanishes, the Weyl calculus shows that $[P_{\eps=0},\psi_j]={\cal O}(h^3)$, 
and since $\eps \geq h^2$, we still get (\ref{6.21.1}). Here we have also used that the subprincipal symbol of 
$\psi_j$ is $0$, $0\leq j\leq N$. 

Near the support of $\psi_j$ it is true that $\Im P_{\eps}\sim \eps$, and an application 
of the semiclassical G\aa{}rding inequality allows us therefore to conclude that 
$$
\left(\Im (P_{\eps}-z)\psi_j u|\psi_j u\right)
\geq \frac{\eps}{{\cal O}(1)}\norm{\psi_j u}^2 -{\cal O}(h^{\infty}) \norm{u}^2.
$$
Here the inner product is taken in $H(\widehat{\Lambda}_{\eps})$. On the other hand, we have 
$$
\left(\Im (P_{\eps}-z)\psi_j u|\psi_j u\right) =
\Im \biggl( (\psi_j (P_{\eps}-z)u|\psi_j u)+([P_{\eps},\psi_j]u|\psi_j u)\biggr),
$$
and since in the operator sense $\psi_j(1-\psi_{j+1})={\cal O}(h^{\infty})$, we see that 
the absolute value of this expression does not exceed 
$$
{\cal O}(1)\norm{(P_{\eps}-z)u}\,\norm{\psi_j u}+{\cal O}(\eps h) \norm{\psi_{j+1}u}^2+{\cal O}(h^{\infty})\norm{u}^2. 
$$
We get 
\begin{eqnarray*}
& & \frac{\eps}{C}\norm{\psi_j u}^2 \leq {\cal O}(1)\norm{(P_{\eps}-z)u}\,\norm{\psi_j u}+
{\cal O}(\eps h) \norm{\psi_{j+1}u}^2+{\cal O}(h^{\infty})\norm{u}^2 \\
& \leq & \frac{\eps}{2C}\norm{\psi_j u}^2 +
\frac{{\cal O}(1)}{\eps}\norm{(P_{\eps}-z)u}^2+{\cal O}(\eps h) \norm{\psi_{j+1}u}^2+{\cal O}(h^{\infty})\norm{u}^2,
\end{eqnarray*} 
and hence, 
$$
\norm{\psi_j u}^2 \leq \frac{{\cal O}(1)}{\eps^2}\norm{(P_{\eps}-z)u}^2 + {\cal O}(h)\norm{\psi_{j+1} u}^2+
{\cal O}(h^{\infty})\norm{u}^2.
$$
Combining these estimates for $j=0,1,\ldots N$, we get
$$
\norm{\psi_0 u}^2 \leq \frac{{\cal O}(1)}{\eps^2 }\norm{(P_{\eps}-z)u}^2+{\cal O}_N(1) h^N \norm{\psi_N u}^2+
{\cal O}(h^{\infty})\norm{u}^2,
$$
and therefore
$$
\norm{\psi_{1,+}u}\leq \frac{{\cal O}(1)}{\eps}\norm{v}+{\cal O}(h^{\infty})\norm{u}. 
$$
The same estimate can be obtained for $\psi_{1,-}u$, microlocally concentrated in a flow invariant 
region where $\Im P_{\eps}\sim -\eps$, and a fortiori such estimates also hold in regions where 
$\Re P_{\eps}\sim 1$ and $\Re P_{\eps}\sim -1$. The bound (\ref{6.21}) follows. 

\vskip 4mm
Write next
\begeq
\label{6.22}
\left(P_{\eps}-z\right)\chi u=\chi v+w,\,\,w=[P_{\eps},\chi]u,
\endeq
where $w$ satisfies 
$$
\norm{w}\leq {\cal O}\left(\frac{h}{\eps}\right)\norm{v}+{\cal O}(h^{\infty})\norm{u}. 
$$
Here we have used (\ref{6.21}) with a cutoff closer to $\widehat{\Lambda}_{0,0}$. 
Applying the operator $U$ of Proposition \ref{Prop6.1} to (\ref{6.22}), we get 
$$
\left(\widehat{P}-z\right)U \chi u=U \chi v +Uw+T_{\infty}u,
$$
where 
$$
T_{\infty}={\cal O}(h^{\infty}): H(\widehat{\Lambda}_{\eps})\rightarrow L^2_{\theta}(\T^2). 
$$
Using an expansion in Fourier series (\ref{6.24}) below, we see that the operator 
$\widehat{P}-z: L^2_{\theta}(\T^2) \rightarrow L^2_{\theta}(\T^2)$ is 
invertible, microlocally in $\abs{\xi} \leq 1/{\cal O}(1)$, with a microlocal inverse of the norm 
${\cal O}(1/ \eps h)$, provided that $z\in R_{C,\eps}$ avoids the discs $\Omega_k(h)$. 
Using also the uniform boundedness of the microlocal inverse $V$ of $U$, we get   
\begeq
\label{6.23}
\norm{\chi u}\leq \frac{{\cal O}(1)}{\eps h} \norm{v}+
{\cal O}(h^{\infty})\norm{u}.
\endeq
Combining (\ref{6.21}) and (\ref{6.23}), we see that when $z\in R_{C,\eps}$ 
is in the 
complement of the union of the $\Omega_k(h)$, the operator 
$$
P_{\eps}-z: H(\widehat{\Lambda}_{\eps})\rightarrow 
H(\widehat{\Lambda}_{\eps})
$$
is injective. Since the ellipticity assumption (\ref{0.6}) 
implies that it is a Fredholm operator of index zero, we know that 
$P_{\eps}-z: H(\widehat{\Lambda}_{\eps})\rightarrow H(\widehat{\Lambda}_{\eps})$ is bijective. 

\vskip 2mm
We shall now let $z$ vary in the disc $\Omega_{k}(h)\subset R_{C,\eps}$, for 
some $k\in \z^2$. We shall show that 
$z\in \Omega_k(h)$ is an eigenvalue of $P_{\eps}$ if and only if 
$z=\widehat{P}(h(k-\theta),\eps,\frac{h}{\eps};h)+r$, where 
$r={\cal O}(h^{\infty})$. In doing so, we shall study a globally 
well-posed Grushin problem for the operator 
$P_{\eps}-z$ in the space $H(\widehat{\Lambda}_{\eps})$. 

As a preparation for that, we shall introduce an auxiliary Grushin problem 
for the operator 
$\widehat{P}-z$, defined microlocally near $\xi=0$ in $T^*\T^2$. From 
(\ref{3.36}), let us recall the 
functions 
$$
e_l(x)=\frac{1}{2\pi} e^{i(l-\theta)x}=\frac{1}{2\pi} 
e^{\frac{i}{h}\left(h\left(l-\frac{k_0}{4}\right)-\frac{S}{2\pi}\right)x},
$$
which form an ON basis for the space $L^2_{\theta}(\T^2)$, so that when 
$u\in L^2_{\theta}(\T^2)$, we have a Fourier series expansion, 
\begeq
\label{6.24}
u(x)= \sum_{l\in {\rm {\bf Z}^2}} \widehat{u}(l-\theta) e_l(x).
\endeq
We also remark that $e_l(x)$ are microlocally concentrated to the region 
of the phase space where
$\xi \sim h\left(l-\frac{k_0}{4}\right)-S/2\pi$.

Introduce rank one operators $\widehat{R}_+: L^2_{\theta}(\T^2)\rightarrow 
\comp$ and 
$\widehat{R}_- : \comp \rightarrow L^2_{\theta}(\T^2)$, given by 
$\widehat{R}_+ u=(u|e_{k})$ and $\widehat{R}_-u_-=u_- e_{k}$. 
Here the inner product in the definition of $\widehat{R}_+$ is taken in 
the space $L^2_{\theta}(\T^2)$. Using (\ref{6.24}), 
it is then easy to see that the operator 
\begeq
\label{6.25}
\widehat{{\cal P}}:= 
\left( \begin{array}{ccc}
\widehat{P}-z & \widehat{R}_- \\
\widehat{R}_+ & 0
\end{array} \right): L^2_{\theta}(\T^2)\times \comp \rightarrow 
L^2_{\theta}(\T^2)\times \comp,
\endeq
defined microlocally near $\xi=0$, has a microlocal inverse there, which 
has the form 
\begeq
\label{6.26}
\widehat{{\cal E}}= 
\left( \begin{array}{ccc}
\widehat{E}(z)& \widehat{E}_+\\
\widehat{E}_- &  \widehat{E}_{-+}(z)
\end{array} \right).
\endeq
The following localization properties can be inferred from the construction 
of $\widehat{{\cal E}}$: if 
$\psi \in C^{\infty}_b (T^* \T^2)$ has its support disjoint from $\xi=0$, 
then it is true 
that $\psi \widehat E_+ = {\cal O}(h^{\infty}): \comp \rightarrow 
L^2_{\theta}$, and 
$\widehat E_- \psi = {\cal O}(h^{\infty}): L^2_{\theta} \rightarrow 
\comp$.
We also find that 
\begeq
\label{6.27}
\widehat{E}_{-+}(z) = 
z-\widehat{P}\left(h(k-\theta),\eps,\frac{h}{\eps};h\right). 
\endeq	
Using (\ref{6.24}), we furthermore see that the following estimates hold true, 
$$
\widehat{E}=\frac{{\cal O}(1)}{\eps h}: L^2_{\theta}(\T^2)\rightarrow L^2_{\theta}(\T^2),
$$
$$
\widehat{E}_+={\cal O}(1): \comp \rightarrow L^2_{\theta}(\T^2), \quad 
\widehat{E}_-={\cal O}(1): L^2_{\theta}(\T^2) \rightarrow \comp, 
$$
$$
\widehat{E}_{-+}={\cal O}(\eps h): \comp \rightarrow \comp,
$$
so that 
\begeq
\label{6.27.1}
\eps h\norm{u}+\norm{u_-}\leq {\cal O}(1) \left (\norm{v}+\eps h\norm{v_+}\right), 
\endeq 
when 
$$
\widehat{{\cal P}} \left( \begin{array}{ccc}
u \\
u_-
\end{array} \right)=\left( \begin{array}{ccc}
v \\
v_+
\end{array} \right).
$$
In (\ref{6.27.1}), the norms of $u$ and $v$ are taken in $L^2_{\theta}(\T^2)$ 
and those of $u_-$ and $v_+$ in $\comp$. 

\vskip 4mm
\noindent
Passing to the case of $P_{\eps}$, we define 
$R_+: H(\widehat{\Lambda}_{\eps})\rightarrow \comp$ and 
$R_-: \comp\rightarrow H(\widehat{\Lambda}_{\eps})$ by 
\begeq
\label{6.28}
R_+ u = \widehat{R}_+ U \chi u = (U\chi u|e_k), \quad R_- u_- = V \widehat{R}_- u_- = u_- V e_{k}. 
\endeq
It is then true that 
\begeq
\label{6.29}
\chi R_- = R_- +{\cal O}(h^{\infty}): \comp \rightarrow H(\widehat{\Lambda_{\eps}}),
\endeq
decreasing the support of $\chi$ if necessary. 
We now claim that for $z\in \Omega_k(h)$, the Grushin problem 
\begin{eqnarray}
\label{6.30}
\cases{\left(P_{\eps}-z\right)u+R_-u_-=v, \cr
R_+ u=v_+ \cr}
\end{eqnarray}
has a unique solution $(u,u_-)\in H(\widehat{\Lambda}_{\eps})\times
\comp$ for every $(v,v_+)\in H(\widehat{\Lambda}_{\eps})\times
\comp$, with an a priori estimate,
\begeq
\label{6.31}
\eps h\norm{u}+\norm{u_-}\leq {\cal O}(1) \left(\norm{v}+\eps h \norm{v_+}\right).
\endeq
Here the norms of $u$ and $v$ are taken in $H(\widehat{\Lambda}_{\eps})$, 
and those of $u_-$ and $v_+$ in $\comp$. 
To verify the claim, we first see that as in (\ref{6.21}), we have 
\begeq
\label{6.32}
\norm{(1-\chi)u}\leq {\cal O}\left(\frac{1}{\eps}\right) \norm{v}+ 
{\cal O}(h^{\infty})\left(\norm{u}+\norm{u_-}\right).
\endeq
Here we have also used (\ref{6.29}). 

Applying $\chi$ to the first equation in 
(\ref{6.30}) we get 
\begin{eqnarray}
\label{6.33}
\cases{\left(P_{\eps}-z\right)\chi u + R_-u_-=\chi v + w + R_{-\infty}u_-, \cr
R_+ u=v_+ , \cr}
\end{eqnarray} 
where $w=[P_{\eps},\chi]u$ satisfies 
$$
\norm{w}\leq {\cal O}\left(\frac{h}{\eps}\right)\norm{v}+{\cal O}(h^{\infty})\left(\norm{u}+\norm{u_-}\right),
$$
and $R_{-\infty}={\cal O}(h^{\infty})$ in the operator norm. 
Applying $U$ to the first equation in (\ref{6.33}) and using (\ref{6.28}), we get 
\begeq
\label{6.34}
\cases{(\widehat{P}-z)U\chi u+\widehat{R}_-u_-=U \chi v+ U w+w_-\cr
\widehat{R}_+ U \chi u=v_+. \cr}
\endeq 
where the $L^2_{\theta}(\T^2)$-norm of $w_-$ is ${\cal O}(h^{\infty})\left(\norm{u}+\norm{u_-}\right)$. 
We therefore get a microlocally well-posed Grushin 
problem for $\widehat{{\cal P}}$ in (\ref{6.25}), and in view of (\ref{6.27.1}) we obtain,  
\begeq
\label{6.35}
\eps h\norm{\chi u}+ \norm{u_-}\leq {\cal O}(1) \left(\norm{v}+\eps h\norm{v_+}\right) +
{\cal O}(h^{\infty})\left(\norm{u}+\norm{u_-}\right). 
\endeq
Combining (\ref{6.32}) and (\ref{6.35}), we get (\ref{6.31}). We have thus also proved that the operator  
\begeq
\label{6.36}
{\cal P}= 
\left( \begin{array}{ccc}
P_{\eps}-z & R_- \\
R_+ & 0
\end{array} \right): H(\widehat{\Lambda}_{\eps})\times \comp \rightarrow 
H(\widehat{\Lambda}_{\eps})\times \comp
\endeq
is injective, for $z\in \Omega_k(h)$. Now ${\cal P}$ is a finite rank 
perturbation of 
$$
\left( \begin{array}{ccc}
P_{\eps}-z & 0 \\
0 & 0
\end{array} \right),
$$
which is a Fredholm operator of index zero. It follows that ${\cal P}$ is 
also Fredholm of index $0$ and hence bijective, 
since we already know that it is injective. The inverse of ${\cal P}$ has 
the form 
\begeq
\label{6.37}
{\cal E}= 
\left( \begin{array}{ccc}
E(z) & E_+ \\
E_- & E_{-+}(z)
\end{array} \right), 
\endeq
and we recall that the spectrum of $P_{\eps}$ in $\Omega_k(h)$ will be the 
set of values $z$ for which $E_{-+}(z)=0$. 

We finally claim that the components $E_+$ and $E_{-+}(z)$ in (\ref{6.37}) 
are given by $E_+=V \widehat{E}_+$, and 
$E_{-+}(z)=\widehat{E}_{-+}(z) = 
z-\widehat{P}\left(h(k-\theta),\eps,\frac{h}{\eps};h\right)$, modulo terms 
that are ${\cal O}(h^{\infty})$. Indeed, we need 
only to check that 
\begeq
\label{6.38}
R_+ V \widehat{E}_+ \equiv 1, \quad \left(P_{\eps}-z\right)V \widehat{E}_+ + R_- \widehat{E}_{-+} \equiv 0,
\endeq
modulo ${\cal O}(h^{\infty})$, and at this stage the verification of 
(\ref{6.38}) is identical to the corresponding computation from 
section 6 of \cite{MeSj}. In particular, we get 
\begeq
\label{6.39}
E_{-+}(z) = z - 
\widehat{P}\left(h(k-\theta),\eps,\frac{h}{\eps};h\right) + {\cal 
O}(h^{\infty}),
\endeq  
and we have now proved the first of our two main results. 

\begin{theo}\label{Th6.2}
Let $F_0$ be a regular value of $\Re \langle{q}\rangle$ viewed as a 
function on $p^{-1}(0)\cap \real^4$. Assume that the Lagrangian manifold 
$$
\Lambda_{0,F_0}: p=0, \Re \langle{q}\rangle=F_0
$$
is connected, and that $T(0)$ is the minimal period of every closed 
$H_p$-trajectory in $\Lambda_{0,F_0}$. When $\alpha_1$ and $\alpha_2$ are 
the fundamental cycles in $\Lambda_{0,F_0}$ with $\alpha_1$ corresponding 
to a closed $H_p$-trajectory of minimal period, we write $S=(S_1,S_2)$ and 
$k_0=(k_0(\alpha_1),k_0(\alpha_2))$ for the actions and Maslov indices of 
the cycles, respectively. Assume furthermore that 
$\eps = {\cal O}(h^{\delta})$, $\delta>0$, is such that $h/{\eps}\ \ll 1$. 
Let $C>0$ be sufficiently large. Then the 
eigenvalues of $P_{\eps}$ in the rectangle 
\begeq
\label{6.40}
\abs{\Re z}<\frac{1}{C}, \quad \abs{\Im z-\eps F_0}<\frac{\eps}{C}
\endeq
are given by 
$$
z_k = \widehat{P}\left( 
h\left(k-\frac{k_0}{4}\right)-\frac{S}{2\pi},\eps,\frac{h}{\eps}; h 
\right),\quad k\in \z^2, 
$$
modulo ${\cal O}(h^{\infty})$. Here $\widehat{P} 
\left(\xi,\eps,\frac{h}{\eps}; h\right)$ is \hol{} in $\xi\in {\rm 
neigh}(0,\comp^2)$, smooth in 
$\eps, \frac{h}{\eps} \in {\rm neigh}(0,\real)$ and has an asymptotic 
expansion in the space of such functions, 
$$
\widehat{P} \left(\xi,\eps,\frac{h}{\eps}; h\right) \sim p(\xi_1)+ 
\eps\left(r_0\left(\xi, \eps, \frac{h}{\eps}\right)+ 
h r_1\left(\xi, \eps, \frac{h}{\eps}\right)+ \ldots \right), \quad 
h\rightarrow 0,
$$
with 
$$
r_0=i\langle{q}\rangle+ {\cal O}(\eps+h/\eps), 
$$
and $r_{\nu}={\cal O}(\eps+h/\eps)$, $\nu \geq 1$. We have exactly one 
eigenvalue for each $k\in \z^2$ such that the corresponding $z_k$ falls 
into the region {\rm (\ref{6.40})}. 
\end{theo}  

\vskip 2mm
Keeping all the general assumptions of the torus case and still taking 
$F_0=0$, we shall next consider the case when the subprincipal symbol of 
the 
unperturbed operator $P_{\eps=0}$ vanishes. It follows then from the 
previous arguments, now making use of the full strength of Theorem 
\ref{Th1.4}, that in this case, microlocally near $\Lambda_{0,0}$, 
\begeq
\label{6.41}
P_{\eps}: H(\Lambda_{\eps G})\rightarrow H(\Lambda_{\eps G})
\endeq
is equivalent to an operator of the form 
\begeq
\label{6.42}
\widetilde{P}_{\eps}(x_2,\xi,\eps;h) \sim \sum_{\nu=0}^{\infty} h^{\nu} 
\widetilde{p}_{\nu}(x_2,\xi,\eps),
\endeq
acting on $L^2_{\theta}(\T^2)$, with 
$$
\widetilde{p}_0(x_2,\xi,\eps)=p(\xi_1)+i\eps\langle{q}\rangle(\xi)+{\cal 
O}(\eps^2), \quad \widetilde{p}_1(x_2,\xi,\eps)=\eps q_1(x_2,\xi,\eps).
$$

In what follows we shall discuss the the range 
\begeq
\label{6.43}
Mh^2 < \eps = {\cal O}(h^{\delta})\quad M \gg 1,\,\, \delta>0. 
\endeq
Recalling the operators $e^{\eps G(x,hD_x,\eps)/h}$ and $U$ from 
(\ref{6.4}) and (\ref{6.10}), respectively, we see, as in the general case, 
that the symbol of $\Im P_{\eps}$ on $H(\Lambda_{\eps G})$ is ${\cal O}(\eps)$, and away 
from any fixed \neigh{} of $\Lambda_{0,0}$ in $\Lambda_{\eps G}$, 
we have $\abs{\Im P_{\eps}(\rho,h)}\sim \eps $, if $\abs{\Re P_{\eps}(\rho,h)}< 1/{\cal O}(1)$. 

We write, as in section \ref{section3}, 
$$
\widetilde{P}(x_2,\xi,\eps,h)=p(\xi_1)+\eps 
\left(r_0\left(x_2,\xi,\eps,\frac{h^2}{\eps}\right)+h 
r_1\left(x_2,\xi,\eps,\frac{h^2}{\eps}\right)
+\ldots\right),
$$
where 
$$
r_0\left(x_2,\xi,\eps,\frac{h^2}{\eps}\right)=i\langle{q}\rangle+{\cal 
O}(\eps)+\frac{h^2}{\eps}\widetilde{p}_2(x_2,\xi,\eps), 
$$
$$
r_1(x_2,\xi,\eps)=q_1(x_2,\xi,\eps)+\frac{h^2}{\eps}\widetilde{p}_3(x_2,\xi,\eps),
\quad r_j(x_2,\xi,\eps)={\cal O}\left(\frac{h^2}{\eps}\right),\,\, j\geq 2.
$$
Using the canonical transformation $\kappa$, generated by the function 
$$
\psi\left(x,\eta, \eps, \frac{h^2}{\eps}\right)=x\cdot \eta+\phi_{\rm 
per}\left(x_2,\eta,\eps, \frac{h^2}{\eps}\right), 
$$
with $\phi_{\rm per}={\cal O}(\eps+\frac{h^2}{\eps})$, constructed in 
section \ref{section3}, we then argue similarly to the general torus case. We thus 
introduce an IR-manifold $\widetilde{\Lambda}\subset 
\widetilde{\T^2}\times \comp^2$ which is an $(\eps+h^2/\eps)$-perturbation 
of $T^* \T^2$,
which agrees with $\kappa(T^* \T^2)$ near $\xi=0$, and further away from 
this set coincides with $T^*\T^2$. When constructing 
$\widetilde{\Lambda}$, 
we first notice that $\kappa(T^*\T^2)$ has the form 
$$
\Im x = G'_{\xi}(\Re (x,\xi)), \quad \Im \xi = - G'_x(\Re(x,\xi)),
$$
where $G=G(x_2,\xi,\eps, \frac{h^2}{\eps})$ is such that 
$$
\partial_{\xi} G, \, \partial_{x_2} G={\cal 
O}\left(\eps+\frac{h^2}{\eps}\right). 
$$
As was observed in section \ref{section3}, the transformation $\kappa$ conserves 
actions, and therefore the smooth function $G$ is single-valued. We may 
assume that 
$$
G={\cal O}\left(\eps+\frac{h^2}{\eps}\right).
$$
If we let $\chi(\xi) \in C_0^{\infty}(\real^2; [0,1])$ be a cutoff 
function with a small support and such that $\chi=1$ in a small \neigh{} 
of $0$, 
we define $\widetilde{\Lambda}$ by 
$$
\Im x=\widetilde{G}'_{\xi}(\Re(x,\xi)), \quad \Im \xi = 
-\widetilde{G}'_{x}(\Re(x,\xi)),\,\, \widetilde{G}(\Re (x,\xi))=\chi(\Re 
\xi)G(\Re (x,\xi)).
$$
We then obtain the desired globally defined IR-manifold 
$\widetilde{\Lambda}$ such that $\Im \xi_1=0$ on $\widetilde{\Lambda}$. 
When acting on 
$H(\widetilde{\Lambda})$, $\widetilde{P_{\eps}}$ is microlocally near 
$\xi=0$ unitarily equivalent to an operator on $L^2(\T^2)$, which has the 
form 
$$
p(\xi_1)+\eps \left(r_0\left(\xi,\eps,\frac{h^2}{\eps}\right)+h 
r_1\left(x_2,\xi,\eps,\frac{h^2}{\eps}\right)+\ldots\right), 
$$
where 
$$
r_0\left(\xi,\eps,\frac{h^2}{\eps}\right)=i\langle{q}\rangle+{\cal 
O}\left(\eps+\frac{h^2}{\eps}\right)
$$
is independent of $x$. 

It follows, as in the general torus case, that on the Bargmann transform 
side, $\widetilde{\Lambda}$ can be described by an 
FBI-weight $\Phi=\Phi(\Im x_1, x_2)$ 
which does not depend on $\Re x_1$. Repeating the previous arguments, we 
obtain therefore a new globally defined Hilbert space 
$H(\widehat{\Lambda})$, associated to an IR-manifold 
$\widehat{\Lambda}\subset \comp ^4$, and a Lagrangian torus 
$\widehat{\Lambda}_{0,0}\subset \widehat{\Lambda}$ such that microlocally 
near $\widehat{\Lambda}_{0,0}$, 
$P_{\eps}: H(\widehat{\Lambda})\rightarrow H(\widehat{\Lambda})$ is 
equivalent to an operator on $L^2_{\theta}(\T^2)$, described in (\ref{3.29}), 
(\ref{3.30}). 

\begin{prop}\label{Prop6.3}
Assume that the subprincipal symbol of $P_{\eps=0}$ vanishes, and consider 
the range $Mh^2 < \eps ={\cal O}(h^{\delta})$ for $M \gg 1$, $\delta>0$. There exists an 
IR-manifold $\widehat{\Lambda}\subset \comp^4$ and a smooth Lagrangian torus 
$\widehat{\Lambda}_{0,0}\subset \widehat{\Lambda}$ such that when 
$\rho\in \widehat{\Lambda}$ is away from a small \neigh{} of 
$\widehat{\Lambda}_{0,0}$ in 
$\widehat{\Lambda}$ and $\abs{\Re P_{\eps}(\rho,h)}<1/C$, for a 
sufficiently large $C>0$, it is true that 
$$
\abs{\Im P_{\eps}(\rho,h)}\sim \eps. 
$$
The manifold $\widehat{\Lambda}$ is $(\eps+h^2/\eps)$-close to $\real^4$ 
and it coincides with $\real^4$ outside a \neigh{} of 
$p^{-1}(0)\cap \real^4$. There exists a \cantransf{} 
$$
\kappa_{\eps}: {\rm 
neigh}(\widehat{\Lambda}_{0,0},\widehat{\Lambda})\rightarrow {\rm 
neigh}(\xi=0, T^*\T^2), 
$$
mapping $\widehat{\Lambda}_{0,0}$ onto $\T^2$, and an elliptic Fourier 
integral operator $U: H(\widehat{\Lambda})\rightarrow L^2_{\theta}(\T^2)$ 
associated to $\kappa_{\eps}$, such that, microlocally near $\widehat{\Lambda}_{0,0}$, 
$UP_{\eps}=\widehat{P}U$. Here 
$$
\widehat{P}=\widehat{P}(hD_x,\eps,\frac{h^2}{\eps};h)
$$
has the Weyl symbol, depending smoothly on $\eps$, $h^2/\eps\in {\rm neigh}(0,\real)$, 
$$
\widehat{P}\left(\xi,\eps,\frac{h^2}{\eps};h\right) \sim p(\xi_1)+\eps 
\sum_{j=0}^{\infty} h^j r_j\left(\xi,\eps,\frac{h^2}{\eps}\right).
$$
We have 
$$
r_0=i\langle{q}\rangle(\xi)+{\cal O}(1) (\eps+h^2/\eps),\quad r_j={\cal 
O}(1), \,\, j\geq 1.
$$
\end{prop}

\vskip 2mm
Repeating the arguments, leading to Theorem \ref{Th6.2}, and using 
Proposition \ref{Prop6.3} 
instead of Proposition \ref{Prop6.1}, we then find first 
that the spectrum of $P_{\eps}$ in a region of the form (\ref{6.19}) is contained in 
the union of disjoint discs of radii $\eps h/{\cal O}(1)$ around the 
quasi-eigenvalues 
$\widehat{P}\left(h(k-\theta), \eps, h^2/ \eps;h\right)$. Furthermore, 
when $z$ varies in such a disc corresponding to $k\in \z^2$, such that 
the corresponding quasi-eigenvalue falls into the region (\ref{6.19}),
an inspection of the previous arguments shows that the Grushin operator 
$$
\left( \begin{array}{ccc}
P_{\eps}-z & R_- \\
R_+ & 0
\end{array} \right): H(\widehat{\Lambda})\times \comp \rightarrow 
H(\widehat{\Lambda})\times \comp
$$
is bijective with the inverse of the norm 
${\cal O}((\eps h)^{-1})$---see (\ref{6.31}) for the precise a priori estimate. 
Here $R_-: \comp \rightarrow H(\widehat{\Lambda})$ and 
$R_+: H(\widehat{\Lambda})\rightarrow \comp$ are defined as in (\ref{6.28}). 
This leads to the following result.  

\begin{theo}\label{Th6.4}
Keep all the assumptions and notation of Theorem {\rm \ref{Th6.2}}, and in addition 
assume that the subprincipal symbol of $P_{\eps=0}$ vanishes. Let 
$\eps = {\cal O}(1) h^{\delta}$ for some fixed $\delta>0$ be such that 
$h^2 \ll \eps$. Then the eigenvalues of $P_{\eps}$ in the rectangle 
$$
\left(-\frac{1}{{\cal O}(1)}, \frac{1}{{\cal O}(1)}\right)+i\eps 
\left(F_0-\frac{1}{{\cal O}(1)}, F_0+\frac{1}{{\cal O}(1)}\right)
$$
are given by 
$$
z_k =\widehat{P}\left(h\left(k-\frac{k_0}{4}\right)-\frac{S}{2\pi},\eps, 
\frac{h^2}{\eps};h\right)+{\cal O}(h^{\infty}), \,\, k\in \z^2.
$$
Here $\widehat{P}(\xi,\eps, h^2/ \eps;h)$ is holomorphic in $\xi \in {\rm 
neigh}(0,\comp^2)$, smooth in 
$\eps$ and $h^2/ \eps\in {\rm neigh}(0,\real)$, and as $h\rightarrow 0$, 
there is an asymptotic expansion 
$$
\widehat{P}\left(\xi,\eps,\frac{h^2}{\eps};h\right)\sim p(\xi_1)+\eps 
\left(r_0\left(\xi,\eps,\frac{h^2}{\eps}\right)+
h r_1\left(\xi,\eps, \frac{h^2}{\eps}\right)+\ldots\right).
$$
We have 
$$
r_0\left(\xi,\eps,\frac{h^2}{\eps}\right)=i\langle{q}\rangle(\xi)+{\cal 
O}\left(\eps+\frac{h^2}{\eps}\right), \quad 
r_j\left(\xi, \eps, \frac{h^2}{\eps}\right)={\cal O}(1), \,\, j\geq 1. 
$$
\end{theo}

\vskip 4mm
We shall now turn to the case II from the introduction. Let us recall from section 1, that if 
$z\in \Spec{P_{\eps}}$ is such that $\abs{\Re z}\leq \delta \rightarrow 
0$, 
then 
\begeq
\label{6.44}
\Im z \in \eps \bigl [ \inf_{\Sigma} \Re \langle{q}\rangle -o(1), 
\sup_{\Sigma} \Re \langle{q}\rangle+o(1)\bigr], \,\, h\rightarrow 0. 
\endeq
Here, as in section 1, we write $\Sigma=p^{-1}(0)\cap \real^4/ 
\exp(\real H_p)$. Our purpose here is to show that the quasi-eigenvalues 
of Propositions \ref{Prop4.4} and \ref{Prop4.5} give, up to ${\cal O}(h^{\infty})$, the actual 
eigenvalues in a set of the form 
$$
\abs{\Re z}\leq \frac{1}{{\cal O}(1)}, \quad \abs{\Im z-\eps F_0} \leq 
\frac{\eps}{{\cal O}(1)}, 
$$
when $F_0\in \{\inf_{\Sigma}\Re \langle{q}\rangle, \sup_{\Sigma} \Re 
\langle{q}\rangle \}$. As we shall see, the analysis here will be parallel 
to the 
torus case just treated, so that in what follows we shall concentrate on 
the new features of the problem, and some of the computations that are 
essentially identical to the ones already performed, will not be repeated. 

In order to fix the ideas, we shall discuss the case when 
$$
F_0=\inf_{\Sigma} \Re \langle{q}\rangle,
$$
and we shall take $F_0=0$. 

\vskip 2mm
Recall from the beginning of this section that the original operator 
$P_{\eps}$ acting on $H(\Lambda_{\eps G})$, is microlocally 
unitarily equivalent to the operator 
\begeq
\label{6.45}
P_{\eps} \sim \sum_{j=0}^{\infty} h^j p_j(x,\xi,\eps),
\endeq
acting on $L^2$ and defined microlocally near $p^{-1}(0)\cap \real^4$, 
with 
$$
p_0=p+i\eps\langle{q}\rangle+{\cal O}(\eps^2),
$$
and the functions $\langle{q}\rangle$ and ${\cal O}(\eps^2)$-term are in 
involution with $p$. Let $\gamma_1,\ldots \gamma_N\subset p^{-1}(0)\cap \real^4$ 
be the finitely many trajectories such that $\Re \langle{q}\rangle=0$ along 
$\gamma_j$, $1\leq j\leq N$. We know that $T(0)$ is the minimal period of each $\gamma_j$, and if 
$\rho_j\in \Sigma$ is the corresponding point, then the Hessian of $\Re\langle{q}\rangle$ at 
$\rho_j$ is positive definite, $1\leq j\leq N$. 
Associated to $\gamma_j$, we have the quantities $S=S(\gamma_j)$ and 
$k_0=k_0(\gamma_j)$, the action along $\gamma_j$ and the Maslov index, 
respectively, defined as in section \ref{section1}, and we recall from~\cite{HeRo} that these 
quantities do not depend on $j$. 

In what follows we shall work microlocally near a fixed critical trajectory, 
say $\gamma_1$. We let $L^2_S (S^1\times \real)$ be the space of locally square integrable 
functions $u(t,x)$ on $\real\times\real$ such that 
$$
\int\!\!\!\int_0^{2\pi} \abs{u(t,x)}^2\,dx\,dt<\infty.
$$
and 
$$
u(t-2\pi,x)=e^{iS/h+ik_0\pi/2}u(t,x).
$$
Applying Theorem \ref{Th1.4} to the \cantransf{} $\kappa$ of 
Proposition \ref{Prop2.1}, we see that there exists an analytic 
microlocally unitary \fourior{} 
$$
U_0: L^2_S(S^1\times \real)\rightarrow L^2(\real^2),
$$
associated to $\kappa$, and defined microlocally from a \neigh{} of 
$\{\tau=x=\xi=0\}$ in 
$T^*\left(S^1\times \real\right)$ to a \neigh{} of $\gamma_1$ in 
$\real^4$, so that we have the two-term Egorov property (\ref{1.3}). 
Combining $\exp(i\eps H_G)$ with $\kappa$, we get a smooth \cantransf{} 
\begeq
\label{6.46}
\kappa_{\eps}: {\rm neigh}\left(\tau=x=\xi=0, T^*\left(S^1\times 
\real\right)\right) \rightarrow {\rm neigh}(\gamma_1, \Lambda_{\eps G}), 
\endeq
where abusing the notation slightly, we write here $\gamma_1\subset 
\Lambda_{\eps G}$ also for the image of $\gamma_1$ under the 
complex \cantransf{} $\exp(i\eps H_G)$. The operator $e^{\eps 
G(x,hD_x,\eps)/h}\circ U_0$ is then associated with $\kappa_{\eps}$ in 
(\ref{6.46}), and an
application of Egorov's theorem shows that, microlocally near $\gamma_1$, 
we get a unitary equivalence between the operator $P_{\eps}$ acting on 
$H(\Lambda_{\eps G})$ and an $h$-\pseudor{} microlocally defined near 
$\tau=x=\xi=0$ in $T^*\left(S^1 \times \real\right)$, with the leading 
symbol 
$$
\widetilde{p}_0(\tau,x,\xi,\eps)=f(\tau)+i\eps\langle{q}\rangle(\tau,x,\xi)+{\cal 
O}(\eps^2), 
$$
independent of $t$. Taking into account an additional conjugation by the 
elliptic operator $e^{iA/h}$, acting on $L^2_S(S^1\times \real)$, with 
$$
A\sim \sum_{k=1}^{\infty} a_k(t,\tau,x,\xi,\eps)h^k,
$$
constructed as a formal power series in $\eps$, $h$ in Proposition 
\ref{Prop2.2}, we 
see that microlocally near $\gamma_1$, the operator 
$P_{\eps}: H(\Lambda_{\eps G})\rightarrow H(\Lambda_{\eps G})$ is 
equivalent to an operator of the form 
\begeq
\label{6.47}
\widetilde{P}_{\eps}(\tau,x,\xi,\eps) \sim \sum_{k=0}^{\infty} h^k 
\widetilde{p}_k(\tau,x,\xi,\eps),
\endeq
acting on $L^2_S(S^1\times \real)$, whose full symbol is independent of 
$t$. We have 
\begeq
\label{6.48}
\widetilde{p}_0=f(\tau)+i\eps\langle{q}\rangle(\tau,x,\xi)+{\cal 
O}(\eps^2), 
\endeq
and 
$$
\Re \langle{q}\rangle(0,x,\xi)\sim x^2+\xi^2
$$
on the real domain.

We shall first consider the general case when the subprincipal symbol of 
the unperturbed operator $P_{\eps=0}$ does not necessarily vanish, and 
in doing so, it will be assumed that 
\begeq
\label{6.49}
h\ll \eps ={\cal O}(1) h^{\delta}, \quad \delta>0.
\endeq
As in section \ref{section4}, we write 
$$
\widetilde{P}_{\eps} 
=f(\tau)+\eps\left(i\langle{q}\rangle(\tau,x,\xi)+{\cal O}(\eps)+
\frac{h}{\eps}\widetilde{p}_1+h\frac{h}{\eps}\widetilde{p}_2+\ldots\right).
$$
According to Proposition \ref{Prop4.3}, there exists a \hol{} \cantransf{} 
$$
\kappa_{\sigma,\eps,\frac{h}{\eps}}: {\rm neigh}(0,\comp^2)\rightarrow 
{\rm neigh}(0,\comp^2), 
$$
depending analytically on $\sigma\in {\rm neigh}(0,\comp)$ and smoothly on 
$\eps, \frac{h}{\eps}\in {\rm neigh}(0,\real)$, such that 
$$
\Im \kappa_{\sigma,\eps,\frac{h}{\eps}}(y,\eta)={\cal 
O}\left(\eps+\frac{h}{\eps}\right),
$$
when $\sigma,y,\eta$ are real, and such that 
$$
\left(\widetilde{p}_0+\eps\frac{h}{\eps}\widetilde{p}_1\right)\left(\sigma, 
\kappa_{\sigma,\eps,\frac{h}{\eps}}(y,\eta)\right)=
f(\sigma)+i\eps g_{\eps, \frac{h}{\eps}}\left(\sigma, 
\frac{y^2+\eta^2}{2}\right). 
$$
Here $g_{\eps, \frac{h}{\eps}}(\sigma,q)$ is an analytic function, 
depending smoothly on $\eps, h/\eps$, for which 
$$
\frac{\partial}{\partial q} \Re g_{\eps,0}(0,0)> 0.
$$
We now lift the family of locally defined canonical transformations 
$\kappa_{\sigma,\eps,\frac{h}{\eps}}$ to a canonical transformation
\begin{eqnarray*}
& & \Xi_{\eps,\frac{h}{\eps}}: {\rm neigh}\left(\Im s=0, \sigma=y=\eta=0, 
T^*\left(\widetilde{S^1}\times \comp\right) \right) \ni (s,\sigma; y,\eta) 
\\
& \mapsto & (t,\tau;x,\xi) \in {\rm neigh}\left(\Im t=0, \tau=x=\xi=0, 
T^*\left(\widetilde{S^1}\times \comp\right) \right)
\end{eqnarray*}
given by 
\begeq
\label{6.50}
\Xi_{\eps,\frac{h}{\eps}}(s,\sigma;y,\eta)=(t,\tau;x,\xi)=(s+h(y,\sigma,\eta),\sigma; 
\kappa_{\sigma,\eps,\frac{h}{\eps}}(y,\eta)).
\endeq 
Here $h(y,\sigma,\eta)$ is uniquely determined up to a function 
$g=g(\sigma)$, and if $\varphi_{\sigma,\eps,\frac{h}{\eps}}(x,y,\theta)$ 
is an analytic 
family of non-degenerate phase functions (in the sense of H\"ormander) 
locally generating the family $\kappa_{\sigma,\eps,\frac{h}{\eps}}$, then 
$$
\Phi_{\eps,\frac{h}{\eps}}(t,x,s,y,\theta,\sigma):=
\varphi_{\sigma,\eps,\frac{h}{\eps}}(x,y,\theta)+(t-s)\sigma
$$
is a non-degenerate phase function with $\theta,\sigma$ as fiber 
variables, such that $\Phi_{\eps,\frac{h}{\eps}}$ generates the graph of 
$\Xi_{\eps,\frac{h}{\eps}}$. 

Associated to $\Xi_{\eps,\frac{h}{\eps}}$, we introduce an IR-manifold 
$\widetilde{\Lambda}\subset T^*\left(\widetilde{S^1}\times \comp\right)$, 
which in a complex \neigh{} of $\tau=x=\xi=0$, is equal to 
$\Xi_{\eps,\frac{h}{\eps}}\left(T^*\left(S^1\times \real\right)\right)$, 
and further 
away from this set agrees with $T^*\left(S^1\times \real\right)$. In the 
intermediate region, we construct $\widetilde{\Lambda}$ in such a way 
that it remains an $(\eps+\frac{h}{\eps})$-perturbation of 
$T^*\left(S^1\times \real\right)$, and so that everywhere on 
$\widetilde{\Lambda}$, it is 
true that 
\begeq
\label{6.51}
(t,\tau; x,\xi)\in \widetilde{\Lambda} \Longrightarrow \tau\in \real. 
\endeq
If we now use the standard FBI-Bargmann transformation, viewed as a 
mapping on $L^2_S(S^1 \times \real)$, 
so that under the associated canonical transformation, $T^*(S^1\times 
\real)$ is mapped to 
$\{ (t,\tau; x,\xi)\in T^* (\widetilde{S^1}\times \comp); (\tau,\xi)=-\Im 
(t,x) \}$, 
then as before we see that after an application of such a transformation, 
the manifold $\widetilde{\Lambda}$ is described by a 
weight function $\Phi=\Phi(\Im t,x)$ which does not depend on $\Re t$. At 
this stage, the situation is similar to the previously 
analyzed torus case, and, in particular, we see again that the form of the 
weight $\Phi(\Im t,x)$ implies that the term $f(\tau)$ in (\ref{6.48}) 
gives no contribution to the imaginary part of the operator. Summing up the 
discussion so far, we arrive to the following result.
 
\begin{prop}\label{Prop6.5}
Make the assumptions of case {\rm II} in the introduction, and assume that 
$$
F_0=\inf_{\Sigma} \Re \langle{q}\rangle=0.
$$
Assume that $\eps={\cal O}(h^{\delta})$, for some $\delta>0$, is such that 
$h\ll \eps$. There exists a closed IR-manifold 
$\Lambda\subset \comp^4$ and finitely many simple closed disjoint 
curves $\gamma_1,\ldots \gamma_N \subset \Lambda$, 
which are $(\eps+h/\eps)$-close to the closed $H_p$-trajectories $\subset p^{-1}(0)\cap \real^4$, 
along which $\Re\langle{q}\rangle=0$, such 
that when $\rho$ is outside a small \neigh{} of $\cup_{j=1}^N \gamma_j$ in $\Lambda$, then
\begeq
\label{6.52}
\abs{\Re P_{\eps}(\rho,h)}\geq \frac{1}{{\cal O}(1)} \quad \wrtext{or} 
\quad \abs{\Im P_{\eps}(\rho,h)} \geq \frac{\eps}{{\cal O}(1)}. 
\endeq
This estimate is true away from an arbitrarily small \neigh{} 
of $\cup_{j=1}^N \gamma_j$, provided that the implicit constant in {\rm (\ref{6.52})} is 
chosen sufficiently large. The manifold $\Lambda$ coincides with $\real^4$ 
away from a \neigh{} of $p^{-1}(0)\cap \real^4$ and is 
$(\eps+h/\eps)$-close to 
$\real^4$ everywhere. For each $j$ with $1\leq j\leq N$, there exists a \cantransf
$$
\kappa_{\eps,j}: {\rm neigh}\left(\gamma_j, \Lambda\right) \rightarrow {\rm 
neigh}\left(\tau=x=\xi=0, T^*(S^1\times \real)\right),
$$
whose domain of definition does not intersect the closure of the union of the domains 
of the $\kappa_{\eps,k}$ for $k\neq j$, and an elliptic \fourior 
$$
U_j = {\cal O}(1): H(\Lambda) \rightarrow L^2_S(S^1\times \real),
$$
associated to $\kappa_{\eps,j}$, such that, microlocally near $\gamma_j$, 
$$
U_j P_{\eps} = \widehat{P}_jU_j. 
$$
Here $\widehat{P}_j=\widehat{P}_j(hD_t, (1/2)(x^2+(hD_x)^2), \eps, 
\frac{h}{\eps};h)$ has the Weyl symbol  
$$
\widehat{P}_j\left(\tau,x,\xi,\eps,\frac{h}{\eps};h\right)=f(\tau)+i\eps 
G_j\left(\tau, \frac{x^2+\xi^2}{2},\eps,\frac{h}{\eps}; h\right), 
$$
with 
$$
G_j\left(\tau,q,\eps,\frac{h}{\eps}; h\right)\sim \sum_{l=0}^{\infty} h^l 
G_{j,l}\left(\tau,q,\eps,\frac{h}{\eps}\right), \quad h\rightarrow 0, 
$$
and $G_{j,l}$ holomorphic in $(\tau,q)\in {\rm neigh}(0,\comp^2)$, smooth in 
$\eps, h/\eps \in {\rm neigh}(0,\real)$. Furthermore, $\Re G_{j,0}(0,0,0,0)=0$ 
and 
$$
\frac{\partial}{\partial q} \Re G_{j,0}(0,0,0,0)>0.
$$
\end{prop}

Take now small open sets $\Omega_j\subset \Lambda$, $1\leq j\leq N$, 
such that $\gamma_j \subset \Omega_j$ and 
$$
\overline{\Omega_j}\cap \overline{\Omega_k}=\emptyset, \quad j\neq k.
$$
Let $\chi_j\in C^{\infty}_0(\Omega_j)$, $0\leq \chi_j \leq 1$, be such that 
$\chi_j=1$ near $\gamma_j$, $1\leq j\leq N$. 
When $z\in \comp$ satisfies 
\begeq
\label{6.53}
\abs{\Re z}\leq \frac{1}{C}, \quad \abs{\Im z}\leq 
\frac{\eps}{C},
\endeq
and $(P_{\eps}-z)u=v$, it follows from (\ref{6.52}) by repeating the 
arguments of the torus case, that  
\begeq
\label{6.54}
\norm{\left(1-\sum_{j=1}^N \chi_j\right)u} \leq {\cal O}\left(\frac{1}{\eps}\right)\norm{v}+
{\cal O}(h^{\infty})\norm{u}.
\endeq
We shall now discuss the setup of the global Grushin problem. Associated with 
each normal form $\widehat{P_j}$, 
$1\leq j\leq N$, we have the quasi-eigenvalues given in Proposition 5.4, 
$$
z(j,k):=f\left(h(k_1-\frac{k_0}{4})-\frac{S}{2\pi}\right)+i\eps 
G_j\left(h(k_1-\frac{k_0}{4})-\frac{S}{2\pi}, 
h\left(k_2+\frac{1}{2}\right),\eps, \frac{h}{\eps};h\right), 
$$
when $1\leq j\leq N$ and $k=(k_1,k_2)\in \z\times \nat$. We also 
introduce an ON system of eigenfunctions of the (formally) 
commuting operators $\widehat{P}_j$, 
$$
e_k(t,x)=\frac{1}{\sqrt{2\pi}} 
e^{\frac{i}{h}\left(h(k_1-\frac{k_0}{4})-\frac{S}{2\pi}\right)t} 
e_{k_2}(x), \quad 
k=(k_1,k_2)\in \z\times \nat,
$$
which forms an ON basis in $L^2_S(S^1\times \real)$. Here 
$e_{k_2}(x)$, $k_2\in \nat$, are the normalized eigenfunctions of 
$1/2(x^2+(hD_x)^2)$ with eigenvalues $(k_2+1/2)h$.

When $1\leq j\leq N$, let 
$$
M_j=\# \left\{ z(j,k), \abs{\Re z(j,k)} < \frac{1}{{\cal O}(1)}, 
\abs{\Im z(j,k)} < \frac{\eps}{{\cal O}(1)} \right\}.
$$
Then $M_j={\cal O}(h^{-2})$ and we let $k(j,1),\ldots k(j,M_j)\in \z\times \nat$ 
be the corresponding half-lattice points. We introduce the auxiliary operator 
$$
R_+: H(\Lambda)\rightarrow \comp^{M_1}\times\ldots\times \comp^{M_N},
$$
given by 
$$
R_+u(j)(l)=(U_j\chi_j u|e_{k(j,l)}), \quad 1\leq j\leq N, \,\,1\leq l\leq M_j.
$$
Here the inner product in the right-hand side is taken in $L^2_S(S^1\times \real)$. Define also 
$$
R_-: \comp^{M_1}\times\ldots\times \comp^{M_N}\rightarrow H(\Lambda),
$$
by 
$$
R_-u_-=\sum_{j=1}^N \sum_{l=1}^{M_j} u_-(j)(l)V_j e_{k(j,l)}.
$$
Here $V_j$ is a microlocal inverse of $U_j$. We then claim that for $z\in \comp$ satisfying (\ref{6.53}), 
with a sufficiently large $C>0$, the Grushin operator 
\begeq
\label{6.55}
{\cal P}= 
\left( \begin{array}{ccc}
P_{\eps}-z & R_- \\
R_+ & 0
\end{array} \right): H(\Lambda)\times \left(\comp^{M_1}\times\ldots\times \comp^{M_N}\right) \rightarrow 
H(\Lambda)\times \left(\comp^{M_1}\times\ldots\times \comp^{M_N}\right)
\endeq
is bijective. Indeed, when $v\in H(\Lambda)$ and $v_+\in \comp^{M_1}\times \ldots \times \comp^{M_N}$, let 
us consider 
\begeq
\label{6.56}
\cases{(P_{\eps}-z)u+ R_-u_-= v, \cr
R_+ u=v_+. \cr}
\endeq 
As in (\ref{6.54}), we get 
$$
\norm{\left(1-\sum_{j=1}^N \chi_j\right)u} 
\leq {\cal O}\left(\frac{1}{\eps}\right)\norm{v}+{\cal O}(h^{\infty})\left(\norm{u}+\norm{u_-}\right).
$$
Applying $\chi_j$ and then $U_j$, $1\leq j\leq N$, to the first equation in (\ref{6.56}), we get 
\begeq
\label{6.57}
\cases{(\widehat{P}_{j}-z)U_j \chi_j u+ \sum_{l=1}^{M_j} u_-(j)(l)e_{k(j,l)} = 
U_j \left(\chi_j v+[P_{\eps},\chi_j]u\right)+R_{\infty}u+R_{-,\infty}(j)u_-, \cr
(U_j \chi_j u|e_{k(j,l)})=v_+(j)(l),\quad 1\leq l\leq M_j,\cr}
\endeq
and here $R_{\infty}=R_{\infty}(j)={\cal O}(h^{\infty})$ and 
$R_{-,\infty}(j)={\cal O}(h^{\infty})$ in the corresponding operator norms. 
For each $j$, $1\leq j\leq N$, we get a 
microlocally well-posed Grushin problem for $\widehat{P}_j-z$ in $L^2_S(S^1\times \real)$, 
with inverse of the norm 
${\cal O}(1/{\eps})$, and the global well-posedness of (\ref{6.56}) follows. 
The inverse ${\cal E}$ of ${\cal P}$ in (\ref{6.55}) has the form 
\begeq
\label{6.58}
{\cal E}= 
\left( \begin{array}{ccc}
{E}(z) & {E}_+\\
{E}_-  &  {E}_{-+}(z)
\end{array} \right),
\endeq
and a straightforward computation shows that 
$$
E_+: \comp^{M_1}\times \ldots \comp^{M_N} \rightarrow H(\Lambda)
$$
modulo ${\cal O}(h^{\infty})$, is given by 
$$
E_+v_+ \equiv \sum_{j=1}^N \sum_{l=1}^{M_j} v_+(j)(l)V_j e_{k(j,l)}=R_-v_+,
$$
and ${E}_{-+}(z)\in {\cal L}\left(\comp^{M_1}\times\ldots\times \comp^{M_N},
\comp^{M_1}\times\ldots\times \comp^{M_N}\right)$ is 
a block diagonal matrix with the blocks $E_{-+}(z)(j)\in {\cal L}(\comp^{M_j},\comp^{M_j})$, $1\leq j\leq N$, 
of the form 
$$
E_{-+}(z)(j)(m,n)\equiv \left(z-z(j,k(j,m))\right)\delta_{mn},\quad 1\leq m\leq n\leq M_j,
$$
modulo ${\cal O}(h^{\infty})$. The computation of eigenvalues near 
the boundary of the band has therefore 
been justified, and we get the second of our two main results. 

\begin{theo}\label{Th6.6}
Assume that 
$$
F_0=\inf_{\Sigma} \Re \langle{q}\rangle
$$
is achieved along finitely many closed $H_p$-trajectories $\gamma_1, 
\ldots \gamma_N \subset p^{-1}(0)\cap\real^4$ of 
minimal period $T(0)$, and that the Hessian of $\Re\langle{q}\rangle$ at 
the corresponding points $\rho_j \in \Sigma$, 
$j=1, \ldots N$, is positive definite. Let us write $S$ and $k_0$ to 
denote the common values of the action and the Maslov index of 
$\gamma_j$, $j=1,\ldots N$, respectively. Assume that $\eps={\cal 
O}(h^{\delta})$ for a fixed $\delta>0$, is such that 
$h \ll \eps$. Let $C>0$ be sufficiently large. Then the eigenvalues of 
$P_{\eps}$ in the set
\begeq
\label{6.59}
\left(-\frac{1}{C}, \frac{1}{C}\right)+i\eps 
\left(F_0-\frac{1}{C},F_0+\frac{1}{C}\right)
\endeq
are given by 
$$
z(j,k)=f\left(h\left(k_1-\frac{k_0}{4}\right)-\frac{S}{2\pi}\right) +i\eps 
G_j\left(h\left(k_1-\frac{k_0}{4}\right)-\frac{S}{2\pi}, 
h\left(\frac{1}{2}+k_2\right), \eps, \frac{h}{\eps};h\right), 
$$
modulo ${\cal O}(h^{\infty})$, when $1\leq j\leq N$ and $(k_1,k_2)\in \z\times \nat$. Here 
$f(\tau)$ is real-valued with $f(0)=0$ and $f'(0)>0$. The function 
$G_j(\tau,q,\eps,h/\eps;h)$, $1\leq j\leq N$, is analytic in $\tau$ and $q$ in a \neigh{} of 
$(0,0)\in \comp^2$, and smooth in $\eps, h/\eps\in {\rm neigh}(0,\real)$, 
and has an asymptotic expansion in the space of such functions, as 
$h\rightarrow 0$, 
$$
G_j\left(\tau,q,\eps,\frac{h}{\eps};h\right)\sim \sum_{l=0}^{\infty} 
G_{j,l}\left(\tau,q,\eps, \frac{h}{\eps},\right) h^l. 
$$
We have $\Re G_{j,0}(0,0,0,0)=F_0$ and 
$$
\frac{\partial}{\partial q} \Re G_{j,0}(0,0,0,0)>0, \quad 1\leq j\leq N. 
$$
\end{theo}

\noindent
\Remark. With obvious modifications, Theorem \ref{Th6.6} describes the eigenvalues 
in the region (\ref{6.59}), when 
$F_0=\sup_{\Sigma}\Re \langle{q}\rangle$, if we assume that $F_0$ is 
attained along finitely many trajectories 
of minimal period $T(0)$, such that the transversal Hessian of $\Re 
\langle{q}\rangle$ along the trajectories is negative definite. 

\vskip 2mm
\noindent
The treatment of the remaining case of the eigenvalues near the boundary 
of the band (\ref{6.59}), when the subprincipal symbol of 
$P_{\eps=0}$ vanishes proceeds in full analogy with the previously 
analyzed torus case. Thus, restricting attention to the region 
$$
Mh^2 < \eps ={\cal O}(h^{\delta}), \quad M \gg 1,
$$
we find that the symbol of $\Im P_{\eps}$, acting on $H(\Lambda_{\eps G})$ 
is ${\cal O}(\eps)$, and away from an arbitrariliy small but fixed 
\neigh{} of $\cup_{j=1}^N \gamma_j$ we have that $\abs{\Im P_{\eps}(\rho)}\geq \eps/C$ 
when we restrict the attention to the region 
$\abs{\Re{P_{\eps}(\rho)}}\leq 1/C$. 

When working microlocally near $\tau=x=\xi=0$ in $T^*(S^1 \times \real)$ 
and simplifying the operator (\ref{6.47}) further, we use Proposition 
\ref{Prop4.3} 
to 
find a \hol{} \cantransf{} 
$$
\kappa_{\sigma, \eps, \frac{h^2}{\eps}}: {\rm neigh}(0,\comp^2)\rightarrow 
{\rm neigh}(0,\comp^2)
$$
depending analytically on $\sigma \in {\rm neigh}(0,\comp)$ and smoothly 
on $\eps, h^2/\eps \in {\rm neigh}(0,\real)$, such that 
$$
\left(\widetilde{p}_0+ \eps \frac{h^2}{\eps} 
\widetilde{p}_2\right)\left(\sigma, \kappa_{\sigma, \eps, 
\frac{h^2}{\eps}}(y,\eta)\right) 
= f(\sigma)+ i\eps g_{\eps, \frac{h^2}{\eps}}\left(\sigma, 
\frac{y^2+\eta^2}{2}\right).
$$
As before, associated to $\kappa_{\sigma,\eps,\frac{h^2}{\eps}}$, we 
construct an IR-submanifold of $T^*(\widetilde{S^1}\times \comp)$ which is 
$(\eps+h^2/\eps)$-close to $T^*(S^1\times \real)$, and which has the 
property that $\tau$ is real along this submanifold. This leads to a new 
IR-manifold $\Lambda \subset \comp^4$ such that on $\Lambda$, 
$\Im{P}_{\eps}$ has a symbol of modulus $\sim \eps$ in the region 
$\abs{\Re P_{\eps}}< 1/C$, when away from the union 
of small \neigh{}s $\Omega_j$ of $\gamma_j \subset \Lambda$, $1\leq j\leq N$. 
In $\Omega_j$, $P_{\eps}$ is equivalent to an 
operator constructed in section \ref{section4}, which has the form 
$$
f(hD_t)+i\eps G_j\left(hD_t, \frac{x^2+(hD_x)^2}{2},\eps, 
\frac{h^2}{\eps};h\right), 
$$
with 
$$
G_j\left(\tau,q,\eps,\frac{h^2}{\eps};h\right)\sim \sum_{l=1}^{\infty} 
G_{j,l}\left(\tau,q,\eps,\frac{h^2}{\eps}\right)h^l.
$$
Again we see that we have a globally well-posed Grushin problem for 
$P_{\eps}-z$ in the $h$-dependent Hilbert space $H(\Lambda)$. The following result complements Theorem 
\ref{Th6.6}. 

\begin{theo}\label{Th6.7}
Make the assumptions of Theorem {\rm \ref{Th6.6}}, and assume in addition that the 
subprincipal symbol of $P_{\eps=0}$ vanishes. Then for $\eps$ in the range 
$$
h^2 \ll \eps < h^{\delta}, \quad \delta>0,
$$
the eigenvalues of $P_{\eps}$ in the set of the form 
$$ 
\left(-\frac{1}{{\cal O}(1)}, \frac{1}{{\cal O}(1)}\right)+i\eps 
\left(F_0-\frac{1}{{\cal O}(1)},F_0+\frac{1}{{\cal O}(1)}\right)
$$
are given by 
$$
f\left(h\left(k_1-\frac{k_0}{4}\right)-\frac{S}{2\pi}\right) + i\eps 
G_j\left(h\left(k_1-\frac{k_0}{4}\right)-\frac{S}{2\pi}, 
h\left(\frac{1}{2}+k_2\right), \eps, \frac{h^2}{\eps};h\right), 
$$
modulo ${\cal O}(h^{\infty})$, when $1\leq j\leq N$ and $(k_1,k_2)\in \z\times \nat$. Here 
$f(\tau)$ is real-valued with $f(0)=0$ and $f'(0)>0$. The function 
$G_j(\tau,q,\eps,h^2/\eps;h)$ for $1\leq j\leq N$, is analytic in $\tau$ and $q$ in a \neigh{} of 
$(0,0)\in \comp^2$, and smooth in 
$\eps, h^2/\eps\in {\rm neigh}(0,\real)$, 
and has an asymptotic expansion in the space of such functions, as 
$h\rightarrow 0$, 
$$
G_j\left(\tau,q,\eps,\frac{h^2}{\eps};h\right)\sim \sum_{l=0}^{\infty} 
G_{j,l}\left(\tau,q,\eps, \frac{h^2}{\eps}\right) h^l, 
$$
where $\Re G_{j,0}(0,0,0,0)=F_0$ and 
$$
\frac{\partial}{\partial q} \Re G_{j,0}(0,0,0,0)>0.
$$
\end{theo}

\section{Barrier top resonances in the resonant case}\label{section6}
\setcounter{equation}{0}

\vspace*{2mm}
\noindent 
Consider 
\begeq
\label{7.1}
P= -h^2 \Delta + V(x), \quad p(x,\xi)=\xi^2 +V(x), \,\, x,\xi \in \real^2, 
\endeq
and let us assume that $V(x)$ is real-valued, and that it extends 
holomorphically to 
a set $\{x \in \comp^2; \abs{\Im x} \leq \langle{\Re x}\rangle/ C \}$, for 
some $C>0$, and tends to $0$ when 
$x\rightarrow \infty$ in that set. The resonances of $P$ can be defined in 
an angle 
$\{z \in \comp; -2\theta_0 < \arg z<0\}$ for some fixed small 
$\theta_0>0$, as the eigenvalues of
$P\bigg|_{e^{i\theta_0}{\bf R}^2}$ in the same region. 

\vskip 2mm
\noindent 
We shall assume that $V(0)=E_0>0$, $\nabla V(0)=0$, and $V''(0)$ is a 
negative definite quadratic form. Assume also that the 
union of trapped $H_p$-trajectories in $p^{-1}(E_0)\cap \real^4$ is 
reduced to $(0,0)\in \real^4$. (We recall that a trapped 
trajectory is a maximal integral curve of the Hamilton vector field $H_p$, 
contained in a bounded set.) We are then 
interested in resonances of $P$ near $E_0$, created by the critical point 
of $V$. After a linear 
symplectic change of coordinates, and a conjugation of $P$ by means of the 
corresponding metaplectic operator, we may 
assume that as $(x,\xi)\rightarrow 0$,   
\begeq
\label{7.2}
p(x,\xi)-E_0 = \sum_{j=1}^2 \frac{\lambda_j}{2} \left(\xi_j^2 
-x_j^2\right) + p_3(x)+p_4(x)+\ldots,\quad \lambda_j>0. 
\endeq
Here $p_j(x)$ is a homogeneous polynomial of degree $j\geq 3$. 

\vskip 2mm
\noindent
For future reference we recall that according to the theory of resonances 
developed in \cite{HeSj}, the resonances of $P$ in a fixed 
$h$-independent \neigh{} of $E_0$ can also be viewed as the eigenvalues of 
$P: H(\Lambda_{t G},1) \rightarrow H(\Lambda_{tG},1)$, 
equipped with the domain $H(\Lambda_{tG},\langle{\xi}\rangle^2)$. Here $G 
\in C^{\infty}(\real^2; \real)$ is an escape function in 
the sense of \cite{HeSj}, $t>0$ is sufficiently small and fixed, and 
$\Lambda_{tG}$ is a suitable IR-deformation of $\real^4$, associated 
with the function $G$. The Hilbert space $H(\Lambda_{tG},1)$ consists of 
all tempered distributions $u$ such that a suitable 
FBI transform $Tu$ belongs to a certain exponentially weighted 
$L^2$-space. We refer to \cite{HeSj} for the original presentation of the 
microlocal theory of resonances, and to \cite{LahMa} for a simplified version 
of the theory, which is applicable in the present setting of 
operators with globally analytic coefficients, converging to the Laplacian 
at infinity. Here we shall only remark that as in 
\cite{KaKe}, the escape function $G$ can be chosen such that $G=x\cdot \xi$ in 
a \neigh{} of $(0,0)$, and such that $H_pG>0$ on 
$p^{-1}(E_0)\setminus\{(0,0)\}$. 

\vskip 2mm
Under the assumptions above, but without any restriction on the dimension 
and without any assumption on the signature of $V''(0)$,
all resonances in a disc around $E_0$ of radius $Ch$ were determined 
in \cite{Sj1}. Here $C>0$ is arbitrarily large and fixed. 
(See also \cite{BrCoDu}.) Specializing the result of \cite{Sj1} to the present 
barrier top case, we may recall that in this disc, 
the resonances are of the form 
\begeq
\label{7.3}
E_0 - i\left(k_1+\half \right)\lambda_1 h - 
i\left(k_2+\half\right)\lambda_2 h +{\cal O}(h^{3/2}), \quad h \rightarrow 
0, 
\quad k=(k_1,k_2) \in \nat^2.
\endeq
Furthermore, in the non-resonant case, i.e. when 
\begeq
\label{7.4}
\lambda \cdot k \neq 0,\quad 0\neq k \in \z^2,
\endeq
a result of Kaidi and Kerdelhu\'e \cite{KaKe} extended \cite{Sj1} to obtain all 
resonances in a disc around $E_0$ of radius $h^{\delta}$, 
for each fixed $\delta>0$ and $h>0$ small enough depending on $\delta$. In 
this case, the resonances are given by asymptotic expansions in integer 
powers of $h$, with the leading term as in (\ref{7.3}). 

Throughout this section we shall work under the following resonant 
assumption, 
\begeq
\label{7.5}
\lambda \cdot k=0, \,\,\wrtext{for some}\,\, 0 \neq k\in \z^2.
\endeq
In this case we shall show how to obtain a description of all the 
resonances in an energy shell of the form 
$$
h^{4/5} \ll \abs{E-E_0} < {\cal O}(1) h^{\delta}, \quad \delta>0, 
$$
provided that we avoid an arbitrarily small half-cubic \neigh{} of 
$E_0-i[0,\infty)$.

\vskip 2mm
\noindent 
The starting point is a reduction to an eigenvalue problem for a scaled 
operator, as in \cite{KaKe}, \cite{MeSj}, \cite{Sj2}. In these works it 
was shown how to adapt the theory of \cite{HeSj} so that $P$ can be realized as 
an operator acting on a suitable $H(\Lambda)$-space, where 
$\Lambda \subset \comp^4$ is an IR-manifold which coincides with $T^* 
\left(e^{i\pi/4}\real^2\right)$ near $(0,0)$, and further away from 
a \neigh{} of this point, it agrees with $\Lambda_{tG}$. Furthermore, 
$\Lambda$ has the property 
that on this manifold, $p-E_0$ is elliptic away from a small \neigh{} of 
$(0,0)$, and this \neigh{} can be chosen arbitrarily 
small, provided that the constant in the elliptic estimate is taken 
sufficiently large. Using a Grushin reduction exactly as in \cite{MeSj}, we may 
and will therefore reduce the study of resonances of $P$ near $E_0$ to an 
eigenvalue problem for $P$ after the complex scaling, 
which near $(0,0)$ is given by $x=e^{i\pi / 4}\widetilde{x}$, $\xi =e^{-i 
\pi/4}\widetilde{\xi}$, $\widetilde{x},\widetilde{\xi}\in \real$. 

Using (\ref{7.2}) and dropping the tildes from the notation, we see that 
the principal symbol of the scaled operator has the form 
\begeq
\label{7.6}
E_0-i \left(p_2(x,\xi) + i e^{3\pi i/4}p_3(x)+ie^{4i\pi /4}p_4(x)+\ldots 
\right), \quad (x,\xi)\rightarrow 0, 
\endeq 
where 
\begeq
\label{7.7}
p_2(x,\xi)=\sum_{j=1}^2 \frac{\lambda_j}{2} \left(\xi_j^2+x_j^2\right) 
\endeq
is the harmonic oscillator. In what follows we shall therefore consider an 
$h$-pseudo\-diffe\-ren\-tial ope\-rator $P$ on $\real^2$, 
microlocally defined near $(0,0)$, with the leading symbol 
\begeq
\label{7.8}
p(x,\xi)=p_2(x,\xi)+i e^{3\pi i/4}p_3(x)+\ldots, \quad (x,\xi)\rightarrow 
0, 
\endeq
and with the vanishing subprincipal symbol. We extend $P$ to be globally 
defined as a symbol of class 
$S^0(\real^4)=C^{\infty}_{\rm b}(\real^4)$, with the asymptotic expansion
$$
P(x,\xi;h) \sim p(x,\xi) + h^2 p^{(2)}(x,\xi)+ \ldots,
$$
in this space, and so that 
$$
\abs{p(x,\xi)} \geq \frac{1}{C}, \quad C>0, 
$$
outside a small \neigh{} of $(0,0)$. 

We shall be interested in eigenvalues $E$ of $P$ with $\abs{E}\sim 
\eps^2$, $0< \eps \ll 1$. It follows from \cite{Sj4} that the 
corresponding eigenfunctions are concentrated in a region where 
$\abs{(x,\xi)}\sim \eps$, and so we introduce the change of variables 
$x=\eps y$, $h^{\delta}\leq \eps \leq 1$, $0< \delta < 1/2$. Then 
$$
\frac{1}{\eps^2} P(x,hD_x;h)=\frac{1}{\eps^2} 
P(\eps(y,\widetilde{h}D_y);h),\quad \widetilde{h}=\frac{h}{\eps^2} \ll 1. 
$$
The corresponding new symbol is
$$
\frac{1}{\eps^2} P(\eps(y,\eta);h) \sim \frac{1}{\eps^2} 
p(\eps(y,\eta))+\eps^2 \widetilde{h}^2 p^{(2)}(\eps(y,\eta))+\ldots,
$$
to be considered in the region where $\abs{(y,\eta)}\sim 1$. The leading 
symbol becomes 
$$
\frac{1}{\eps^2} p(\eps(y,\eta))=p_2(y,\eta)+i\eps e^{3\pi 
i/4}p_3(y)+{\cal O}(\eps^2), 
$$
for $(y,\eta)$ in a fixed \neigh{} of $(0,0)$. 

Now the resonant assumption (\ref{7.5}) implies that the $H_{p_2}$-flow is 
periodic on $p_2^{-1}(E)$, for 
$E\in {\rm neigh}(1,\real)$, with period $T>0$ which does not depend on 
$E$. For $z\in {\rm neigh}(1,\comp)$, we shall then 
discuss the invertibility of 
$$
1/\eps^2 P(x,hD_x;h)-z
$$ 
in the range of $\eps$, dictated by Theorem \ref{Th6.4}, and using 
$\widetilde{h}$ as the new 
semiclassical parameter. Indeed, all the assumptions of that theorem are 
satisfied in a fixed \neigh{} of $(0,0)$, and outside such a 
\neigh{}, we have ellipticity which guarantees the invertibility there. 

\begin{prop}\label{Prop7.1}
Assume that {\rm (\ref{7.5})} holds. When $p_3$ is a homogeneous 
polynomial of degree $3$ on $\real^2$, we let 
$\langle{p_3}\rangle$ denote the average of $p_3$ along the trajectories 
of the Hamilton vector field of 
$p_2$ in {\rm (\ref{7.7})}, and assume that $\langle{p_3}\rangle$ is not 
identically zero. Let $F_0 \in \real$ be a regular value of 
$\cos(3 \pi /4)\langle{p_3}\rangle$ restricted to $p_2^{-1}(1)$, and 
assume that $T$ is the minimal period of the 
$H_{p_2}$-trajectories in the manifold $\Lambda_{1,F_0}$ given by 
$$
\Lambda_{1,F_0}: p_2 = 1, \cos \left(\frac{3\pi}{4}\right) 
\langle{p_3}\rangle = F_0. 
$$
Assume that $\Lambda_{1,F_0}$ is connected. Let $\eps$ satisfy 
\begeq
\label{7.9}
h^{2/5} \ll \eps = {\cal O}(1) h^{\delta},\quad \delta>0.
\endeq 
Then for $z$ in the set 
\begeq
\label{7.9.1}
\left[1-\frac{1}{{\cal O}(1)}, 1+\frac{1}{{\cal O}(1)}\right] + i \eps 
\left[F_0 - \frac{1}{{\cal O}(1)}, 
F_0 + \frac{1}{{\cal O}(1)}\right], 
\endeq
the operator $\eps^{-2} P(x,hD_x;h)-z: L^2 \rightarrow L^2$ is 
non-invertible precisely when $z=z_k$ for some $k\in \z^2$, where the numbers $z_k$ satisfy  
$$
z_k = \widehat{P}\left(\widetilde{h}(k-\frac{\alpha}{4})-\frac{S}{2\pi}, 
\eps, \frac{\widetilde{h}^2}{\eps}; \widetilde{h}\right)+
{\cal O}(h^{\infty}),\quad \widetilde{h}=\frac{h}{\eps^2}. 
$$
Here $\widehat{P}\left(\xi,\eps,\frac{\widetilde{h}^2}{\eps}; 
\widetilde{h}\right)$ has an expansion as $\widetilde{h}\rightarrow 0$, 
$$
\widehat{P}\left(\xi,\eps,\frac{\widetilde{h}^2}{\eps}; 
\widetilde{h}\right)\sim p_2(\xi_1) + 
\eps \sum_{j=0}^{\infty} \widetilde{h}^j 
r_j\left(\xi,\eps,\frac{\widetilde{h}^2}{\eps}\right), 
$$
where 
$$
r_0 = ie^{3\pi i/4}\langle{p_3}\rangle(\xi)+ {\cal 
O}\left(\eps+\frac{\widetilde{h}^2}{\eps}\right). 
$$
The coordinates $\xi_1=\xi_1(E)$ and $\xi_2=\xi_2(E,F)$ are the 
norma\-lized ac\-tions of 
$$
\Lambda_{E,F}: p_2=E, \cos\left(\frac{3\pi}{4}\right) 
\langle{p_3}\rangle=F,
$$
for $E\in {\rm neigh}(1,\real)$, $F\in {\rm neigh}(F_0,\real)$, given by 
\begeq
\label{7.10}
\xi_j = \frac{1}{2\pi} \left(\int_{\gamma_j(E,F)} 
\eta\,dy-\int_{\gamma_j(1,F_0)} \eta \,dy\right), \quad j=1,2,
\endeq
with $\gamma_j(E,F)$ being fundamental cycles in $\Lambda_{E,F}$, such 
that $\gamma_1(E,F)$ corresponds to a closed $H_{p_2}$-trajectory of 
minimal 
period $T$. Furthermore, 
\begeq
\label{7.11}
S_j =\int_{\gamma_j(1,F_0)} \eta\,dy, \quad j=1,2,\,\,\, S=(S_1,S_2),  
\endeq
and $\alpha \in \z^2$ is fixed.  
\end{prop}

\Remark. In the case when the compact manifold $\Lambda_{1,F_0}$ has finitely many connected 
components $\Lambda_j$, $1 \leq j \leq M$, with each $\Lambda_j$ being diffeomorphic to a torus, 
the set of $z$ in (\ref{7.9.1}) for which the operator $\eps^{-2} P(x,hD_x;h)-z$ is 
non-invertible agrees with the union of the quasi-eigenvalues constructed for each component, 
up to an error which is ${\cal O}(h^{\infty})$. In the following discussion, for simplicity 
it will be tacitly assumed that $\Lambda_{1,F_0}$ is connected.  

\vskip 4mm 
The reduction by complex scaling together with the scaling argument above 
and Proposition \ref{Prop7.1} allows us to describe the 
resonances $E$ of the operator (\ref{7.1}) in the set 
\begeq
\label{7.12}
h^{4/5} \ll \abs{E-E_0} ={\cal O}(1) h^{\delta}, \quad \delta>0, 
\endeq
by 
\begeq
\label{7.13}
E=E_0 -i\eps^2 
\widehat{P}\left(\widetilde{h}\left(k-\frac{\alpha}{4}\right)-\frac{S}{2\pi}, 
\eps, \frac{\widetilde{h}^2}{\eps}; \widetilde{h}\right)+
{\cal O}(h^{\infty}), 
\endeq
where we choose $\eps>0$ with $\abs{E-E_0}/\eps^2 \sim 1$. The description 
(\ref{7.13}) is valid provided that we exclude sets of 
the form 
\begeq
\label{7.14}
E\in \comp, \,\, \abs{\Re E - E_0 - F_0 \abs{\Im E}^{3/2}} < 
\frac{1}{{\cal O}(1)} \abs{\Im E}^{3/2}, 
\endeq
from the domain (\ref{7.12}). Here $F_0$ varies over the set of critical 
values of $\cos(3\pi/4)\langle{p_3}\rangle$ restricted to $p_2^{-1}(1)$. 
Indeed, writing $E=E_0-i\eps^2 z$, we see that the set (\ref{7.14}) in the 
$E$--plane corresponds to the set 
$\abs{\Im z-\eps F_0} < \eps/{\cal O}(1)$ in the $z$--plane. It is also 
clear that when 
$F_0 \in \{ \inf_{p_2^{-1}(1)} \cos(3\pi/4)\langle{p_3}\rangle, 
\sup_{p_2^{-1}(1)} \cos(3\pi/4)\langle{p_3}\rangle\}$, an 
application of Theorem \ref{Th6.7} will allow us to extend a 
description of the resonances to a set of the form (\ref{7.14}), 
provided that the assumptions of that theorem are satisfied. In what follows, 
we shall content ourselves by discussing an explicit example. 

\vskip 2mm
\noindent
Our starting point will be deriving an expression for $\langle{p_3}\rangle$. Consider 
$$
p_2(x,\xi)=\sum_{j=1}^2 \frac{\lambda_j}{2}(x_j+\xi_j^2), \quad \lambda_j>0, 
$$
where the $\lambda_j$ satisfy (\ref{7.5}). In order to describe the $H_{p_2}$-flow, it is convenient to 
introduce the action-angle variables $I_j\ge 0$, $\tau_j\in{\bf R}/2\pi {\bf Z}$ for $p_2$, given by
\begeq
\label{7.15}
x_j=\sqrt{2I _j}\cos \tau _j,\ \xi _j=-\sqrt{2I _j}\sin \tau_j.
\endeq
Then $p_2 =\sum \lambda _j I_j$ and the Hamilton flow is given by
${\bf R}\ni t\mapsto (I(t),\tau(t))$, with $I(t)=I(0)$,
$\tau(t)=\tau(0)+t\lambda $, $\lambda =(\lambda _1,\lambda _2)$. In the 
original coordinates, this gives
\begeq
\label{7.16}
{\cases{x_j(t)=\sqrt{2I_j(0)}\cos (\tau_j(0)+\lambda _jt)\cr
\xi_j(t)=-\sqrt{2I_j(0)}\sin (\tau_j(0)+\lambda _jt)},}
\endeq
and we get a combination of two rotations in $(x_j,\xi _j)$, $j=1,2$, with 
minimal periods $2\pi /\lambda _j$ (except in the degenerate
cases when one of the $(x_j,\xi _j)$ vanishes). Avoiding the totally 
degenerate case when $I=0$, we get trajectories with
\begin{itemize}
\item minimal period $2\pi /\lambda _2$ when $I_1(0)=0$,
\item minimal period $2\pi /\lambda_1$ when $I_2(0)=0$, 
\item minimal period $T=-k_2^02\pi /\lambda _1=k_1^02\pi/ \lambda _2$, 
when both $I_1(0)$ and $I_2(0)$ are $\ne 0$.
\end{itemize}
Here we let $k^0 = (k^0_1,k^0_2)$ be the point satisfying (\ref{7.5}), 
which has minimal norm and positive first component. The integers $k$
in (\ref{7.5}) are equally spaced on the straight line $\lambda^{\perp}$, 
and it will be convenient to represent them in the form 
$nk^0$, $n\in \z\setminus\{0\}$. 

We shall now compute the averages $\langle x^\alpha \rangle $ along the 
$H_{p_2}$-trajectories of a monomial 
$x^\alpha =x_1^{\alpha _1}x_2^{\alpha _2}$. Using (\ref{7.16}), we get
\begin{eqnarray}
\label{7.17}
& & \langle{x^{\alpha}}\rangle = I(0)^{\alpha \over 2}2^{\vert \alpha 
\vert \over 2}{1\over T} 
\int_0^T (\cos (\tau_1(0)+\lambda _1t))^{\alpha _1}(\cos 
(\tau_2(0)+\lambda _2t))^{\alpha _2}dt \\ \nonumber 
& = & \frac{I(0)^{\alpha \over 2}}{2^{\vert \alpha \vert \over 2}} {1\over 
T}
\int_0^T (e^{i(\tau_1(0)+\lambda_1t)}+e^{-i(\tau_1(0)+\lambda 
_1t)})^{\alpha _1}
(e^{i(\tau_2(0)+\lambda_2t)}+e^{-i(\tau_2(0)+\lambda _2t)})^{\alpha _2}dt. 
\end{eqnarray}
Here the integrand can be developed with the binomial theorem, 
$$
\sum_{k_1=0}^{\alpha _1}\sum_{k_2=0}^{\alpha _2}\pmatrix{\alpha_1\cr 
k_1}\pmatrix{\alpha_2\cr k_2}
e^{i((2k_1-\alpha _1)\tau_1(0)+(2k_2-\alpha _2)\tau_2(0))}
e^{i((2k_1-\alpha _1)\lambda _1+(2k_2-\alpha _2)\lambda _2)t}, 
$$
and only the terms with $(2k_1-\alpha _1)\lambda_1+(2k_2-\alpha _2)\lambda 
_2=0$ can give a non-vanishing
contribution to the integral. This means that $2k-\alpha =nk^0$ for some 
$n\in{\bf Z}$, i.e. $\alpha +nk^0=2k$ with 
$0\le k\leq \alpha $ componentwise. We get
\begeq
\label{7.18}
{\langle x^\alpha \rangle ={I(0)^{\alpha /2}\over 2^{\vert\alpha \vert /2}}
\sum_{\alpha +nk_0=2k\atop 0\le k\le \alpha }\pmatrix{\alpha _1\cr 
k_1}\pmatrix{\alpha _2\cr k_2}
\cos ((2k_1-\alpha _1)\tau_1(0)+(2k_2-\alpha _2)\tau_2(0)), }
\endeq
where it is understood that $n\in{\bf Z}$, $k\in{\bf N}^2$, and where we 
notice that if $\alpha +nk_0=2k,\, 0\le k\le \alpha $,
then $\widetilde{k}:=\alpha -k$ also participates in the sum, since $0\le 
\widetilde{k}\le \alpha $ and $\alpha
-nk_0=2\widetilde{k}$. Also notice that the cosine in (\ref{7.18}) can be 
written in the form $\cos (nk_0\cdot \tau(0))$. In order to 
find the 
non-vanishing terms in (\ref{7.18}), we consider the "line"  ${\bf Z}\ni 
n\mapsto \alpha +nk_0\in{\bf Z}^2$. The points on this line 
in the rectangle $([0,2\alpha_1]\times [0,2\alpha _2])\cap{\bf N}^2$ with 
even coordinates correspond to the terms in (\ref{7.18}). 

\smallskip
\par
\noindent 
\it Example 1. \rm Let $k^0=(1,-1)$, correponding for instance to $\lambda 
=(1,1)$. In this case the two components of $\alpha $
must have the same parity.
\par\noindent For $\alpha =(2,0)$ we have only one term with
$n=0$, $k=(1,0)$, and $\langle x_1^2 \rangle=I_1(0) $.
\par\noindent For $\alpha =(0,2)$ we get similarly $\langle x_2^2
\rangle=I_2(0) $.
\par\noindent For $\alpha =(1,1)$ we get two terms with
$n=1,k=(1,0)$ and $n=-1,k=(0,1)$ respectively, and $\langle
x_1x_2\rangle =\sqrt{I_1(0)I_2(0)}\cos (\tau_1(0)-\tau_2(0))$.
\par\noindent For $\vert \alpha \vert =3$ we get no non-vanishing terms.
\par\noindent For $\alpha =(4,0)$ we have one term with $n=0$, $k=(2,0)$ 
and we get $\langle x_1^4\rangle ={3\over 2}I_1(0)^2$.
\par\noindent For $\alpha =(0,4)$ we get similarly, $\langle x_2^4\rangle 
={3\over 2}I_2(0)^2$.
\par\noindent For $\alpha =(2,2)$ we get one term with $n=2,k=(2,0)$ and 
one with $n=-2,k=(0,2)$, We also have a term
with $n=0,k=(1,1)$, and this leads to $\langle x_1^2x_2^2\rangle 
=I_1(0)I_2(0)(1+{1\over 2}\cos 2(\tau_1(0)-\tau_2(0)))$.

\medskip 
\noindent 
It follows from Example 1 that Proposition 7.1 does not apply when 
$\lambda={\rm Const.}(1,1)$, since in this case 
$\langle{p_3}\rangle \equiv 0$. We shall therefore consider a different 
choice of the resonant frequencies. 

\vskip 2mm
\noindent 
\it Example 2. \rm Let us take $k^0=(2,-1)$, corresponding for instance to 
$\lambda =(1,2)$, and let 
$\vert \alpha \vert =3$. For $\alpha =(3,0),(0,3),(1,2)$ it follows from 
(\ref{7.18}) that $\langle x^\alpha \rangle =0$.
For $\alpha =(2,1)$ we get two terms, one with $n=1,k=(2,0)$ and one with 
$n=-1,k=(0,1)$. It follows that 
\begeq
\label{7.19}
\langle x_1^2x_2\rangle=2^{-1/2}I_1(0)I_2(0)^{1/2}\cos 
(2\tau_1(0)-\tau_2(0)).
\endeq

\medskip 

For future reference, we shall also describe how the averages 
$\langle{x^{\alpha}}\rangle$ can be computed after a suitable 
complex linear change of symplectic coordinates. Introduce 
$$
\cases{y={1\over \sqrt{2}}(x-i\xi )\cr \eta ={1\over i\sqrt{2}}(x+i\xi 
)},\ \cases{x={1\over \sqrt{2}}(y+i\eta )\cr
\xi ={i\over \sqrt{2}}(y-i\eta )}.
$$
In these coordinates $p=\sum_{j=1}^2 i\lambda _jy_j\eta _j$, and 
$$
\exp (tH_p)(y,\eta )=(e^{it\lambda_1 
}y_1,e^{it\lambda_2}y_2,e^{-it\lambda_1 }\eta _1,e^{-it\lambda_2 }\eta _2),
$$
so that 
$$
\langle y^\alpha \eta ^\beta \rangle ={1\over T}\int_0^T e^{i\lambda \cdot 
(\alpha -\beta )t}dt y^\alpha \eta ^\beta
=\cases{y^\alpha \eta ^\beta \hbox{ if }\lambda \cdot (\alpha -
\beta)=0,\cr 0\hbox{ otherwise}. }
$$
We apply this to
$$
x^\alpha ={1\over 2^{\vert \alpha \vert/2}}\sum_{0\le k\le \alpha }
\pmatrix{\alpha \cr k}y^k(i\eta)^{\alpha -k},
$$ 
and get
\begeq
\label{7.20}
{\langle x^\alpha \rangle =2^{-\vert \alpha \vert 
}\sum_{\alpha+nk_0=2k\atop 0\le k\le\alpha }
\pmatrix{\alpha \cr k}(x-i\xi)^k(x+i\xi )^{\alpha -k}.}
\endeq
As before we check that for each term present there is also the
complex conjugate. 

The computations of Examples 1 and 2 can be written like (\ref{7.20}). We 
shall only do it for the last example with 
$k^0=(2,-1),\alpha =(2,1)$:
\begeq
\label{7.21}
{\langle x_1^2x_2\rangle ={1\over 4}\Re ((x_1+i\xi _1)^2(x_2-i\xi
_2))={1\over 4}(x_1^2x_2+2x_1\xi _1\xi _2-x_2\xi _1^2).}
\endeq
We may assume that $\lambda =(1,2)$, so that
\begeq
\label{7.22}
p_2={1\over 2}(x_1^2+\xi _1^2)+(x_2^2+\xi _2^2),
\endeq
and we may then check directly that $H_{p_2}\langle x_1^2x_2\rangle =0$.

From (\ref{7.21}) and (\ref{7.22}) it is clear that $dp_2$ and 
$d\langle x_1^2x_2\rangle $ are linearly \indep{} 
except on some set of measure 0. When computing the critical points of 
$\langle{x_1^2x_2}\rangle$ on $p_2^{-1}(1)$, we shall first make use 
of the $(I,\tau)$-coordinates. From (\ref{7.19}) we recall that 
\begeq
\label{7.23}
{
p_2=I_1+2I_2,\ \sqrt{2}\langle x_1^2x_2\rangle =I_1 I_2^{1\over 2}\cos 
(2\tau_1-\tau_2). }
\endeq
It follows from the Hamilton equations that $\theta :=2\tau_1-\tau_2$ is 
invariant under the $H_{p_2}$-flow, and we can therefore 
work in the coordinates $I_1,I_2,\theta $. We have 
\begeq
\label{7.24}
{dp_2=dI_1+2dI _2,} {\sqrt{2}d\langle x_1^2x_2\rangle =(I_2^{1\over 2}\cos 
\theta )dI_1+{1\over 2}I_1I_2^{-{1\over2}}
(\cos \theta )d I_2 -I_1 I_2^{1\over 2}(\sin \theta)d\theta .}
\endeq
If $\theta \not\in \pi {\bf Z}$, $I_1,I_2\ne 0$, we have $\partial _\theta 
\langle x_1^2x_2\rangle \ne 0$, and hence the
differentials are linearly \indep{}. Still with $I_1,I_2\ne 0$, let 
$\theta \in \pi {\bf Z}$, so that $\cos \theta =\pm 1$. Then 
the differentials are linearly dependent iff 
$$
0=\det \pmatrix{1 & 2\cr I_2^{1/2} &{1\over 2}I_1 I_2^{-{1\over 2}}},
$$
i.e. iff 
$$
I_1=4 I_2.
$$
This gives two closed trajectories inside the energy surface $p_2=1$ and 
the corresponding values for $\langle x_1^2x_2\rangle $:
\begeq
\label{7.25}
I_1={2\over 3},\,I_2={1\over 6},\ 2\tau_1-\tau_2=0;\ \ \langle 
x_1^2x_2\rangle ={1\over 3\sqrt{3}}, 
\endeq
and 
\begeq 
\label{7.26}
I_1={2\over 3},\,I_2={1\over 6},\ 2\tau_1-\tau_2=\pi ;\ \ \langle 
x_1^2x_2\rangle ={-1\over 3\sqrt{3}}.
\endeq

When $I_1=0$ or $I_2=0$, the question of linear independence of the 
differentials should be analyzed directly
in the $(x,\xi )$-coordinates (or $(y,\eta )$-coordinates), and here we 
shall use (\ref{7.21}). On the plane $I_1=0$, corresponding to 
$x_1=\xi _1=0$, we have $d\langle x_1^2x_2\rangle =0$, so here we have 
linear dependence, with the corresponding critical value 
$\langle x_1^2x_2 \rangle =0$. On the plane $I_2=0$, corresponding to 
$x_2=\xi _2=0$, we have
$$
\cases{ d\langle x_1^2x_2\rangle ={1\over 4}(x_1^2-\xi
_1^2)dx_2+{1\over 2}x_1\xi _1d\xi _2,\cr dp_2=x_1dx_1+\xi _1d\xi
_1,}
$$ 
and these differentials are \indep{}, since we avoid the point $x=\xi =0$.

\vskip 2mm 
\noindent
We shall now look at the nature of the critical points of $\langle 
x_1^2x_2\rangle $, when viewed as a \fu{} on
$\Sigma:= p_2^{-1}(1)/\exp ({\bf R}H_p)$. For the trajectories found in 
(\ref{7.25}) and (\ref{7.26}), we use $\theta$ and $I_2$ as 
local coordinates on $\Sigma$, and using (\ref{7.23}) together with 
$I_1=1-2I_2$, we get for $\theta =k\pi $, $k=0,1$, $I_2=1/6$ and 
$f=\sqrt{2}\langle x_1^2x_2\rangle $, 
$$
\partial_\theta \partial _{I_2} f=0,\ \partial_\theta^2 
f=-(1-2I_2)I_2^{1\over 2}(-1)^k,\ 
\partial_{I_2}^2f=-(-1)^k \left(\frac{1}{4} I_2^{-{3\over 2}}+{3\over 
2}I_2^{-{1\over 2}}\right).
$$
For $k=0$ we therefore have a \nondeg{} maximum and for $k=1$ we get a 
\nondeg{} minimum.

For the third trajectory, given by 
\begeq
\label{7.27}
x_1=\xi _1=0, \quad x_2^2+\xi _2^2=1, 
\endeq 
we use that $\langle x_1^2x_1\rangle$ vanishes to the second order there, 
and hence that the transversal Hessian in 
$p_2^{-1}(1)$ can be identified with the free Hessian \wrt{} $x_1,\xi _1$, 
which is given by the matrix
$$
{1\over 2}\pmatrix{x_2 &\xi _2\cr \xi _2 & -x_2}.
$$
The \ev{}s are ${1\over 2}$ and $-{1\over 2}$. Thus we have a \nondeg{} 
saddle point.

\vskip 2mm
\noindent
We summarize the discussion above in the following proposition. 

\begin{prop}\label{Prop7.2} 
Let 
$$
p_2(x,\xi)=\frac{1}{2}\left(x_1^2+\xi_1^2\right)+(x_2^2+\xi_2^2).
$$ 
Then the $H_{p_2}$-flow is periodic in $p^{-1}_2(E)$, for $E\in {\rm 
neigh}(1,\real)$, with period $T=2\pi$. If 
$$
p_3(x)=a_{3,0}x_1^3+a_{1,2}x_1x_2^2+x_1^2x_2+a_{0,3}x_2^3, 
$$
then we have 
$$
\langle{p_3}\rangle(x,\xi)=\frac{1}{4} 
\left(x_1^2x_2+2x_1\xi_1\xi_2-x_2\xi_1^2\right).
$$
The differential of $\langle{p_3}\rangle$, restricted to $p^{-1}_2(1)$, 
vanishes along three closed 
$H_{p_2}$-trajecto\-ries, given by {\rm (\ref{7.25})}, {\rm (\ref{7.26})}, 
and {\rm (\ref{7.27})}. These critical trajec\-tories are non-de\-ge\-ne\-rate 
in the 
sense that the transversal Hessian of $\langle{p_3}\rangle$ is 
non-degenerate. The set of the critical values of $\langle{p_3}\rangle$ is 
$\{\pm (3\sqrt{3})^{-1},0\}$, and the maximum and the minimum of 
$\langle{p_3}\rangle$ are attained along the trajectories {\rm 
(\ref{7.25})} and 
{\rm (\ref{7.26})}, respectively. The transversal Hessian of 
$\langle{p_3}\rangle$ along {\rm (\ref{7.27})} has the signature 
$(1,-1)$. 
The minimal period of the trajectories in {\rm (\ref{7.25})} and 
{\rm (\ref{7.26})} is equal to $T=2\pi$, and the minimal period in 
{\rm (\ref{7.27})} is $\pi$. Let finally $F_0$ be a regular value of 
$\langle{p_3}\rangle$ restricted to 
$p^{-1}_2(1)$. Then the minimal period of every closed 
$H_{p_2}$-trajectory in the Lagrangian manifold 
$$
\Lambda_{1,F_0}: p_2=1, \langle{p_3}\rangle=F_0
$$
is equal to $T=2\pi$. 
\end{prop}

We now return to the operator $P$ with principal symbol $p$ in 
(\ref{7.1}). Under the general assumptions from the beginning of this 
section, we shall assume that as $(x,\xi)\rightarrow 0$, we have 
$$
p(x,\xi)-E_0 = \frac{1}{2}(\xi_1^2-x_1^2) + (\xi_2^2-x_2^2)+p_3(x) +{\cal 
O}(x^4),
$$
where 
$$
p_3(x)=a_{3,0}x_1^3+a_{1,2}x_1x_2^2+x_1^2x_2+a_{0,3}x_2^3.
$$
Let us write $A_1=-(3\sqrt{6})^{-1}$, $A_2=(3\sqrt{6})^{-1}$, and $A_3=0$. 

\begin{prop}\label{Prop7.3}
The resonances of $P$ in the domain 
\begeq
\label{7.28}
\bigl\{ z\in \comp;\,\, h^{4/5}\ll \abs{z-E_0} = {\cal O}(1) h^{\delta} 
\bigr\} \setminus \bigcup_{j=1}^3 
\{z; \abs{\Re z-E_0-A_j \abs{\Im z}^{3/2}}< \eta \abs{\Im z}^{3/2}\},
\endeq 
where $\delta$, $\eta>0$ are arbitrary but fixed, are given by 
\begeq
\label{7.29}
\sim 
E_0-i \left(h(k_1-\alpha_1/4)+\eps^3 \sum_{j=0}^{\infty} h^j \eps^{-2j} 
r_j\left(\frac{h}{\eps^2}\left(k-\frac{\alpha}{4}\right)-
\frac{S}{2\pi},\eps,\frac{h^2}{\eps^5}\right)\right),
\endeq
with 
$$
r_0\left(\xi, \eps, \frac{h^2}{\eps^5}\right)=ie^{3\pi 
i/4}\langle{p_3}\rangle(\xi)+{\cal O}\left(\eps+\frac{h^2}{\eps^5}\right),
$$
$$
r_j\left(\xi, \eps, \frac{h^2}{\eps^5}\right)={\cal 
O}\left(\eps+\frac{h^2}{\eps^5}\right), \quad j \geq 1
$$
analytic in $\xi\in {\rm neigh}(0,\comp^2)$, and smooth in $\eps, h^2/\eps 
\in {\rm neigh}(0,\real)$. We have 
$k=(k_1,k_2)\in \z^2$, $S=(S_1,S_2)$ with $S_1=2\pi$, and 
$\alpha=(\alpha_1, \alpha_2)\in \z^2$ is fixed, and we
choose $\eps>0$ with $\abs{E-E_0}\sim \eps^2$. The resonances in the set 
\begeq
\label{7.30}
\bigl \{z\in \comp,\,\, \abs{\Re z-E_0-A_1 \abs{\Im z}^{3/2}}< \eta 
\abs{\Im z}^{3/2}\bigr \}\,\, \wrtext{and}\,\, 
h^{4/5}\ll \abs{z-E_0} = {\cal O}(1) h^{\delta},
\endeq
are given by $E_0$ plus 
\begeq
\label{7.31}
\frac{1}{i}\left(h\left(k_1-\frac{\alpha_1}{4}\right)+i\eps^3 
\sum_{j=0}^{\infty} h^j \eps^{-2j} 
G_j\left(\frac{h}{\eps^2}\left(k_1-\frac{\alpha_1}{4}\right)-1,
\frac{h}{\eps^2}\left(k_2+\frac{1}{2}\right), 
\eps,\frac{h^2}{\eps^5}\right)\right),
\endeq
with $(k_1,k_2)\in \z \times \nat$, $\alpha_1 \in \z$, and 
$\abs{E-E_0}\sim \eps^2$. The function $G_0(\tau,q,\eps,h^2/\eps^5)$ 
is such that $\Re G(0,0,0,0)=A_1$ and $\frac{\partial}{\partial q}\Re 
G_0(0,0,0,0)>0$. 
An analogous description of resonances is valid in the domain {\rm 
(\ref{7.30})} with $A_1$ replaced by $A_2$.
\end{prop}

Here in (\ref{7.29}) we have also used that when expressed in terms of 
the action coordinates from (\ref{7.10}), it is 
true that $p_2(\xi_1)=\xi_1+1$. 

\Remark. If we replace $r_j(\xi,\eps, h^2/\eps^5)$ in (\ref{7.29}) by 
$r_j(\xi+S/2\pi, \eps, h^2/\eps^5)$, then we get 
$$ 
\sim E_0-i \left(h(k_1-\alpha_1/4)+\eps^3 \sum_{j=0}^{\infty} h^j 
\eps^{-2j} 
r_j\left(\frac{h}{\eps^2}\left(k-\frac{\alpha}{4}\right),\eps,\frac{h^2}{\eps^5}\right)\right).
$$
Now let us notice that the choice of $\eps$ is not unique, and replacing 
$\eps$ by $\lambda \eps$, with $\lambda \sim 1$,
does not affect the resonances. It follows therefore that 
\begeq
\label{7.32}
r_j(\xi,\eps,\tau)=\lambda^{3-2j} r_j\left(\frac{\xi}{\lambda^2}, \lambda 
\eps, \frac{\tau}{\lambda^5}\right).
\endeq
Using this, we define 
$$
r_j(\xi,1,\tau)=\eps^{3-2j} r_j\left(\frac{\xi}{\eps^2},\eps, 
\frac{\tau}{\eps^5}\right),
$$
when $\abs{\xi}\sim \eps^2$ and $\abs{\tau} \leq {\cal O}(\eps^5)$. Then 
(\ref{7.29}) becomes 
$$
\sim E_0-i \left(h(k_1-\alpha_1/4)+\sum_{j=0}^{\infty} h^j 
r_j\left(h\left(k-\frac{\alpha}{4}\right),1,h^2\right)\right).
$$

\appendix
\section{Function spaces and FBI-trans\-forms on ma\-ni\-folds}
\setcounter{equation}{0}
\par Let $X$ be a compact analytic \mfld{} of dimension $n$. In
this section we first review some parts of section 1 in
\cite{Sj5} about how to define global FBI-transforms on $X$,
and function spaces associated to certain IR-deformations of the
real cotangent space. After that we shall perform  Bargmann type
\tf{}s which allow us to view the above mentioned function spaces,
microlocally in a bounded frequency region, as weighted spaces of
\hol{} \fu{}s. The theory in \cite{Sj5} is an adaptation to the case
of compact \mfld{}s of the one in \cite{HeSj} and this as well as
the Bargmann \tf{} below are closely related to similar ideas and
techniques, developed in \cite{BoSj}, \cite{Bo}, \cite{Sj6}, \cite{Zw},
\cite{GoLe}.

\par We equip $X$ with some analytic Riemannian metric so that we
have a distance $d$ and a volume density $dy$. Let $\phi (\alpha
,y)$ be an analytic \fu{} on $\{ (\alpha ,y)\in T^*X\times X;\,
d(\alpha _x,y)< 1/C\}$ (using the notation $\alpha =(\alpha
_x,\alpha _\xi )$, $\alpha _x\in X$, $\alpha _\xi \in T^*_{\alpha
_x}X$) with the following two properties (A) and (B):
\medskip
\par\noindent (A) $\phi $ has a \hol{} extension to a domain of
the form
\begeq
\label{1}
{\{ (\alpha ,y)\in T^*\widetilde{X}\times \widetilde{X};\, \vert
\Im
\alpha _x\vert ,\vert \Im y\vert  <{1\over C},\, \vert \Re \alpha _x -\Re y\vert <{1\over C},\,
\vert \Im \alpha _\xi \vert <{1\over C}\vert \langle \alpha _\xi
\rangle \vert \}}
\endeq 
and satisfies $\vert \phi \vert \le {\cal O}(1)\vert \langle
\alpha _\xi \rangle \vert $ there.

\medskip

\par Here $\widetilde{X}$ is some complexification of $X$ and
$T^*\widetilde{X}$ denotes the cotangent space in the sense of
complex \mfld{}s with pointwise fiber spanned by the pointwise
(1,0)-forms.  We write $\langle \alpha _\xi \rangle =\sqrt{1+\alpha
_\xi ^2}$ with $\alpha _\xi ^2$ defined by means of the dual
metric, and as below, we shall often give statements in local
coordinates whenever convenient and leave to the reader to check
that the statements make sense globally. Notice that by the
Cauchy inequalities,
\ekv{2}
{\partial _{\alpha _x}^k\partial _{\alpha _\xi }^\ell \partial
_y^m\phi ={\cal O}_{k,\ell ,m}(1)\vert \langle \alpha _\xi
\rangle \vert ^{1-\vert \ell \vert },}
in a set of the form (\ref{1}), with a slightly increased constant $C$.

\par The second assumption is \medskip
\par\noindent (B) $\phi (\alpha ,\alpha _x)=0$, $(\partial _y\phi
)(\alpha ,\alpha _x)=-\alpha _\xi $, $\Im (\partial _y^2\phi
)(\alpha ,\alpha _x)\sim \vert \langle \Re \alpha _\xi \rangle
\vert I$.
\medskip

\par By Taylor's formula, we have 
\ekv{3}
{\phi (\alpha ,y)=\alpha _\xi \cdot (\alpha _x-y)+{\cal
O}(1)\langle \alpha _\xi \rangle \vert \alpha _x-y\vert ^2,}
and on the real domain, for $d(\alpha _x,y)\le 1/C$, with $C$
\sufly{} large, we have:
\ekv{4}
{\Im \phi (\alpha ,y)\sim \langle \alpha _\xi \rangle (\alpha
_x-y)^2.}

\par The following example was found in a joint discussion with M.
Zworski: Let $\exp _x:T_xX\to X$ be the geodesic exponential map.
Then we can take 
\ekv{5}
{\phi (\alpha ,y)=-\alpha _\xi \cdot \exp_{\alpha
_x}^{-1}(y)+{i\over 2}\langle \alpha _\xi \rangle d(\alpha
_x,y)^2.}

\par Let $\Lambda \subset T^*\widetilde{X}$ be
a closed I-Lagrangian \mfld{} which is close to $T^*X$ in the
$C^\infty $-sense and which coincides with this set outside a
compact set. Recall that "I-Lagrangian" means Lagrangian for the
real symplectic form $-\Im \sigma $, where $\sigma =\sum d\alpha
_{\xi _j}\wedge d\alpha _{x_j}$ is the standard complex
symplectic form. This means that if we choose (analytic)
coordinates $y$ in $X$ and let $(y,\eta )$ be the corresponding
canonical coordinates on $T^*X$ and $T^*\widetilde{X}$, then
$\Lambda $
is of the form $\{ (y,\eta )+iH_G(y,\eta );\, (y,\eta )\in
T^*X\}$ for some real-valued smooth \fu{} $G(y,\eta )$ which is
close to 0
 in the $C^\infty $-sense and has compact support in $\eta $.
Here $H_G$ denotes the Hamilton field of $G$. Since $\Lambda $ is
close to $T^*X$, it is also R-symplectic in the sense that the
restriction to $\Lambda $ of $\Re \sigma $ is non--degenerate. (We
say that $\Lambda $ is an IR-\mfld.) It follows that 
$${{d\alpha _\vert}_{\Lambda }={d\alpha
_{x_1}\wedge..\wedge d\alpha _{x_n}\wedge d\alpha _{\xi _1}\wedge
..\wedge d\alpha _{\xi _n}} _\vert}_{\Lambda }={1\over
n!}{{\sigma ^n}_\vert}_{\Lambda }$$
is a real non-vanishing $2n$-form on $\Lambda $, that we view as
a positive density. 

\par We also need some symbol classes. A smooth \fu{} $a(x,\xi
;h)$, defined on $\Lambda $ or on a suitable \neigh{} of $T^*X$ in
$T^*\widetilde{X}$ is said to be of class $S^{m,k}$, if 
\ekv{6}
{\partial _x^p\partial _\xi ^qa={\cal O}(1) h^{-m}\langle \xi
\rangle ^{k-q}.}

\par A formal classical symbol $a\in S_{{\rm cl}}^{m,k}$
is of the form $a\sim h^{-m}(a_0+ha_1+...)$ where $a_j\in
S^{0,k-j}$ is \indep{} of $h$. Here and in the following, we let
$0<h\le h_0$ for some \sufly{} small $h_0>0$. When the domain of
definition is real or equal to $\Lambda $, we can find a
realization of $a$ in $S^{m,k}$ (denoted by the same letter $a$)
so that 
$$a-h^{-m}\sum _0^Nh^ja_j\in S^{-(N+1)+m,k-(N+1)}.
$$
When the domain of definition is a complex domain, we say that
$a\in S^{m,k}_{{\rm cl}}$ is a formal classical analytic symbol
($a\in S_{{\rm cla}}^{m,k}$) if $a_j$ are \hol{} and satisfy
\ekv{7}
{\vert a_j\vert \le C_0C^j(j!)\vert \langle \xi \rangle \vert
^{k-j}.}
It is then standard, that we can find a realization $a\in S^{m,k}$
(denoted by the same letter $a$) such that
\begin{eqnarray}
\label{8}
& & \partial _x^k\partial _\xi ^\ell \overline{\partial }_{x,\xi }
a={\cal O}_{k,\ell}(1)e^{-\vert \langle \xi \rangle \vert/{Ch}}, \\ \nonumber 
& & \vert a-h^{-m}\sum_{0\le j\le \vert \langle \xi \rangle \vert
/C_0h}\vert \le {\cal O}(1)e^{-\vert \langle \xi \vert \rangle
/C_1h},
\end{eqnarray}
where in the last estimate $C_0>0$ is \sufly{} large and
$C,C_1>0$ depend on $C_0$. We will denote by $S_{{\rm cl}}^{m,k}$
and $S_{{\rm cla}}^{m,k}$
also the classes of realizations of classical symbols.
We say that a classical (analytic) symbol $a\sim
h^{-m}(a_0+ha_1+...)$ is elliptic, if $a_0$ is elliptic, so that
$a_0^{-1}\in S^{0,-k}$. Take such an elliptic $a(\alpha ,y;h)\in
S_{{\rm cla}}^{{3n\over 4},{n\over 4}}$ and put 
\ekv{9}
{Tu(\alpha ;h)=\int e^{{i\over h}\phi (\alpha ,y)}a(\alpha
,y;h)\chi (\alpha _x,y)u(y)dy,}
where $\chi $
is smooth with support close to the diagonal and equal to 1 in a
\neigh{} of the same set.

\par According to \cite{Sj5} there exists
$b(\alpha ,x;h)\in S^{{3n\over 4},{n\over 4}}_{{\rm cla}}$, such
that if 
\ekv{10}
{
Sv(x)=\int_{T^*X}e^{-{i\over h}\phi ^*(x,\alpha )}b(\alpha
,x;h)\chi (\alpha _x,x)v(\alpha )d\alpha , }
then
\ekv{11}
{STu=u+Ru,}
where $R$ has a distribution kernel $R(x,y;h)$ satisfying 
\ekv{12}
{
\vert \partial _x^\alpha \partial _y^\ell R\vert \le
C_{k,\ell}e^{-{1\over C_0 h}}. }
Here we denote in general by $f^*$, the \hol{} extension of the
complex conjugate of $f$.

\par With $\Lambda $ as above, we put
\ekv{13}
{T_\Lambda u={Tu_\vert}_{\Lambda },} 
and define $S_\Lambda v$ by (\ref{10}), but with $T^*X$ replaced by
$\Lambda $. Then,
\ekv{14}
{S_\Lambda T_\Lambda u=u+R_\Lambda u,}
where $R_\Lambda $ satisfies (\ref{12}) (with a slightly larger $C_0$
and under the assumption that $\Lambda $ is \sufly{} close to
$T^*X$). In fact, using Stokes' formula and the exponential
decrease of $\overline{\partial }$ of the symbols involved, we see
that $S_\Lambda T_\Lambda $ coincides up to an exponentially
small error with $ST$.

\par Since $\Lambda $ is I-Lagrangian, we can find locally a
real-valued smooth function $H(\alpha )$ on $\Lambda $, such that 
\ekv{15}
{
dH=-\Im {(\alpha _\xi \cdot d\alpha _x)_\vert }_\Lambda .
}
Indeed, $-\Im (\alpha _\xi \cdot d\alpha _x)$ is a primitive of
$-\Im \sigma $
and the latter vanishes on $\Lambda $, so the \rhs{} of (\ref{15}) is
closed.

\par We assume:
\ekv{16}
{\hbox{The \e{} (\ref{15}) has a global solution }H\in C^\infty
(\Lambda ;{\bf R}).}
Notice that this property is equivalent to 
\ekv{17}
{
\Im \int_\gamma  (\alpha _\xi \cdot d\alpha _x)=0,\hbox{ for all
closed curves }\gamma \subset\Lambda . }
When (\ref{16}) is fulfilled, $H$ is well--defined up to a constant, and
we shall always choose $H$ to be zero for large $\alpha _\xi $.

\par As in \cite{Sj5} we notice that (\ref{16}) is fulfilled in the
case of IR-\mfld{}s generated by a weight $G\in C^\infty
(T^*\widetilde{X};{\bf R})$ in the following way: Let
$H_G=H_G^{\Im \sigma }$ be the Hamilton field of $G$ \wrt{} $\Im
\sigma $, and assume that $G=0$ in the region where $\vert \alpha
_\xi \vert $ is large. Then for $t$ real with $\vert t\vert $
small enough, we can consider the IR-\mfld{} $\Lambda _t=\exp
(tH_G)(\Lambda _0)$, where $\Lambda _0=T^*X$. Then we get (\ref{16})
with $H=H_t$ given by 
\ekv{18}
{
H_t=\int_0^t (\exp (s-t)H_G)^*(G+\langle H_G,\omega \rangle)ds, }
where $\omega =-\Im (\alpha _\xi \cdot d\alpha _x)$

\par The \fu{} $H$ appears naturally in connection with 
$T_\Lambda $. We have $d_\alpha \phi =\alpha _\xi \cdot d\alpha
_x+{\cal O}(\vert \alpha _x-y\vert )$, so $(d_\alpha \phi )(\alpha
,\alpha _x)=\alpha _\xi \cdot d\alpha _x$ and 
\ekv{19}
{
{-\Im (d_\alpha \phi )(\alpha ,\alpha _x)_\vert }_\Lambda
=d_\alpha H. }
\smallskip
\par\noindent \it Definition. \rm For $m\in{\bf R}$, put
\ekv{20}
{
H(\Lambda ;\langle \alpha _\xi \rangle ^m)=\{ u\in {\cal
D}'(X);\, T_\Lambda u\in L^2(\Lambda ;e^{-2H/h}\vert \langle
\alpha _\xi \rangle \vert ^{2m}d\alpha )\} . }
\smallskip

\par When $\Lambda =T^*X$ we get the usual $h$-Sobolev spaces,
and in particular the case $m=0$ just gives $L^2(X)$. For general
$\Lambda $ we get the same spaces, but the equivalence of the norm 
\ekv{21}
{\Vert u\Vert _{H(\Lambda ,\langle \alpha _\xi \rangle ^m)}=\Vert
T_\Lambda u\Vert _{L^2(\Lambda ;e^{-2H/h}\vert \langle \alpha
_\xi \rangle \vert ^{2m}d\alpha )}}
with the $h$-$m$--Sobolev norm $\Vert u\Vert _{H(T^*X,\langle
\alpha _\xi \rangle ^m)}$ is no longer uniform \wrt{} $h$, in
general.

\par Recall from \cite{Sj5} that if we choose another
FBI-\tf{} $\widetilde{T}$ of the same type as $T$ but with
different phase $\widetilde{\phi }$ and amplitude
$\widetilde{a}$, then for $\Lambda $ close enough to $T^*X$, the
definition (\ref{20}) does not change if we replace $T$ by
$\widetilde{T}$, and we get a new norm which is equivalent to the 
previous one, uniformly \wrt{} $h$. This follows from a fairly 
explicit description of $\widetilde{T}_\Lambda T_\Lambda ^{-1}$.

\par We also know that $Tu=T_{T^*X}u$ and $T_\Lambda u$ satisfy
compatibility conditions similar to the Cauchy-Riemann \e{}s for
\hol{} \fu{}s. For the analysis in the most interesting region
where $\xi $ is \bdd{}, it will be convenient to work with \tf{}s
which are
\hol{} up to exponentially small errors, and for that we make a
different choice of $T$, and take an FBI-\tf{} as in
\cite{Sj6}, now with a global choice of phase (cf \cite{Bo}, 
\cite{GoLe}, \cite{Zw}).
\par The \fu{} $d(x,y)^2$ is analytic in a \neigh{} of the
diagonal in $X\times X$, so we can consider it as a \hol{}
\fu{} in a region 
$$\{ (x,y)\in\widetilde{X}\times \widetilde{X};\, {\rm
dist\,}(x,y)<{1\over C},\, \vert \Im x\vert ,\vert \Im y\vert
<{1\over C}\}.$$
Put 
\ekv{22}
{\phi (x,y)=i\lambda d(x,y)^2,}
where $\lambda >0$ is a constant that we choose large enough,
depending on the size of the \neigh{} of the zero section in
$T^*X$, that we wish to cover. 

\par For $x\in \widetilde{X}$, $\vert \Im x\vert <1/C$, put 
\ekv{23}
{{\cal T}u(x;h)=h^{-{3n\over 4}}\int e^{{i\over h}\phi (x,y)}\chi
(x,y) u(y)dy,\ u\in {\cal D}'(X),}
where $\chi $ is a smooth cut-off function with support in $\{
(x,y)\in \widetilde{X}\times X;\, \vert \Im x\vert <1/C,\,
d(y,y(x))<1/C\}$. Here $y(x)\in X$ is the point close to $x$,
where $X\ni y\mapsto -\Im \phi (x,y)$ attains its non--degenerate
maximum. We have the following facts (\cite{Sj6}):\smallskip

\par The function $\Phi _0(x)=-\Im \phi (x,y(x))$, $x\in
\widetilde{X}$, $\vert \Im x\vert <1/C$, is \stpsh{} and is of the
order of magnitude $\sim \vert \Im x\vert ^2$.\smallskip

\par $\Lambda _{\Phi _0}:=\{ (x,{2\over i}\partial \Phi _0)\in
T^*\widetilde{X}\}$ is an IR-\mfld{} given by $\Lambda _{\Phi
_0}=\kappa _{\cal T}(T^*X)$, where $\kappa _{{\cal T}}$
is the complex canonical \tf{} associated to ${\cal T}$, given by
$(y,-\phi _y'(x,y))\mapsto (x,\phi '_x(x,y))$. Here and in the 
following, we identify $\widetilde{X}$ with its intersection with a
tubular \neigh{} of
$X$ which is \indep{} of the choice of $\lambda $ in
(\ref{22}).\smallskip

\par If $L^2_{\Phi _0}=L^2(\widetilde{X};e^{-2\Phi _0/h}L(dx))$,
for $L(dx)$ denoting a choice of Lebesgue measure (up to a
non-vanishing continuous factor), then ${\cal T}={\cal
O}(1):L^2(X)\to L^2_{\Phi _0}$, $\overline{\partial }_x{\cal
T}={\cal O}(e^{-1/Ch}):L^2(X)\to L^2_{\Phi _0}$. This means that
up to an exponentially small error ${\cal T}u$ is \hol{}
for $u\in L^2(X)$ (and even for $u\in{\cal D}'(X)$). A natural
choice of Lebesgue measure might be $(n!)^{-1}\vert \pi
_*({\sigma _\vert}_{\Lambda _{\Phi _0}})^n\vert $, where $\pi
:\Lambda _{\Phi _0}\to \widetilde{X}$ is the natural projection.
\smallskip

\par Unitarity: Modulo exponentially small errors and
microlocally, ${\cal T}$ is unitary $L^2(X)\to L
^2(\widetilde{X}; a_0e^{-2\Phi _0/h}L(dx))$, where $L(dx)$ is
chosen as indicated above, and $a_0(x;h)$ is a positive elliptic
analytic symbol of order $0$. \smallskip

\par Let $\Lambda \subset T^*\widetilde{X}$ be an IR-\mfld{} as
before, satisfying (\ref{16}) (or the equivalent condition (\ref{17}). Then
$\kappa _{\cal T}(\Lambda )=\Lambda _\Phi $, where $\Phi =\Phi
_\Lambda $, can be normalized by the requirement that $\Phi
=\Phi _0$ near the \bdy{} of $\widetilde{X}$. (Here is where we
have to choose $\lambda $ large enough, depending on $\Lambda $.
In the applications, for a given elliptic \op{}, $\Lambda $
and $T^*X$ will coincide outside a fixed compact \neigh{} of the
zero section, and the whole study will be carried out with a
fixed $\lambda $.)\smallskip

\par Let $\Omega \subset T^*X$ be the open \neigh{} of the
0-section, given by $\pi _x\kappa _{{\cal T}}\Omega
=\widetilde{X}$ and view also $\Omega $ as a subset of $\Lambda $
in the natural sense, assuming that $T^*X$ and $\Lambda $ coincide
in a \neigh{} of the closure of the complement of $\Omega $. If
$\chi \in C_0^\infty (\Omega )$, then the norm $\Vert u\Vert
_{H(\Lambda ,\langle \alpha _\xi \rangle ^m)}$ is equivalent to
the norm 
$$\Vert {\cal T}u\Vert _{L^2_\Phi }+\Vert (1-\chi )T_\Lambda
u\Vert _{L^2(\Lambda ;e^{-2H/h}\vert \langle \alpha _\xi \rangle
\vert ^{2m}d\alpha )}$$
\ufly{} \wrt{} $h$.

\end{document}